\documentclass[preprint,11pt]{elsarticle}

\usepackage{amssymb}
\usepackage{pgfplotstable}
\usepackage{pgfplots}
\usepackage{tikz}
\usepackage{amsmath}
\usepackage{amssymb}
\usepackage{epsfig}
\usepackage{hyperref}
\usepackage{url}
\usepackage{stackengine}
\usepackage{booktabs}
\usepackage{color}
\usepackage{tabularray}
\usepackage[colorinlistoftodos]{todonotes}

\newcommand{\old}[1]{}
\newcommand{\new}[1]{#1}

\newcommand{\eps}{\varepsilon}
\newcommand{\beq}{\begin{equation}}
\newcommand{\eeq}{\end{equation}}
\newcommand{\bea}{\begin{eqnarray}}
\newcommand{\eea}{\end{eqnarray}}

\newcommand{\bg}{\boldsymbol{g}}
\newcommand{\bG}{\boldsymbol{G}}

\newcommand{\bbf}{\boldsymbol{f}}
\newcommand{\br}{\boldsymbol{r}}

\newcommand{\half}{\mbox{$\frac{1}{2}$}}
\newcommand{\dif}{{\rm d}}

\newcommand{\teps}{{\eps_s}}
\newcommand{\tN}{\tilde{N}}
\newcommand{\tQ}{\tilde{Q}}

\journal{Comput. Methods Appl. Mech. Eng.}

\begin{document}

\begin{frontmatter}


\title{Shear-flexible geometrically exact beam element based on finite differences}

\author[label1]{Milan Jir\'{a}sek}

\author[label1,label3]{Martin  Hor\'{a}k}

\author[label2]{Emma La Malfa Ribolla\corref{cor1}}
\ead{emma.lamalfaribolla@unipa.it}

\author[label2]{Chiara Bonvissuto}

\cortext[cor1]{Corresponding author. Tel.: +39 09123896570}

\address[label1]{Department of Mechanics, Faculty of Civil Engineering,\\ Czech Technical University in Prague, Czech Republic}
\address[label2]{Department of Engineering, University of Palermo, Italy}
\address[label3]{Institute of Information Theory and Automation,\\ Czech Academy of Sciences, Czech Republic}

\begin{abstract}
The proposed two-dimensional geometrically exact beam element extends our previous work by including the effects of shear distortion, and also of distributed forces and moments acting along the beam. The general flexibility-based formulation exploits the kinematic equations combined with the inverted sectional equations and the integrated form of equilibrium equations. The resulting set of three first-order differential equations is discretized by finite differences and the boundary value problem is converted into an initial value problem using the shooting method. Due to the
special structure of the governing equations, the scheme remains explicit even though the first derivatives are approximated by central differences, leading to high accuracy. 
The main advantage of the adopted approach is that the error can be efficiently reduced 
by refining the computational grid used for finite differences at the element level while keeping the number of global degrees of freedom low. 
The efficiency is also increased by dealing directly with the global centerline coordinates and sectional inclination with respect to global axes as the primary unknowns at the element level, thereby avoiding transformations between local and global coordinates.  
\\
Two formulations of the sectional equations, \old{referred to as the Reissner and Ziegler models} \new{namely the widely used Reissner model and a less common version referred to as the Ziegler model}, are presented and compared. In particular, stability of an axially loaded beam/column is investigated and the connections to the Haringx and Engesser stability theories are discussed. 
Both approaches are tested in a series of numerical examples, which illustrate (i) high accuracy with quadratic convergence when the spatial discretization is refined, (ii) easy modeling of variable stiffness along the element (such as rigid joint offsets), (iii) efficient and accurate characterization of the buckling and post-buckling behavior.
\end{abstract}


\begin{keyword}
geometrically nonlinear beam \sep shear flexible beam \sep large rotations  \sep  shooting method \sep buckling \sep stability
\end{keyword}

\end{frontmatter}


\section{Introduction}\label{sec:intro}
In a wide range of applications and scientific areas, including mechanical and structural engineering, biomedical engineering, and soft robotics, slender beam- or rod-like elements are essential constituents. Examples are lattice-based metamaterials, which exploit structural topology to achieve certain favorable mechanical and physical properties, high-tensile industrial webbings, hard-magnetic soft rods and biomedical prostheses. The need to predict the mechanical behavior of these flexible structures has stimulated the development of increasingly efficient and robust computational tools based on the pioneering finite-strain beam models.\\
Long slender bodies can be modeled using various rod \cite{antman2005} and beam theories.
To introduce the terminology used throughout the paper, several pivotal studies on large-strain beam models are mentioned below. The origins date back to Euler's formulation \cite{euler} based on the assumptions of (i) perpendicularity between the beam centerline and cross sections and (ii) axial inextensibility. The Kirchhoff model \cite{kirchhoff1859} enhanced the formulation to account for axial strain and torsion, while still neglecting the shear strain. 
Under the assumptions of small displacements and rotations, the effects of shear were included
in the Timoshenko-Ehrenfest beam theory \cite{Tim21}.
Much later,
Reissner \cite{reissner1972} proposed the first full formulation that accounts for the shear distortion without kinematic restrictions. 
Based on Reissner's work, Simo \cite{simo1985} developed a formulation known as the geometrically exact beam theory, which is truly consistent in the sense that the kinematic quantities and relations are systematically derived from the 3D continuum theory, while the constitutive law is postulated. Various numerical implementations were based on finite element approaches, corotational methods \cite{Belytschko1977,Crisfield1990,Oran1976},  Absolute Nodal Coordinate formulations \cite{Shabana1998,Gerstmayr2013}, Lagrangian approaches \cite{Pagani2018}, finite differences\cite{Jung2011} and iso-geometric schemes\cite{Tasora2011,Marino2016}.

It should be emphasized that the Reissner–Simo beam model accounting for transverse shear deformation can lead to shear locking for displace\-ment-based finite element formulations\cite{Pai2013,Pai2014}. Substantial research work related to the shear-deformable beam model has been conducted and, even recently, new  theoretical\cite{Mohyeddin2014,Beiranvand2023} as well as computational\cite{Adnan2022,DiRe2024} developments have been proposed.
In a recent study, Kim et al. \cite{KimMehAtt20} have listed various possible choices of the internal forces in shear-flexible models. 


In our recent papers, we have presented a geometrically exact formulation for the planar behavior of beams that are initially straight \cite{HorMalJir23} or curved\cite{JirMalHor21}. Locally, the beam is treated by finite differences, but in the context of structural analysis it plays the same role as a traditional finite \old{beam} element. The advantage is that the accuracy of the numerical approximation can be increased by refinement of the finite difference grid while keeping the number of global degrees of freedom constant. A locking-free behavior is observed, which is closely related to the fact that no a priori chosen shape functions for the kinematic approximation are needed. 

The present study extends the beam element from\cite{JirMalHor21} by incorporating the effects of (i) shear distortion and (ii) distributed external forces and moments acting along the beam. In addition to the theoretical enhancement, the numerical scheme is also improved by a modified choice of the primary variables (centerline coordinates instead of displacements), which eliminates the need to convert variables from local to global reference systems.

In the following, a brief outline of the paper content is provided. In Section \ref{sec:beam model}, two models which differ in the choice of the basic deformation variables are \new{briefly} presented and the fundamental equations are \old{derived} \new{summarized (their detailed derivation can be found in \ref{appA})}. The first model (Section \ref{secreis}), referred to as the \emph{Reissner} model, is characterized by stretching and shear deformation variables work-conjugate to internal forces perpendicular and parallel to the deformed section, respectively. In this case, axial stiffness is activated only by the change of distance between the neighboring sections and not by parallel sliding of the sections. This choice finds its corresponding buckling model in the work of Haringx \cite{Haringx1942}, who studied the stability of helical springs. The second model (Section \ref{seczieg}), referred to as the \emph{Ziegler} model, is characterized by stretching and shear deformation variables work-conjugate to internal forces parallel and perpendicular to the deformed centerline, respectively. In this case, axial stretching occurs also during parallel sliding of sections. Such approach can be considered
as an extension of the original Engesser model for stability of a shear-deformable column \cite{Engesser1891} to the axially
compressible case, as presented by Ziegler \cite{Zie82}.

For both models, Section 3 shows how to treat the fundamental equations numerically in a general case by an efficient procedure that exploits the idea of the shooting method. In Section 4, numerical simulations of several examples are performed and compared with existing solutions from the literature \new{and with the analytical solution from \ref{appB}}. 
\old{
Section 5 is devoted to an analytical study in the framework of stability showing that the Reissner model predicts a certain finite critical load even in uniaxial tension, while the Ziegler model does not give any bifurcation in tension. On the other hand, in compression, the Reissner critical stress tends to infinity as the slenderness approaches zero, while the Ziegler model has a finite limit on the critical stress. These results are numerically confirmed and discussed in Section 6.} 
\new{Additional examples related to bifurcations from
the straight state and to the subsequent post-critical response 
are presented in Section~5, where the numerical results are compared
with analytical solutions derived in \ref{appC}.}
 Finally, Section \old{7} \new{6} summarizes the conclusions and discusses possible extensions. 


Throughout the paper, we will consistently use the following terminology:
\begin{itemize}
    \item{\bf Euler} refers to the 2D beam with an inextensible centerline and with sections perpendicular to the centerline, which deforms only by {\bf bending}.
    \item{\bf Kirchhoff} is used for the beam with sections perpendicular to the centerline, which deforms by {\bf axial stretching} and {\bf bending} (in 3D also by torsion).
    \item{\bf Reissner} and {\bf Ziegler} are two formulations that, in addition to {\bf axial stretching} and {\bf bending,} incorporate also the effect of {\bf shear} distortion. 
\end{itemize}
In particular, it is worth noting that we consider the Kirchhoff beam as axially extensible, even though this is not always recognized in the literature. 
For instance, Kirchhoff rods are in \cite{antman2005} 
considered as the 3D extension of Euler's elastica, with
torsion included but axial deformation still neglected.
However, 
Kirchhoff's original paper \cite{kirchhoff1859} explicitly includes a term that accounts for changes in the centerline length.

\section{Shear-flexible beam model} \label{sec:beam model}
\subsection{Basic assumptions and variables}

In the Timoshenko beam theory\cite{Tim21}, the condition of normality between the section and the deformed centerline (i.e., the Euler-Bernoulli hypothesis) is relaxed, and
 the transverse shear strain is no longer zero. 
 In one of his classic papers, Reissner\cite{reissner1972} presented a shear-flexible planar beam model that can be considered as a geometrically exact
 extension of the Timoshenko approach to large strains,
 since the relations between the deformation variables and
 the kinematic variables characterizing the displacements and rotations do not exploit any simplifications based on the
 assumed smallness of certain quantities.

\begin{figure}[h]
\hskip -20mm
    \begin{tabular}{cc}
    (a) & (b) \\
\includegraphics[scale=0.55]{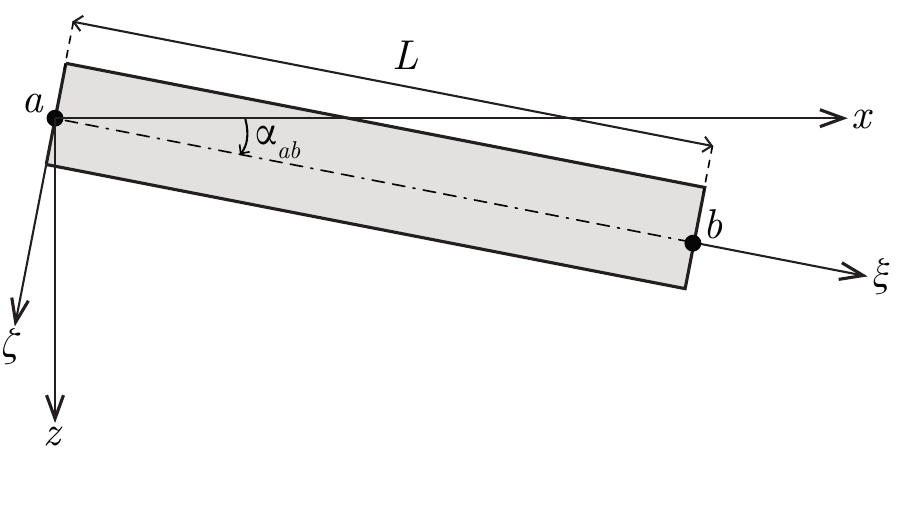}
&
 \includegraphics[scale=0.55]{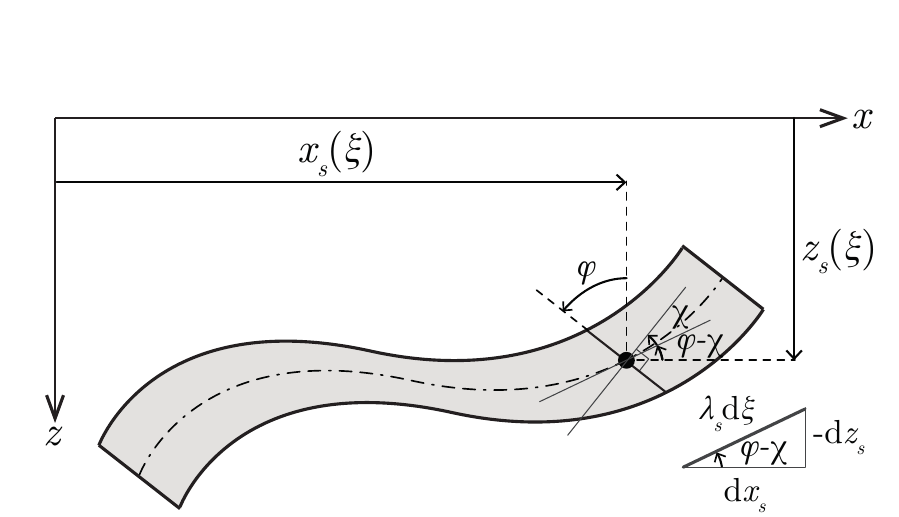}
    \end{tabular}
\caption{Kinematics of the geometrically exact beam model in the global $x$-$z$ plane: (a)~undeformed straight beam and (b) deformed beam.} 
\label{T1}
\end{figure}

Consider an initially straight \new{planar} beam of length $L$\old{ with a constant cross section}, which connects joints $a$ and $b$ \new{located in the global $x$-$z$
plane of the  global coordinate system, as indicated in Fig.~\ref{T1}a. A generic section initially located 
at distance $\xi$ from the left end is in the
deformed state
 described by (i) coordinates $x_s$ and $z_s$ of its centroid and (ii) inclination angle $\varphi$
 between the section (assumed to remain plane) and the 
 global $z$-axis.
The angle between the tangent to the deformed centerline 
and the normal
to the section, denoted here as $\chi$, represents
the shear distortion. 
Based on the geometry of the infinitesimal triangle shown in Fig.~\ref{T1}b, it is possible to express the shear angle
\beq \label{eqn_chinew}
\chi= \varphi - \arctan\frac{-\dif z_s}{\dif x_s}= \varphi + \arctan\frac{z_s'}{x_s'}
\eeq 
and also
the centerline stretch 
\beq\label{eqn_lamsnew}
\lambda_s = \frac{1}{\dif\xi}\sqrt{\dif x_s^2+\dif z_s^2}= \sqrt{x_s'^2+z_s'^2}
\eeq
Primes denote here derivatives with respect to $\chi$.
The bending part of deformation is characterized by the curvature,
usually defined as
\beq \label{ee5}
\kappa = \frac{\dif\varphi}{\dif\xi}=\varphi'
\eeq 
}

For a Kirchhoff beam\new{, the shear distortion is neglected}, the infinitesimal segment is 
decomposed into individual fibers parallel to the centerline,
and each fiber is considered to be under uniaxial tension. 
The stretch and thus also the strain in each fiber can be linked
to the centerline stretch and curvature by a simple relation.
Therefore, one can start from a given uniaxial hyperelastic stress-strain law, characterized by an expression for the strain energy density (per unit volume) as a function of strain, and then integrate over the section and construct the corresponding
strain energy density per unit length,
considered as a function of the axial strain and curvature. 
In this way, the sectional equations that link the deformation
variables to the internal forces can be deduced in a consistent way
from the selected constitutive model; see for instance \cite{JirMalHor21,ger2008,simo1986}.

Extension of this approach to shear-flexible beam models is in principle possible but is not straightforward, because 
it is no longer sufficient to work with a uniaxial stress-strain law.
Unless the application area is restricted to small-strain cases, 
the effects of the longitudinal fiber strain and the 
shear strain become coupled, and the transversal normal strain needs to
be eliminated based on a suitable assumption regarding the stress state\old{ (analogous to but more complicated than the assumption of vanishing
transversal normal stress for a fiber under uniaxial tension with
vanishing shear)}. The resulting mathematical problem is tractable
in closed form only in rare cases. Therefore, a simple alternative
route is frequently taken, and the strain energy density is directly
postulated on the level of an infinitesimal beam segment, as a function
of properly selected deformation variables.
\new{Two formulations of this kind are presented in detail in~\ref{appA}, and the resulting equations are summarized in the next two sections.}

\subsection{Reissner model}\label{secreis}

\old{A substantial part of this section has been moved to the appendix.}

\new{
The sectional equations in the form originally suggested by Reissner \cite{reissner1972} can be written as
\beq\label{reis-s} 
N=EA\,\eps, \hskip 10mm Q=GA_s\,\gamma, \hskip 10mm M=EI\,\kappa
\eeq 
They link the deformation variables
\beq\label{reis-k}
\varepsilon =\lambda_s \cos{\chi}-1, \hskip 10mm
\gamma = \lambda_s \sin{\chi}, \hskip 10mm \kappa=\varphi'
\eeq
to the work-conjugate internal forces,
with the normal force $N$ understood as the component
perpendicular to the section, and the shear force $Q$ as the component
parallel to the section. Symbol $M$ denotes the bending moment,
and $EA$, $GA_s$ and $EI$ are the sectional stiffnesses related
to axial stretching, shear distortion and bending.
It is useful to remark that the deformation variables $\eps$
and $\gamma$ defined in (\ref{reis-k}) can be directly linked to the kinematic
variables $x_s$, $z_s$ and $\varphi$, making use of relations
(\ref{eqn_chinew})--(\ref{eqn_lamsnew}). In particular,
in~\ref{appA1} it is shown that
\beq\label{reis-kk}
x_s' = (1+\eps)\cos\varphi+\gamma\sin\varphi, \hskip 10mm
z_s' = \gamma\cos\varphi-(1+\eps)\sin\varphi, \hskip 10mm \varphi'=\kappa
\eeq 
where the last equation is just copied from (\ref{reis-k}c).

A variational approach elaborated in detail in~\ref{appA1}
leads to equilibrium equations in the differential form (\ref{eq26y}-\ref{eq28y})
and to boundary conditions (\ref{ee26}-\ref{ee35}). After integration, the
equilibrium equations read
\bea \label{ee36new}
N &=& -\left(X_{ab}+P_x\right)\cos\varphi + \left(Z_{ab}+P_z\right)\sin\varphi \\
Q &=& -\left(X_{ab}+P_x\right)\sin\varphi - \left(Z_{ab}+P_z\right)\cos\varphi \label{ee37new}
\\
\label{ee41new}
M &=& -M_{ab} + X_{ab}(z_s-z_a)-Z_{ab}(x_s-x_a) + M_p
\eea 
where
\bea 
P_x(\xi) &=& \int_0^\xi p_x(s)\,\mathrm{d}s
\\
P_z(\xi) &=& \int_0^\xi p_z(s)\,\mathrm{d}s
\\
\label{e41new}
M_p(\xi) &=& -\int_0^\xi m(s)\,{\rm d}s + \int_0^\xi P_x(s)\, z_s'(s)\, \textrm{d}s  -  \int_0^\xi P_z(s)x_s'(s)\, \textrm{d}s
\eea
are functions that can be deduced from the prescribed distributed
load intensities $p_x$, $p_z$ and $m$. In the above,
$X_{ab}$, $Z_{ab}$ and $M_{ab}$
are constants that represent the global components of the force and the moment between the left end of the beam and the joint to which it is attached (see Fig.~\ref{f:defbeam}), and $x_a$ and $z_a$ are the current coordinates of the centroid of the left end section, i.e., $x_a=x_s(0)$ and $z_a=z_s(0)$. 

Equations (\ref{reis-s}) and (\ref{reis-kk})--(\ref{ee41new}) can be combined into a set of three first-order
differential equations that form the basis of the 
numerical scheme to be developed next. Such combined equations 
 (\ref{ee51})--(\ref{ee53}) presented in~\ref{appA1} are rather lengthy, and so the algorithm actually implements 
individual equations (\ref{reis-s}) and (\ref{reis-kk})--(\ref{ee41new})
as intermediate steps, to keep the procedure transparent. 
}

\subsection{Ziegler model}\label{seczieg}
\old{A substantial part of this section has been moved to the appendix.}

\new{
The path to the formulation that we refer to as the Ziegler model is described in~\ref{appA2}. The selected deformation 
$\teps = \lambda_s-1$, $\chi$ and $\kappa$
are linked to the internal forces by sectional equations
\beq 
\tN = EA\,\teps ,\hskip 10mm
\tQ =  GA_s\,\chi ,\hskip 10mm
M  = EI\,\kappa \label{ee19xx}
\eeq
The difference compared to the Reissner model is that the 
axial stretching mode is characterized by axial strain $\teps$ and the shear mode by the shear angle $\chi$,
which leads to work-conjugate internal forces $\tN$ and $\tQ$
that represent the components in the direction parallel and perpendicular to the centerline (and not to the section). 
Moreover, $\tQ$ is actually the force perpendicular to the
centerline, $Q^*$, multiplied by the axial stretch $\lambda_s$;
see \ref{appA2} for a detailed explanation.
The variational approach then leads to differential equilibrium
equations which, after integration, have the form
\bea \label{ee36znew}
\tN &=& -\left(X_{ab}+P_x\right)\cos(\varphi-\chi) + \left(Z_{ab}+P_z\right)\sin(\varphi-\chi) \\
Q^* &=& -\left(X_{ab}+P_x\right)\sin(\varphi-\chi) - \left(Z_{ab}+P_z\right)\cos(\varphi-\chi) \label{ee37znew}
\\
\label{ee41znew}
M &=& -M_{ab} + X_{ab}(z_s-z_a)-Z_{ab}(x_s-x_a) + M_p
\eea
where $Q^*=\tQ/\lambda_s$ is the true shear force. 
Relations (\ref{reis-kk}) valid for the Reissner model are now replaced by
\beq \label{zieg-kk}
x_s' = (1+\teps) \cos(\varphi-\chi) ,\hskip 5mm
z_s'  =-(1+\teps) \sin(\varphi-\chi),\hskip 5mm \varphi'=\kappa
\eeq 
and the resulting set of three first-order differential equations
is given by (\ref{ee110})--(\ref{ee111})  and (\ref{ee137}) in~\ref{appA2}.
}


\section{Numerical procedures}

\subsection{Basic concept}

\subsubsection{Evaluation of generalized end forces}

\begin{figure}[h]
    \centering
    \includegraphics[width=0.6\linewidth]{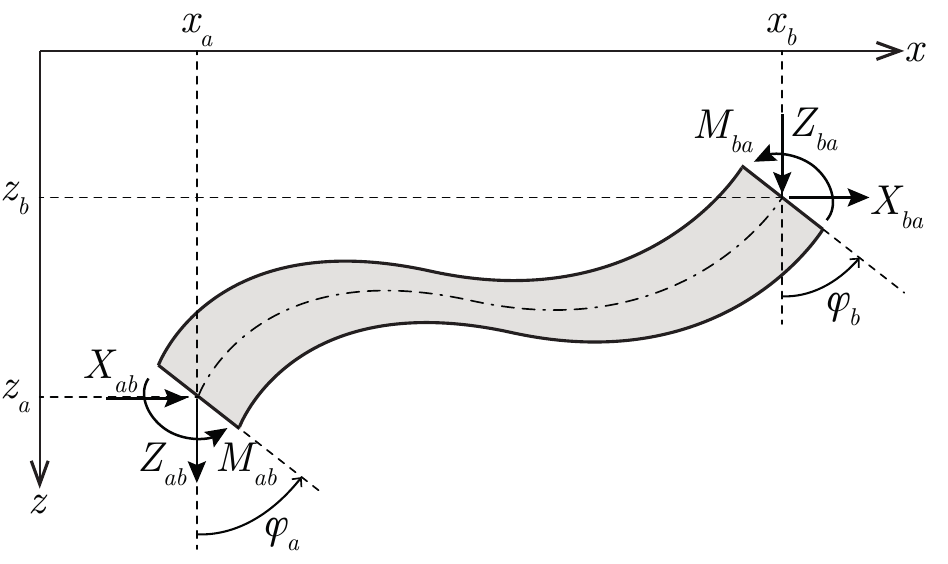}
    \caption{Deformed beam, generalized end coordinates and generalized end forces}
    \label{f:defbeam}
\end{figure}

We have shown that, for the assumptions adopted, the description of
a shear-flexible beam in the $x$-$z$ plane can be reduced to a set
of three first-order differential equations for three primary
unknown functions, $x_s$, $z_s$ and $\varphi$, which characterize the deformed centerline and
the sectional inclination. The numerical solution of these equations
will be based on the finite difference scheme. For the solution
to be unique, three first-order differential equations require three additional conditions, but for
our problem we have a total of six boundary conditions, because
the end \old{displacements and rotations} \new{coordinates and inclinations} are supposed to be given  \new{(see Fig.~\ref{f:defbeam})}.
At the same time, the governing differential equations contain
not only the primary unknown functions and given constants
(such as the sectional stiffnesses) but also the left-end forces, $X_{ab}$ and $Z_{ab}$,
and the left-end moment, $M_{ab}$. Their values need to be determined so that the six
kinematic boundary conditions can be satisfied simultaneously.
This will be achieved by using a shooting method. Three kinematic
boundary conditions at the left end will be used as initial conditions
for integration of three first-order differential equations, 
and the left-end forces
and moment will be considered as additional unknown variables
(not functions). Their values will first be estimated (for instance,
set to zero, or taken from the previous step or iteration at the structural
level), and then iteratively adjusted until the numerically computed
right-end coordinates and inclination get sufficiently close to their
prescribed values.

To formalize the outlined procedure, let us introduce the column
matrices of generalized left-end forces, $\bbf_a=(X_{ab},Z_{ab},M_{ab})^T$,
generalized right-end forces, $\bbf_b=(X_{ba},Z_{ba},M_{ba})^T$,
generalized left-end coordinates, $\br_a=(x_a,z_a,\varphi_a)^T$ and
generalized right-end coordinates, $\br_b=(x_b,z_b,\varphi_b)^T$.
Their meaning is apparent from Fig.~\ref{f:defbeam}. 
Notice that the force components and coordinates are taken with respect to the global coordinate system, and that the inclinations $\varphi_a$ and 
$\varphi_b$ are understood as the angles by which the
respective section deviates from the global $z$-axis counterclockwise.
They differ from the rotations of the structural joints to which the
end sections are attached by the initial deviation of the sections
from the vertical direction. More specifically, if joints $a$ and $b$
rotate by $\Phi_a$ and $\Phi_b$, and the beam centerline in the undeformed
state deviates by angle $\alpha_{ab}$ clockwise from the global $x$-axis \new{(as indicated in Fig.~\ref{T1}a)},
then we have
\bea 
\varphi_a &=& \Phi_a-\alpha_{ab} \\
\varphi_b &=& \Phi_b-\alpha_{ab}
\eea 
Since, for a given beam, $\alpha_{ab}$ is a fixed constant,
the increments or variations of $\varphi$ and $\Phi$ are the same,
only the total values differ by a constant. Incorporation of the 
initial beam position directly into the beam-level algorithm avoids
the need for transformation of forces, displacements and coordinates
between the global and local (beam-related) coordinate systems.

A numerical algorithm, to be described later, can integrate the
governing differential equations along the beam and compute the
generalized coordinates of the right end, $\br_b$, based on the
specified generalized coordinates of the left end, $\br_a$,
and generalized left-end forces, $\bbf_a$. This defines a mapping
$\bg$ such that
\beq \label{ee140}
\br_b = \bg(\bbf_a,\br_a)
\eeq 
If the generalized coordinates at both ends are given, condition 
(\ref{ee140}) can be interpreted as a set of three nonlinear algebraic
equations with three unknowns---the generalized forces assembled in $\bbf_a$.
This set of equations is solved by the Newton-Raphson method,
starting from an initial guess, $\bbf_a^{(0)}$. In a generic
iteration number $k$, equations  (\ref{ee140}) are replaced by their
linearized version
\beq \label{ee141}
\bg(\bbf_a^{(k-1)},\br_a)+ \bG(\bbf_a^{(k-1)},\br_a)\,\Delta\bbf_a^{(k)} = \br_b 
\eeq 
where
\beq 
\bG(\bbf_a,\br_a) = \frac{\partial\bg(\bbf_a,\br_a)}{\partial\bbf_a}
\eeq 
is the Jacobi matrix of partial derivatives of function $\bg$ with
respect to the first set of arguments, $\bbf_a$. Components of this
matrix are computed simultaneously with the evaluation of $\bg$
using consistent linearization of the algorithm around the current solution.
The iterative correction 
\beq \label{ee143}
\Delta\bbf_a^{(k)} = \bG^{-1}(\bbf_a^{(k-1)},\br_a)\left(\br_b-\bg(\bbf_a^{(k-1)},\br_a)\right)
\eeq 
is then added to the current approximation of $\bbf_a$, and an improved
approximation 
\beq \label{ee144}
\bbf_a^{(k)} = \bbf_a^{(k-1)} + \Delta\bbf_a^{(k)}
\eeq 
is evaluated. The iterative loop runs until the norm of the
residual $\br_b-\bg(\bbf_a^{(k)},\br_a)$
becomes negligibly small (a weighted norm is used, because two components of the residual have the dimension of length while the last one is dimensionless). When convergence is achieved, the most
recent approximation $\bbf_a^{(k)}$ represents the generalized
forces acting on the left end of the beam, which can be directly 
incorporated into the global equilibrium conditions of joint $a$ on the
structural level. The generalized right-end forces, $\bbf_b^{(k)}$, are
provided among the output variables of the algorithm that integrates the fundamental equations along the beam. Alternatively, they can be evaluated from the equilibrium equations
written for the whole beam, which lead to
\bea 
X_{ba} &=&  - X_{ab} - P_z(L) \\
Z_{ba} &=&  - Z_{ab} - P_z(L) \\
M_{ba} &=&  -M_{ab}+X_{ab}(z_b-z_a)-Z_{ab}(x_b-x_a)+M_p(L)
\label{e109}
\eea 
The evaluation of $M_p(L)$ using formula (\ref{e41new}) with $\xi=L$ would require
numerical integration, and so it is better to provide the value as
one of the output variables of the numerical algorithm, which needs
to compute $M_p$ along the way anyway. In the absence of distributed
loading, $M_p$ vanishes. 

\subsubsection{Evaluation of element stiffness matrix}
\label{sec:3.1.2}

For the solution of the joint equilibrium equations at the structural level, the beam element must be able to provide not only the end forces, but also the element tangent stiffness matrix,
which links infinitesimal increments of generalized end displacements to infinitesimal increments of generalized end forces.  
Suppose that $\br_a$ and $\br_b$ are the prescribed generalized
end coordinates and $\bbf_a$ are the corresponding generalized left-end
forces, calculated from equation (\ref{ee140}). Infinitesimal increments of these variables are linked by the linearized form
of (\ref{ee140}), which can be written as
\beq \label{ee140x}
\dif\br_b = \bG(\bbf_a,\br_a)\,\dif\bbf_a + \bG_r(\bbf_a,\br_a)\,\dif\br_a
\eeq 
where $\bG$ is the already defined Jacobi matrix of partial derivatives
$\partial\bg/\partial\bbf_a$, and 
\beq \label{eq:bgr}
\bG_r = \frac{\partial\bg}{\partial\br_a}
\eeq 
is another matrix that contains partial derivatives of mapping $\bg$ with respect to the
second set of arguments---the generalized left-end displacements.

For simplicity, we will omit the arguments $\bbf_a$ and $\br_a$
in the subsequent derivations. 
The (infinitesimal) increment of left-end forces that would be caused
by increments $\dif\br_a$ and $\dif\br_b$ of the generalized
end coordinates (which coincide with the increments of displacements
and rotations of the joints to which the beam element is attached)
are easily expressed from (\ref{ee140x}) as
\beq \label{ee140y}
 \dif\bbf_a= \bG^{-1} \left(\dif\br_b- \bG_r\,\dif\br_a\right) = -\bG^{-1} \bG_r\,\dif\br_a+\bG^{-1} \dif\br_b
\eeq 
Consequently, matrix $-\bG^{-1} \bG_r$ represents the upper left
$3\times 3$ block of the tangent element stiffness matrix
that links the column $(\dif\bbf_a^T,\dif\bbf_b^T)^T$ (increments of generalized 
end forces) to the column $(\dif\br_a^T,\dif\br_b^T)^T$ (increments of generalized
end coordinates or displacements), and matrix $\bG^{-1}$ is the
upper right $3\times 3$ block of the tangent element stiffness matrix. In this way, we obtain the first three rows of the stiffness matrix, and then, based
on equilibrium equations written for the whole beam, we can
construct the remaining three rows. 

The Jacobi matrix $\bG$ is computed by the end-force evaluation
algorithm, because it is needed for the iterative process
described by (\ref{ee143})--(\ref{ee144}). For evaluation 
of the tangent stiffness, we also need matrix $\bG_r$, but not
all of its entries have to be computed numerically. 
If the left-end coordinates $x_a$ and $z_a$ are changed by
$\dif x_a$ and $\dif z_a$ and everything else is kept fixed,
then the resulting right-end coordinates $x_b$ and $z_b$ will
also be changed by $\dif x_a$ and $\dif z_a$, and the resulting
right-end inclination $\varphi_b$ will remain unchanged. 
Physically, one can think of a   rigid body
translation applied to the deformed beam, which does not
disturb the validity of the basic equations. Consequently,
the first two columns of matrix $\bG_r$ are $(1,0,0)^T$
and $(0,1,0)^T$, and the first two columns of matrix
$\bG^{-1}\bG_r$ will be just copies of the first two columns
of $\bG^{-1}$. 

Let us now imagine what happens if the left-end inclination
$\varphi_a$ is changed by $\dif\varphi_a$. One might 
expect a similar effect as if the whole beam is rotated by
$\dif\varphi_a$, keeping the same deformed shape. However,
this is true only if the forces acting on the beam rotate as well.
For a beam subjected to distributed loading, this cannot be
easily emulated, and the third column of $\bG_r$ needs to
be evaluated numerically. Only if the beam is not loaded by
forces distributed along its length, the third column can be
constructed analytically.
If we rotate the whole deformed beam and simultaneously also
the end forces by $\dif\varphi_a$ counterclockwise, then  the generalized coordinates
at the right end change by $\dif x_b=(z_b-z_a)\,\dif\varphi_a$,
$\dif z_b=-(x_b-x_a)\,\dif\varphi_a$ and $\dif\varphi_b=\dif\varphi_a$, as an effect of a rigid
rotation of the beam, during which the deformed shape remains the same. The same change can be thought of as resulting from the effect of rotation of the left end section
at constant (non-rotating) left-end forces, combined
with the effect of a modification of the global components
of left-end forces by changes that would
arise during rotation of these forces. 
Since the rotation $\dif\varphi_a$ is infinitesimal, 
counterclockwise rotation of the horizontal force contributes
by $-X_{ab}\dif\varphi_a$ to the vertical component
while counterclockwise rotation of the vertical force contributes
by $Z_{ab}\dif\varphi_a$ to the horizontal component.
Using the formalism of equation (\ref{ee140x}), this decomposition can be described by
\beq \label{ee151}
\left(\begin{array}{c} z_b-z_a \\ x_a-x_b \\ 1 \end{array}\right)\dif\varphi_a =
\bG
\left(\begin{array}{r} Z_{ab} \\ -X_{ab} \\ 0  \end{array}\right)\dif\varphi_a
+
\bG_r
\left(\begin{array}{c} 0 \\ 0 \\ 1 \end{array}\right)\dif\varphi_a
\eeq 
Cancelling $\dif\varphi_a$ and regrouping the terms, we find out that
 the third column of $\bG_r$ could be evaluated by summing the 
column $(z_b-z_a,x_a-x_b,1)^T$, the first column
of $\bG$ multiplied by $-Z_{ab}$, and the second column of $\bG$
multiplied by $X_{ab}$.
We do not even have to evaluate this expression, because what is actually needed is the third column of $-\bG^{-1}\bG_r$,
and from (\ref{ee151}) multiplied from the left by $\bG^{-1}$  we can deduce that
\beq 
-\bG^{-1}\bG_r
\left(\begin{array}{c} 0 \\ 0 \\ 1 \end{array}\right) = 
\bG^{-1}\left(\begin{array}{c} z_a-z_b \\ x_b-x_a \\ -1 \end{array}\right) +
\left(\begin{array}{r} Z_{ab} \\ -X_{ab} \\ 0  \end{array}\right)
\eeq 
The left-hand side is the third column of $-\bG^{-1}\bG_r$,
i.e., of the upper left $3\times 3$ block of the element
stiffness matrix, and the right-hand side provides a rule
for its evaluation: Take the column $(Z_{ab},-X_{ab},0)^T$
and add a linear combination of the columns of $\bG^{-1}$
with coefficients $z_a-z_b$, $x_b-x_a$ and $-1$.
All this is true in the absence of distributed external forces. 

If nonzero external forces acting on the beam are prescribed, the first two columns
of $-\bG^{-1}\bG_r$ still correspond to the first two
columns of $\bG^{-1}$ with reverted signs, but the third
column of $\bG_r$
must be computed numerically by evaluating the effect
of $\dif\varphi_a$ on $\dif x_b$, $\dif z_b$ and $\dif\varphi_b$, and then multiplication from the left by $-\bG^{-1}$ leads to the third
column of $-\bG^{-1}\bG_r$.

As already explained, matrices $-\bG^{-1}\bG_r$ and $\bG^{-1}$
together fill the upper half of the $6\times 6$ stiffness
matrix, i.e., the first three rows. Based on force equilibrium of the whole beam, it is easy to see that the
fourth row is equal to minus the first row, and the fifth
row is equal to minus the second row. 
The sixth row can be deduced from the moment equilibrium condition, which is somewhat more tricky. 
Differentiation of (\ref{ee41}) leads to
\bea
\nonumber
\dif M_{ba} &=& -\dif M_{ab} + \dif X_{ab}(z_b-z_a)+ X_{ab}(\dif z_b-\dif z_a)-\dif Z_{ab}(x_b-x_a)+\\
&-& Z_{ab}(\dif x_b-\dif x_a) + \dif M_p(L)
\eea
In the absence of distributed external loads, $M_p(L)$
vanishes and its change $\dif M_p(L)$ as well, 
and the formula can be translated into a prescription
for the sixth row of the element stiffness matrix:
Initialize by $(Z_{ab},-X_{ab},0,-Z_{ab},X_{ab},0)$,
add the first row multiplied by $z_b-z_a$ and
the second row multiplied by $x_a-x_b$, and subtract the third row.

In the presence of external loads, the contribution 
of $\dif M_p(L)$ must be added. From the algorithm in Section~\ref{sec:3.2.3} it will be clear that
the sensitivities of $M_p$ at
individual grid points (to be introduced in the next section) to infinitesimal changes $\dif X_{ab}$,
$\dif Z_{ab}$, $\dif M_{ab}$ and $\dif\varphi_a$ must be computed anyway, because they affect the evaluation of $\dif\varphi$ at individual grid points. 
When the solution process reaches the right end
of the beam, these sensitivities evaluated for
the last grid point represent the derivatives
of $M_p(L)$ with respect to the left-end generalized
forces and to the left-end inclination
(notice that changes of the left-end coordinates have no effect on
$M_p$). The derivatives with respect to the left-end generalized
forces are then used as three coefficients that multiply 
the first, second and third row of the stiffness matrix in order to obtain 
the contribution of distributed loads 
to the last row of the stiffness matrix.
The fourth coefficient, i.e., the derivative 
of $M_p(L)$ with respect to $\varphi_a$,
is simply
added to the stiffness entry in the third column
of the last row, i.e., to $K_{63}$.

\subsection{Finite difference scheme}

Recall that the governing differential equations to be solved are described by (\ref{ee51})--(\ref{ee53}) for the Reissner model and by
(\ref{ee110})--(\ref{ee111}) and (\ref{ee137})
for the Ziegler model. 
The spatial derivatives of unknown functions will be replaced by suitable finite
difference approximations.
Interestingly, for both models, the equations that
describe the spatial derivatives of $x_s$ and $z_s$
have right-hand sides that can be evaluated exclusively
from $\varphi$, and the equation  
that
describes the spatial derivative of $\varphi$
has a right-hand side that can be evaluated exclusively
from $x_s$ and $z_s$. Consequently, we can use
central differences for the approximation of the 
first derivatives while keeping the algorithm explicit.
Let us divide the interval $[0,L]$ into $N$ equally
sized segments of length $\Delta\xi=L/N$, and
define the coordinates of grid points
\beq 
\xi_i = i\,\Delta\xi, \hskip 10mm i=0,1,2,\ldots N
\eeq 
and midpoints
\beq 
\xi_{i-1/2} = (i-1/2)\,\Delta\xi, \hskip 10mm i=1,2,\ldots N
\eeq 
If the governing equations are presented in the symbolic form 
\bea 
x_s' &=& f_x(\varphi) \\
z_s' &=& f_z(\varphi) \\
\varphi' &=& f_\varphi(x_s,z_s) 
\eea 
then their finite difference approximations are 
\bea 
\frac{x_i-x_{i-1}}{\Delta\xi} &=& f_x(\varphi_{i-1/2}), \hskip 10mm i=1,2,\ldots N \\
\frac{z_i-z_{i-1}}{\Delta\xi} &=& f_z(\varphi_{i-1/2}), \hskip 10mm i=1,2,\ldots N \\
\frac{\varphi_{i+1/2}-\varphi_{i-1/2}}{\Delta\xi} &=& f_\varphi(x_i,z_i), \hskip 10mm i=1,2,\ldots N-1 
\eea 
where $x_i$ and $z_i$ are numerical values of functions $x_s$ and $z_s$ at the grid points, and $\varphi_{i-1/2}$
are numerical values of function $\varphi$ at the midpoints.
To initiate the algorithm, $x_0$ and $z_0$ are set to the values $x_a$ and $z_a$ equal to the current coordinates
of the left end centroid, and $\varphi_{1/2}$ is evaluated
from the given $\varphi_0=\varphi_a$ by adopting 
a forward Euler formula
\beq 
\frac{\varphi_{1/2}-\varphi_{0}}{\Delta\xi/2} =f_\varphi(x_0,z_0)
\eeq 
for the first half-step. In a similar fashion,
the value of $\varphi_N$ at the end of the interval is 
calculated from a backward Euler formula
\beq 
\frac{\varphi_{N}-\varphi_{N-1/2}}{\Delta\xi/2} = f_\varphi(x_N,z_N)
\eeq 
adopted for the last half-step. Even this formula
is explicit, because the values of $x_N$ and $z_N$
are already known and $\varphi_N$ is not needed on the
right-hand side.

Functions $f_x$, $f_z$ and $f_\varphi$ could be directly replaced
by the complicated expressions on the right-hand sides
of (\ref{ee51})--(\ref{ee53}), or of
(\ref{ee110})--(\ref{ee111}) and (\ref{ee137}).
However, for a clean numerical implementation, it is preferable to evaluate individual terms sequentially, which 
allows for easier physical interpretation of the intermediate
results. For example, for the Reissner model
we get back to the individual
equations from which (\ref{ee51})--(\ref{ee53}) was composed, noting that the differential character
of the equations comes from (\ref{reis-kk}),
while the remaining relations (\ref{reis-s}) and (\ref{ee36new})--(\ref{ee41new}) are algebraic.

\subsection{Reissner model}

\subsubsection{Evaluation of right-end generalized coordinates}
\label{sec:3.2.1}

The basic algorithm needed for the implementation is the numerical
evaluation of function $\bg(\bbf_a,\br_a)$, 
which takes as input the left-end generalized forces
$\bbf_a=(X_{ab},Z_{ab},M_{ab})^T$ and left-end generalized
coordinates $\br_a=(x_a,z_a,\varphi_a)^T$ and provides
as output the right-end generalized
coordinates $\br_b=(x_b,z_b,\varphi_b)^T$.

If the beam is loaded by distributed forces, it is useful to precompute and store the partial load resultants
\old{(the equations have been compacted by using a generic subscript $r$)}
\bea \label{ee162}
P_{\new{r},1/2} &=& \frac{p_r(0)+p_r(\xi_{1/2})}{4} \Delta\xi
\\
P_{\new{r},i+1/2} &=& P_{r,i-1/2} + \frac{p_r(\xi_{i-1/2})+p_r(\xi_{i+1/2})}{2} \Delta\xi, \hskip 5mm i=1,2,\ldots N-1 \nonumber \\
\\
P_{\new{r},N} &=& P_{r,N-1/2} + \frac{p_r(L)+p_r(\xi_{N-1/2})}{4} \Delta\xi
\label{ee167}
\eea 
\new{where $r\in\{x,z\}$.}
If the distributed forces are fixed, this needs be done only once, even before the incremental-iterative
solution of the beam structure is started.
If the distributed forces are increased proportionally to a load factor
which is varied in each step, the pre-computed partial resultants 
represent reference values and in the algorithm below need to be 
multiplied by the load factor. It is worth noting that one could easily incorporate concentrated forces acting on the beam as well. Their global components would be simply added to $P_{x,i+1/2}$ and $P_{z,i+1/2}$ when the corresponding coordinate 
$\xi_{i+1/2}$ becomes for the first time greater than the coordinate of the section at which the concentrated force is applied. 

The algorithm for evaluation of $\bg(\bbf_a,\br_a)$ can be organized as follows:

\begin{enumerate}
    \item Set initial values $x_0=x_a$, $z_0=z_a$, $\varphi_0=\varphi_a\equiv\Phi_a-\alpha_{ab}$, $M_{p0}=0$, and
    $M_0=-M_{ab}$.
    \item For $i=1,2,\ldots N$ evaluate
\bea \label{e234mj3}
\varphi_{i-1/2} &=& \varphi_{i-1} + \frac{M_{i-1}}{EI}\frac{\Delta\xi}{2} \\
c&=& \cos\varphi_{i-1/2} \\
s&=& \sin\varphi_{i-1/2} \\
\label{N135}
N_{i-1/2} &=&  -c\left(X_{ab}+P_{x,i-1/2}\right) +s\left(Z_{ab}+P_{z,i-1/2}\right) \\ 
\label{Q136}
Q_{i-1/2} &=& - s\left(X_{ab}+P_{x,i-1/2}\right)  -c\left(Z_{ab}+P_{z,i-1/2}\right) \\
\eps_{i-1/2} &=& \frac{N_{i-1/2}}{EA} \\
\gamma_{i-1/2} &=& \frac{Q_{i-1/2}}{GA_s} \label{ee177} \\
\label{dx139}
\Delta x &=& \left(c(1+\eps_{i-1/2})+s\gamma_{i-1/2}\right)\Delta\xi \\
\label{dz140}
\Delta z &=& \left(c\gamma_{i-1/2}-s(1+\eps_{i-1/2})\right)\Delta\xi \\
x_i &=& x_{i-1}+\Delta x\\
z_i &=& z_{i-1}+\Delta z \\
M_{pi} &=& M_{p,i-1}-m(\xi_{i-1/2})\,\Delta\xi+P_{x,i-1/2}\,\Delta z - P_{z,i-1/2}\,\Delta x \nonumber \\
\\
M_i &=& -M_{ab}+X_{ab}\left(z_i-z_a\right)-Z_{ab}\left(x_i-x_a\right) +M_{pi}\label{e77}  \\ 
\label{e240mj3}
\varphi_i &=& \varphi_{i-1/2} + \frac{M_i}{EI}\frac{\Delta\xi}{2}
\eea 
\item The resulting coordinates and inclination at the right end are $x_b=x_N$, $z_b=z_N$ and $\varphi_b=\varphi_N$. It is also useful to provide $M_{pN}$ as part of the output.
\end{enumerate}
Notice that $c$ and $s$ are values of the cosine and sine
of angle $\varphi_{i-1/2}$, which change in each cycle
but for brevity are not denoted by subscript $i-1/2$.
These auxiliary variables are also used in the loop that
computes the Jacobi matrix, to be presented in the next section.

It is also worth noting that the present algorithm could 
be adjusted to accommodate the Kirchhoff beam model as a limit
case with an infinite shear stiffness, resulting into a
zero shear measure $\gamma$. Evaluation of 
$Q_{i-1/2}$ in (\ref{Q136}) and of
$\gamma_{i-1/2}$
in (\ref{ee177}) would be skipped, and the terms with  $\gamma_{i-1/2}$
would be deleted from (\ref{dx139})--(\ref{dz140}). Instead of introducing
branches in the algorithm or writing a special code for
the Kirchhoff case, it is even simpler to implement
the operation in (\ref{ee177}) as multiplication of the shear force
by the sectional shear compliance $C_s\equiv 1/GA_s$, and to consider
this compliance (instead of the stiffness $GA_s$) as one of the
input parameters describing the beam. 
The results for a Kirchhoff beam can then be obtained
without touching the algorithm, simply by using
$C_s=0$ on input.
One can even go further and implement the fraction
in (\ref{N135}) as the multiplication of the normal force by the
sectional axial compliance $C_a\equiv 1/EA$,
and the fractions in (\ref{e234mj3}) and (\ref{e240mj3})
as the multiplication of the bending moment by the
sectional flexural compliance $C_b\equiv 1/EI$.
Then, by using $C_s=0$ and $C_a=0$, we obtain the Euler elastica, and by using also $C_b=0$, we obtain a perfectly
rigid beam.
In general, the sectional properties can vary along the beam, 
and by setting $C_s=0$, $C_a=0$ and $C_b=0$ in certain
numerical segments, we can easily combine flexible
parts of the beam with rigid
segments, e.g., near the joints, which might be useful
to obtain a more realistic frame model for 3D-printed
beam assemblies, as demonstrated in \cite{DiRe2024}.
Only minor changes in the code are needed, consisting
in the added control of the values of sectional compliances
at individual grid points. The actual sequence of operations in the 
algorithm does not need to be changed.

\subsubsection{Evaluation of Jacobi matrix}
\label{sec:3.2.3}

The Jacobi matrix is needed for the iterative evaluation of the generalized end forces on the element level as well as for the construction of the element stiffness matrix. 
The entries of the Jacobi matrix $\bG(\bbf_a,\br_a)$ are evaluated numerically using the linearized form of the computational scheme.
Simultaneously, the derivatives $\partial\bg(\bbf_a,\br_a)/\partial\varphi_a$ are computed, since they are needed
to construct the element stiffness matrix, as described
in Section~\ref{sec:3.1.2}.

Suppose that the input values $X_{ab}$, $Z_{ab}$ and $M_{ab}$ are changed by infinitesimal increments
$\dif X_{ab}$, $\dif Z_{ab}$ and $\dif M_{ab}$.
Linearization of equations (\ref{e234mj3})--(\ref{e240mj3}) around the currently considered solution leads to the following algorithm:
\begin{enumerate}
    \item 
    Set one of the values of $\dif X_{ab}$,
    $\dif Z_{ab}$, $\dif M_{ab}$ or $\dif\varphi_a$ to 1 and the others to 0, depending on which column of the extended Jacobi matrix
    is to be computed (the fourth column now corresponds to  $\partial\bg/\partial\varphi_a$).
    Set initial values $\dif x_0=0$, $\dif z_0=0$, $\dif M_{p0}=0$, and
    $\dif M_0=-\dif M_{ab}$. 
    \item For $i=1,2,\ldots N$ evaluate
\bea 
\dif\varphi_{i-1/2} &=& \dif\varphi_{i-1} + \frac{\dif M_{i-1}}{EI}\frac{\Delta\xi}{2} \\
\dif N_{i-1/2} &=&  -c\,\dif X_{ab} +s\,\dif Z_{ab} -Q_{i-1/2}\,\dif\varphi_{i-1/2}\\ 
\dif Q_{i-1/2} &=& - s\,\dif X_{ab}  -c\,\dif Z_{ab} +N_{i-1/2}\,\dif\varphi_{i-1/2}\\
\dif\eps_{i-1/2} &=& \frac{\dif N_{i-1/2}}{EA} \\
\label{ee189}
\dif\gamma_{i-1/2} &=& \frac{\dif Q_{i-1/2}}{GA_s} \\
\dif\Delta x &=& \left(c\,\dif\eps_{i-1/2}+s\,\dif\gamma_{i-1/2}-\left(s\,(1+\eps_{i-1/2})-c\,\gamma_{i-1/2}\right)\dif\varphi_{i-1/2}\right) \Delta\xi\\
\dif\Delta z &=& \left(c\,\dif\gamma_{i-1/2}-s\,\dif\eps_{i-1/2}- \left(s\,\gamma_{i-1/2}+c\,(1+\eps_{i-1/2})\right)\dif\varphi_{i-1/2}\right)\Delta\xi\\
\dif x_i &=& \dif x_{i-1}+\dif\Delta x\\
\dif z_i &=& \dif z_{i-1}+\dif\Delta z \\
\dif M_{pi} &=& \dif M_{p,i-1}+P_{x,i-1/2}\,\dif\Delta z - P_{z,i-1/2}\,\dif\Delta x
\\ 
\dif M_i &=& -\dif M_{ab}+\dif X_{ab}(z_i-z_a)+X_{ab}\,\dif z_i-\dif Z_{ab}(x_i-x_a) -Z_{ab}\,\dif x_i+\dif M_{pi}\nonumber\\ \label{e77e}  \\ 
\dif\varphi_i &=& \dif\varphi_{i-1/2} + \frac{\dif M_i}{EI}\frac{\Delta\xi}{2}
\eea 
    \end{enumerate}
For $\dif X_{ab}=1$ and $\dif Z_{ab}=\dif M_{ab}=\dif\varphi_a=0$, the final values of
$\dif x_N$, $\dif z_N$ and $\dif\varphi_N$ obtained by running this 
algorithm represent entries of the first column of the Jacobi matrix $\bG$.
In a similar fashion, an evaluation working with $\dif X_{ab}=0$, $\dif Z_{ab}=1$ and $\dif M_{ab}=\dif\varphi_a=0$ provides the second column, and an evaluation working with $\dif X_{ab}=\dif Z_{ab}=\dif\varphi_a=0$ and $\dif M_{ab}=1$ provides the third column. In practice, all three columns of 
the Jacobi matrix are computed in a single loop, along with the value
of function $\bg$. Also, in the presence of distributed loads,
the algorithm is simultaneously executed with  $\dif X_{ab}=\dif Z_{ab}=\dif M_{ab}=0$ and $\dif\varphi_a=1$, leading to the third
column of matrix $\bG_r$, which contains the derivatives  $\partial\bg/\partial\varphi_a$. Everything is done within a single loop, and so the values of internal forces $N_{i-1/2}$ and $Q_{i-1/2}$
and deformation variables $\eps_{i-1/2}$ and $\gamma_{i-1/2}$
are those calculated by the algorithm presented in Section~\ref{sec:3.2.1}, and $c$ and $s$ have the meaning of the cosine and sine of $\varphi_{i-1/2}$.
Since they are reused, the corresponding formulas are not repeated here. All variables with identifiers starting by the symbol ``d''
are in fact arrays with 3 or 4 values, depending on whether we 
compute the columns of $\bG$ only, or also the third column of $\bG_r$.

\subsection{Ziegler model}\label{sec:3.4}

\subsubsection{Iterative scheme for shear angle evaluation}

For evaluation of internal forces at midpoint $\xi_{i-1/2}$,
we need the midpoint shear angle $\chi_{i-1/2}$, which is defined implicitly as the solution of nonlinear equation
(\ref{eq92}), in which $P_1$ is replaced by $X_{ab}+P_x$, $P_2$ by $Z_{ab}+P_z$, $\varphi$ by $\varphi_{i-1/2}$ and $\chi$ by $\chi_{i-1/2}$. Since the midpoint rotation $\varphi_{i-1/2}$ can be evaluated independently of $\chi_{i-1/2}$, and 
the forces $X_{ab}$, $Z_{ab}$, $P_x$ and $P_z$ are given
(in the context of this algorithm), this equation can be 
presented in the form
\beq 
F\left(\chi_{i-1/2}\right) = 0
\eeq 
where function
\bea 
F(\chi) = GA_s\chi - \left(1+\frac{f_N(\chi)}{EA}\right)f_Q(\chi)
\eea 
is defined using auxiliary functions
\bea 
f_N(\chi) &=& -(X_{ab}+P_x)\cos(\varphi_{i-1/2}-\chi) + (Z_{ab}+P_z) \sin(\varphi_{i-1/2}-\chi) \nonumber\\ \\
f_Q(\chi) &=& -(X_{ab}+P_x)\sin(\varphi_{i-1/2}-\chi) - (Z_{ab}+P_z) \cos(\varphi_{i-1/2}-\chi)\nonumber\\
\eea 
It is easy to see that $\dif f_N/\dif\chi=f_Q$ and
$\dif f_Q/\dif\chi=-f_N$, and so
\beq
F'(\chi)\equiv\frac{\dif F(\chi)}{\dif \chi} = GA_s + f_N(\chi)\left(1+\frac{f_N(\chi)}{EA}\right)- \frac{f_Q^2(\chi)}{EA}
\eeq 
The iterative process typically starts from the first guess $\chi^{(0)}=\chi_{i-3/2}=$ shear angle at the midpoint of the previous spatial segment (for $i=1$, the initial guess is set to zero). Improved approximations are calculated
using the recursive formula
\beq 
\chi^{(j)} = \chi^{(j-1)} - \frac{F(\chi^{(j-1)})}{F'(\chi^{(j-1)})}, \hskip 5mm j=1,2,\ldots
\eeq 
and it is terminated when the absolute value of $F(\chi^{(j)})$
becomes sufficiently small. The converged value $\chi^{(j)}$
is then used as $\chi_{i-1/2}$, and we can proceed to the
evaluation of other quantities at the midstep. In fact,
the values of $f_N(\chi^{(j)})$ and $f_Q(\chi^{(j)})$
calculated during the evaluation of $F(\chi^{(j)})$
are directly the internal forces at midstep, $\tilde N_{i-1/2}$
and $Q^*_{i-1/2}$, and the already evaluated expression
$1+f_N(\chi^{(j)})/EA$ is the centerline stretch, $\lambda_{i-1/2}$, from which the corresponding
increments of primary variables $x_s$ and $z_s$ can be evaluated.

\subsubsection{Evaluation of right-end displacements and rotation}
\label{sec:3.3.2}

If the beam is loaded by distributed forces, it is useful to precompute and store the partial load resultants using the
same formulae (\ref{ee162})--(\ref{ee167}) as for the Reissner model. This is done only once. Afterwards, the following algorithm is invoked.

\begin{enumerate}
    \item Set 
     initial values $x_0=x_a$, $z_0=z_a$, $\varphi_0=\varphi_a\equiv\Phi_a-\alpha_{ab}$, $M_{p0}=0$ and
    $M_0=-M_{ab}$. Set the tolerance $\epsilon_{tol}$ of the internal iteration loop.
\item For $i=1,2,\ldots N$ evaluate
\bea \label{ee215}
\varphi_{i-1/2} &=& \varphi_{i-1} + \frac{M_{i-1}}{EI}\frac{\Delta\xi}{2} \\ 
\chi_{i-1/2} &=& \chi_{i-3/2} \hskip 5mm (\mbox{or } \chi_{i-1/2}=0 \mbox{ if } i=1 )\\
    \mbox{repeat:}\hskip 10mm\nonumber\\
    c &=& \cos(\varphi_{i-1/2}-\chi_{i-1/2}) \label{ee217}\\
     s &=& \sin(\varphi_{i-1/2}-\chi_{i-1/2}) \\
     \tilde N_{i-1/2} &=& -c\left(X_{ab}+P_{x,i-1/2}\right)+s\left(Z_{ab}+P_{z,i-1/2}\right) \\
    Q^*_{i-1/2} &=& -s\left(X_{ab}+P_{x,i-1/2}\right)-c\left(Z_{ab}+P_{z,i-1/2}\right) \\
    \lambda_{i-1/2} &=& 1+\frac{\tilde N_{i-1/2}}{EA} \\
    F &=& GA_s\chi_{i-1/2} - \lambda_{i-1/2} Q^*_{i-1/2} \\  \nonumber
    &&\mbox{if } \; \vert F \vert < \epsilon_{tol} \mbox{ then go to (\ref{ee225})} \\
    F'&=& GA_s + \lambda_{i-1/2} \tilde N_{i-1/2} - \frac{Q_{i-1/2}^{*2}}{EA} \\
    \chi_{i-1/2} &=& \chi_{i-1/2}-F/F' \label{ee224}\\ \nonumber
      \mbox{until}\hskip 10mm && \mbox{the maximum number of internal}\\  && \mbox{iterations is exceeded (report failure)}\nonumber
\eea
\bea\label{ee225}
\Delta x &=& c\lambda_{i-1/2} \,\Delta\xi
\\
\Delta z &=& -s\lambda_{i-1/2}\,\Delta\xi
\\
x_i &=& x_{i-1}+\Delta x\\
z_i &=& z_{i-1}+\Delta z \\
M_{pi} &=& M_{p,i-1}-m(\xi_{i-1/2})\,\Delta\xi+P_{x,i-1/2}\,\Delta z - P_{z,i-1/2}\,\Delta x
\\
M_i &=& -M_{ab}+X_{ab}\left(z_i-z_a\right)-Z_{ab}\left(x_i-x_a\right) +M_{pi}  \\ 
\varphi_i &=& \varphi_{i-1/2} + \frac{M_i}{EI}\frac{\Delta\xi}{2}
\label{ee231}
\eea 
\item The resulting coordinates and rotation at the right end are $x_b=x_N$, $z_b=z_N$ and $\varphi_b=\varphi_N$. It is also useful to provide $M_{pN}$ as part of the output.
\end{enumerate}

\subsubsection{Evaluation of Jacobi matrix}

Suppose that the input values $X_{ab}$, $Z_{ab}$ and $M_{ab}$ are changed by infinitesimal increments
$\dif X_{ab}$, $\dif Z_{ab}$ and $\dif M_{ab}$.
Linearization of the algorithm described in Section~\ref{sec:3.3.2} around the currently considered solution is straightforward, but special
attention needs to be paid to the implicitly defined shear angle
$\chi_{i-1/2}$, which is evaluated
in the internal
iterative loop (\ref{ee217})--(\ref{ee224}). The converged solution
is affected by changes in $X_{ab}$ and $Z_{ab}$, and also by the change in $\varphi_{i-1/2}$. Consistent linearization of condition
$F=0$ leads to
\beq \label{ee232}
GA_s\,\dif\chi_{i-1/2} - Q_{i-1/2}^*\,\dif\lambda_{i-1/2}-\lambda_{i-1/2}^*\,\dif Q_{i-1/2}^* = 0
\eeq
in which
\bea \label{ee233}
\dif\lambda_{i-1/2} &=& \frac{\dif \tilde N_{i-1/2}}{EA} \\
\dif \tilde N_{i-1/2}
\label{ee234}
&=& -c\,\dif X_{ab}+s\,\dif Z_{ab} -Q_{i-1/2}^*\left(\dif\varphi_{i-1/2}-\dif\chi_{i-1/2}\right)
\\
\dif Q_{i-1/2}^* 
\label{ee235}
&=& -s\,\dif X_{ab}-c\,\dif Z_{ab} +\tilde N_{i-1/2}\left(\dif\varphi_{i-1/2}-\dif\chi_{i-1/2}\right)
\eea 
Substituting (\ref{ee233})--(\ref{ee235}) into (\ref{ee232}) and solving for $\dif\chi_{i-1/2}$, we get
\beq 
\dif\chi_{i-1/2} = \scalebox{0.8}{$\dfrac{\left(-c\,Q^*_{i-1/2}/EA-s\,\lambda_{i-1/2}\right)\dif X_a+\left(s\,Q^*_{i-1/2}/EA-c\,\lambda_{i-1/2}\right)\dif Z_a+\left(\lambda_{i-1/2}\tilde N_{i-1/2}-Q^{*2}_{i-1/2}/EA\right)\dif\varphi_a}{GA_s+\lambda_{i-1/2}\tilde N_{i-1/2}-Q^{*2}_{i-1/2}/EA}$}
\eeq
In combination with the linearized form of the explicit expressions
(\ref{ee215}) and (\ref{ee225})--(\ref{ee231}), the following algorithm is obtained:
\begin{enumerate}
    \item 
    Set one of the values of $\dif X_{ab}$,
    $\dif Z_{ab}$, $\dif M_{ab}$ or $\dif\varphi_a$ to 1 and the others to 0, depending on which column of the Jacobi matrix
    is to be computed.
    Set initial values $\dif x_0=0$, $\dif z_0=0$, $\dif M_{p0}=0$, and
    $\dif M_0=-\dif M_{ab}$.
\item For $i=1,2,\ldots N$ evaluate
\bea 
\dif\varphi_{i-1/2} &=& \dif\varphi_{i-1} + \frac{\dif M_{i-1}}{EI}\frac{\Delta\xi}{2} \\ 
    c &=& \cos\left(\varphi_{i-1/2}-\chi_{i-1/2}\right) \\
     s &=& \sin\left(\varphi_{i-1/2}-\chi_{i-1/2}\right) \\
     \nonumber
 \dif\chi_{i-1/2} &=& \scalebox{0.8}{$\dfrac{\displaystyle\left(-\frac{cQ^*_{i-1/2}}{EA}-s\lambda_{i-1/2}\right)\dif X_a+\left(\frac{sQ^*_{i-1/2}}{EA}-c\lambda_{i-1/2}\right)\dif Z_a+\left(\lambda_{i-1/2}\tilde N_{i-1/2}-\frac{Q^{*2}_{i-1/2}}{EA}\right)\,\dif\varphi_a}{\displaystyle GA_s+\lambda_{i-1/2}\tilde N_{i-1/2}-\frac{Q^{*2}_{i-1/2}}{EA}}$}\\
 \\
\dif\tilde N_{i-1/2} &=&  -c\,\dif X_{ab}+s\,\dif Z_{ab} -Q_{i-1/2}^*\left(\dif\varphi_{i-1/2}-\dif\chi_{i-1/2}\right)
\\
\dif\lambda_{i-1/2}&=&\frac{\dif\tilde N_{i-1/2}}{EA}
\\
\dif\Delta x &=& \left(c\,\dif\lambda_{i-1/2}-s\lambda_{i-1/2}\left(\dif\varphi_{i-1/2}-\dif\chi_{i-1/2}\right)\right)\Delta\xi
\\
\dif\Delta z &=& \left(-s\,\dif\lambda_{i-1/2}-c\lambda_{i-1/2}\left(\dif\varphi_{i-1/2}-\dif\chi_{i-1/2}\right)\right)\Delta\xi
\\
\dif x_i &=& \dif x_{i-1}+\dif\Delta x\\
\dif z_i &=& \dif z_{i-1}+\dif \Delta z \\
\dif M_{pi} &=& \dif M_{p,i-1}+P_{x,i-1/2}\,\dif\Delta z - P_{z,i-1/2}\,\dif\Delta x
\\
\nonumber
\dif M_i &=& -\dif M_{ab}+\dif X_{ab}\left(z_i-z_a\right)+X_{ab}\,\dif z_i-\dif Z_{ab}\left(x_i-x_a\right) -Z_{ab}\,\dif x_i+\dif M_{pi}  \\  \\
\dif\varphi_i &=& \dif\varphi_{i-1/2} + \frac{\dif M_i}{EI}\frac{\Delta\xi}{2}
\eea 
\item For $\dif X_{ab}=1$ and $\dif Z_{ab}=\dif M_{ab}=\dif \varphi_0=0$, the final values of
$\dif x_N$, $\dif z_N$ and $\dif\varphi_N$ obtained by running this 
algorithm represent the entries of the first column of the Jacobi matrix $\bG$.
In a similar fashion, an evaluation working with $\dif X_{ab}=0$, $\dif Z_{ab}=1$ and $\dif M_{ab}=\dif\varphi_a=0$ provides the second column, and an evaluation working with $\dif X_{ab}=\dif Z_{ab}=\dif\varphi_a=0$ and $\dif M_{ab}=1$ provides the third column. In practice, all three columns of 
the Jacobi matrix are computed in a single loop, along with the value
of function $\bg$. Also, in the presence of distributed loads,
the algorithm is simultaneously executed with  $\dif X_{ab}=\dif Z_{ab}=\dif M_{ab}=0$ and $\dif\varphi_a=1$, leading to the third
column of matrix $\bG_r$.
\end{enumerate}

\section{Examples---bending with shear}
\label{sec:examples}

The accuracy and efficiency of the proposed numerical schemes will be tested using several 
benchmark examples.
Unless stated otherwise, the examples consider 
an isotropic linear elastic material with
Young's modulus 
denoted as $E$ and Poisson's ratio set to $\nu=0.25$, 
which leads
to the shear modulus $G=E/(2(1+\nu))=0.4\,E$.
The cross section is taken as rectangular, 
characterized by width $b$ and depth $h$,
with the shear reduction factor
set to the usual value of 5/6 (even though Cowper's approach\cite{Cowper1966} would lead to a slightly different value), which provides the sectional shear stiffness
$GA_s=5GA/6=EA/3$ and dimensionless parameter
$\Gamma\equiv GA_s/EA=1/3$. 
Furthermore, the bending stiffness can be expressed
as $EI=Ebh^3/12=EAh^2/12$. 
A characteristic length (typically the beam span) is denoted as $L$.
To eliminate the dependence of the results on the specific
choice of input parameters,
we normalize forces by $EI/L^2$ and lengths by $L$,
which means that moments are normalized by $EI/L$. 
Since the Poisson ratio $\nu=0.25$ is considered as fixed, there remains
a single dimensionless input parameter $h/L$, i.e., the
depth-to-span ratio. The resulting deflections and other displacements
are normalized by $L$.

\subsection{Simply supported beam, concentrated force}

\begin{figure}[t]
\centering
 \includegraphics[scale=1.0]{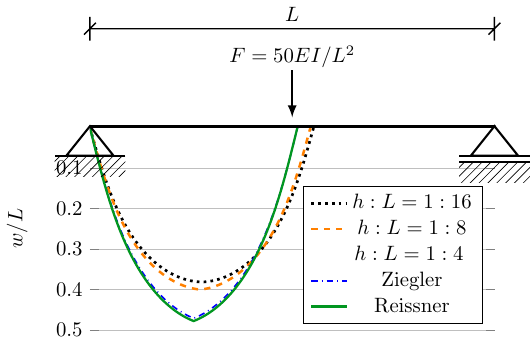}
    \caption{Simply supported beam loaded by a concentrated force $F$ at midspan: undeformed beam with supports and loading, and deformed centerline shapes at load level $F=50\,EI/L^2$ evaluated for beams of various depth-to-span ratios using the Reissner model
    and, for $h:L=1:4$, also the Ziegler model}
    \label{ex4.1_defshapes}
\end{figure}

Numerical examples start from a simply supported beam of span $L$, loaded
by a vertical concentrated force $F$ at midspan, as indicated in Fig.~\ref{ex4.1_defshapes}.
The prescribed load is converted into the dimensionless
factor $FL^2/(EI)$, and the resulting midspan deflection is reported in its
dimensionless form $w/L$. 
To explore the effect of the beam shape (slenderness), let us analyze
three cases with $h/L$=1/4, 1/8, and 1/16.

\new{\subsubsection{Deformed shapes}}

\begin{figure}[ht]
    \centering
    \begin{tabular}{cc}
    (a) & (b) 
    \\
    \includegraphics[scale=0.8]{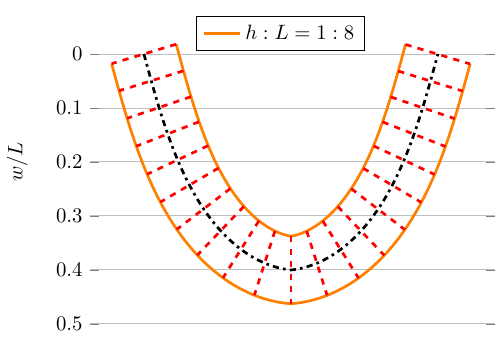}
    &  \includegraphics[scale=0.8]{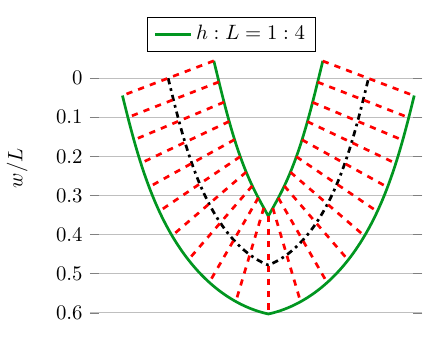}
    \end{tabular}
\caption{Simply supported beam loaded by a concentrated force $F$ at midspan: the deformed centerline and \old{rotation of}  selected sections \old{obtained with the Reissner model} at load level $F=50\,EI/L^2$ evaluated for beams of depth-to-span ratios 1:8 and 1:4 using the Reissner model. } 
\label{ex4_1_beams}
\end{figure}

The problem is symmetric with respect to the vertical axis (except for the nonsymmetry in the horizontal support, which could be easily eliminated), and so it could be solved using only one half of the beam discretized by one
single element. Equivalent results are obtained by 
discretizing the whole beam into two equally sized elements,
each of them of length $L/2$, divided into a certain number of integration segments.
In this case, the concentrated force is applied at the
node located at the midspan. Fig.~\ref{ex4.1_defshapes} shows the deformed 
centerlines for the three considered geometries
(depth-to-span ratios) at the same normalized load level, obtained using the Reissner model (dotted, dashed and dash-dotted curves).
Since the load is normalized by the bending stiffness,
the largest deflection is found for the stocky beam
with $h/L=1/4$, for which the relative importance
of the shear deformation is higher than for more slender beams. A kink in the middle of the deflection curve is
clearly visible in this case. \new{As illustrated by Fig.~\ref{ex4_1_beams},} the sectional inclination
$\varphi$ \old{(not seen in the figure)} varies smoothly
and remains differentiable, because its derivative 
is the curvature, proportional to the bending moment,
which is continuous. However, the slope
of the centerline differs from $\varphi$ by the shear
angle $\chi$, which is related to the shear force and
has a jump at the section under the applied concentrated force. Such a jump is present in all the cases, but it is
less conspicuous for slender beams.
For comparison, Fig.~\ref{ex4.1_defshapes} also contains the 
deformed centerline obtained with the Ziegler model for the
stocky beam with $h/L=1/4$ (solid curve). The difference between the Reissner
and Ziegler models is seen to be very small in this example.
For slender beams, it is almost negligible, and only the 
Reissner curves are shown.

\new{\subsubsection{Load-displacement diagrams}}

\begin{figure}[t]
\hskip -20mm
    \begin{tabular}{cc}
    (a) & (b) \\    
\includegraphics[scale=1]{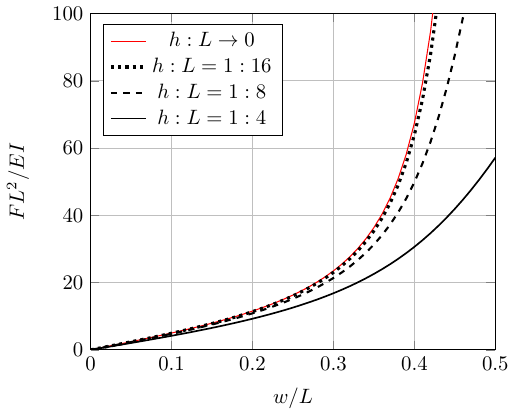}
&
\includegraphics{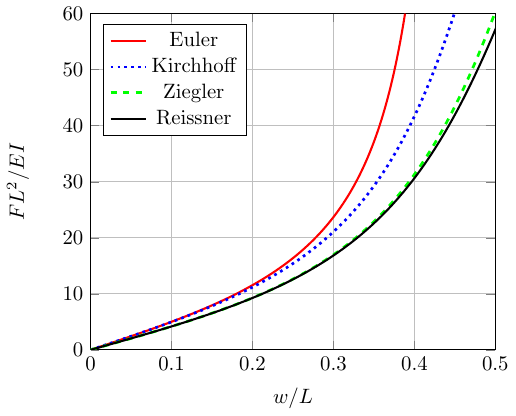}
    \end{tabular}
\caption{Simply supported beam loaded by a concentrated force at midspan: normalized load-displacement curves (a) for the Reissner model and various depth-to-span ratios,
(b) for depth-to-span ratio 1:4 and \old{both} \new{various} models} 
\label{ex4_1}
\end{figure}

Load-displacement diagrams covering a wide range up to midspan deflections close to one half of the span
are plotted in Fig.~\ref{ex4_1}.
In part (a), the Reissner model is used and the results are
presented for various depth-to-span ratios.
The curve obtained in the limit of $h/L\to 0$ is also included.
This limit corresponds to the Euler model with zero axial
deformation and zero shear deformation, enforced by using $1/EA=0$ and $1/GA_s=0$, as described at the end of Section~\ref{sec:3.2.1}). 
For a slender beam with $h/L=1/16$, the full Reissner
model with axial and shear flexibility included
gives only a slight increase of the deflection
(dashed curve). As expected, the increase is much more pronounced for stocky beams, especially for the depth-to-span ratio 1:4 (black curve).


 The load-displacement
curves  plotted in
Fig.~\ref{ex4_1}a correspond to the Reissner model. 
It turns out that for slender beams, such as those with $h/L=1/16$ and 1/8, the Ziegler model
gives load-displacement curves that are graphically indistinguishable from those for the Reissner formulation. 
Visible differences are found for the stocky beam
with $h/L=1/4$, as seen in Fig.~\ref{ex4_1}b.
In addition to the full Reissner and Ziegler formulations,
the curves are plotted for the Kirchhoff model with shear deformation
neglected but axial deformation still taken into account,
and also for the Euler model, which neglects both effects. 
The comparison shows that, in the highly nonlinear range,
both effects (axial and shear deformation) significantly 
increase the resulting deflections. The Reissner formulation leads to
slightly larger deflections than the Ziegler model.

\new{\subsubsection{Initial stiffness}}

According to the linear Timoshenko theory, the \new{normalized} mid-span deflection would be
\old{ 
\beq 
w = \frac{FL^3}{48EI} + \frac{FL}{4GA_s}
\eeq 
which can be written in terms of the dimensionless quantities as
\beq 
\frac{w}{L} = \frac{FL^2}{EI}\left(  \frac{1}{48} + \frac{EI}{4GA_sL^2}\right) = \frac{FL^2}{EI} \frac{1}{48}\left(  1 + \frac{3h^2}{L^2}\right)
\eeq 
}
\new{
\beq \label{eq200}
\frac{w}{L} = \frac{FL^2}{48EI}+\frac{F}{4GA_s} = \frac{FL^2}{EI} \frac{1}{48}\left(  1 + \frac{3h^2}{L^2}\right)
\eeq 
}
This means that the dimensionless structural stiffness (initial slope of the \old{dimensionless}
load-displacement diagram \new{that links the normalized deflection to the normalized force $FL^2/EI$}) \old{would be} \new{is} 48 for the Kirchhoff beam with no shear
deformation, while for the Timoshenko beam it is reduced by the factor
$1+3h^2/L^2$, leading to the dimensionless stiffnesses of 40.421,
45.851 and 47.444 for $h/L=1/4$, $1/8$ and $1/16$, respectively.

The initial slope of the load-displacement curves from Fig.~\ref{ex4_1} corresponds to the stiffness of the
Timoshenko beam. \old{, which is available in an analytical form,
based on (\ref{eq200}).} This gives us an opportunity to
evaluate the numerical error induced by discretization
in space. The initial structural stiffness extracted from the
tangent stiffness matrix calculated for the initial
undeformed state is reported in Table~\ref{tab:eigenval2x}.
The results have been evaluated for the Reissner model,
but the initial stiffness computed
for the Ziegler model would be exactly the same.
As the number of integration segments is increased,
the numerically obtained value quickly converges
to the theoretical limit, no matter whether the beam is slender ($h/L=1/16$) or stocky ($h/L=1/4$). Already for 8 integration segments per element, the error
is below 1~\%, and it decreases in proportion to the 
square of the segment size. 

\new{\subsubsection{Discretization error in nonlinear range}}

A comparable (even slightly lower) error level is observed 
in the nonlinear range, as documented in Table~\ref{tab:convergence2},
which shows the deflection obtained for a fixed, rather high load.
The reported deflection values refer to the state 
after convergence of global equilibrium iterations (with a sufficiently strict convergence criterion), but they still depend
on the spatial resolution. For both models, Reissner and Ziegler, the relative accuracy
is virtually the same when the same spatial resolution is used, even though the reference values
obtained in the limit are slightly different.
The \old{``almost exact''} \new{``converged''} 
reference values have been computed by using 10,000 integration
segments, which is overkill (six-digit accuracy would be obtained already with 500 segments). Let us emphasize that the refinement of the
grid inside the element does not affect the global number
of degrees of freedom, which remains equal to 6 in the
present example (can be reduced to 3 if symmetry is exploited, and even to 2 if the midspan deflection is prescribed instead of the midspan force).

 \begin{table}[t]
     \centering
  \caption{Simply supported beam loaded by a concentrated force at midspan and discretized by two elements: convergence of the initial structural stiffness\\}
     \label{tab:eigenval2x}
     \begin{tabular}{rllll}
     \hline
number of & \multicolumn{2}{c}{$h/L=1/4$} & \multicolumn{2}{c}{$h/L=1/16$} \\
segments& $FL^3/(w\,EI)$ & error $[\%]$ & $FL^3/(w\,EI)$ & error $[\%]$\\
          \hline
2	 &  36.5714  &  9.52 & 42.2268& 11.0 \\
4	 & 39.3846  &  2.56 & 46.0225&  3.00 \\
8	 &  40.1569  & 0.65  &  47.0805  & 0.77 \\
16	 &  40.3547  &  0.16 &  47.3526  & 0.19  \\
32	 &  40.4044  &  0.041 &  47.4211  & 0.048  \\
64	 &  40.4169  &  0.010 & 47.4383   &  0.012 \\
128	 &  40.4200  &   0.0026 &  47.4426  & 0.0030  \\
\old{$\infty$} \new{analytical} &  40.4210  &  & 47.4440  & \\
\hline
     \end{tabular}
 \end{table}

  \begin{table}[bth]
  \caption{\textcolor{black}{
Simply supported beam loaded by a concentrated force at midspan
and discretized by two elements: convergence of the mid-span displacement obtained for a given load $F=50\,EI/L^2$ as the number of integration segments per element is increased\\}}
     \label{tab:convergence2}
     \hskip -20mm
     \begin{tabular}{rllllllll}
     \hline
& \multicolumn{4}{c}{Reissner model} &   \multicolumn{4}{c}{Ziegler model}  \\
number of & \multicolumn{2}{c}{$h/L=1/4$} & \multicolumn{2}{c}{$h/L=1/16$} & \multicolumn{2}{c}{$h/L=1/4$} & \multicolumn{2}{c}{$h/L=1/16$}\\
segments & $w/L$ & error $[\%]$ & $w/L$ & error $[\%]$  & $w/L$ & error $[\%]$ & $w/L$ & error $[\%]$\\
          \hline
2	 & 0.506722   &   5.99 & 0.420842 & 10.28 & 0.499423 & 5.98 & 0.420817 & 10.28\\
4	 & 0.486911 &  1.84 &  0.392607 &  2.88 & 0.479664 & 1.79 & 0.392582  & 2.88\\
8	 & 0.480365  & 0.48    & 0.384369  & 0.72 & 0.473427 & 0.46 & 0.384345 & 0.72\\
16	 &  0.478647 & 0.12  & 0.382297   &  0.18 & 0.471796 & 0.12 & 0.382273 & 0.18\\
32	 &  0.478214  & 0.030  &  0.381779  &  0.045 & 0.471385 & 0.029 & 0.381755 & 0.045\\
64	 &  0.478105  & 0.007  &  0.381650  &  0.011 & 0.471282 & 0.007 & 0.381626 & 0.011\\
128	 &  0.478078  &  0.0019 &   0.381617 &  0.0028 & 0.471256 & 0.0019 & 0.381594 & 0.0028\\
\new{converged} & 0.478069  &   &  0.381607 & & 0.471247 &  & 0.381583\\
          \hline
     \end{tabular}
 \end{table}

It is worth noting that the quadratic nature of the 
convergence rate (as the spatial grid is refined) is independent of the slenderness and,
for a given number of integration
points, 
the error for a slender beam is only slightly larger
than for a stocky beam. This indicates that the present formulation does not suffer \old{by} \new{from} shear locking. 
\new{For an extremely slender beam with $h/L=1/64$,
the relative error in the mid-span deflection at load level $F=50\,EI/L^2$ computed using 16 integration segments per element would be $0.19~\%$, which is almost the same as for $h/L=1/16$
(the converged value of $w/L$ is 0.375826 while the approximate value for 16 segments is 0.376523).}

In fact, the algorithm is, after a minor change, applicable to the
Kirchhoff formulation that neglects shear, without losing
the nice convergence properties. The minimal modification
is to replace formula (\ref{ee177}) by $\gamma_{i-1/2}=0$
and formula (\ref{ee189}) by $\dif\gamma_{i-1/2}=0$;
see also the discussion at the end of Section~\ref{sec:3.2.1}.

\new{\subsubsection{Comparison with results from the literature}}

\new{
Owing to symmetry, the simply supported beam considered here is equivalent to a cantilever of span $L/2$, loaded at its free end
by a vertical force $F/2$. For this elementary problem, certain
reference solutions are available. Analytical solutions
were derived for the Reissner model with neglected axial compliance
by Batista \cite{Batista2016}, and the results presented in  Table~3
of \cite{Batista2016} apply to our beam with $F=20$, $L=2$ and $EI=10$. The results were reported for selected values of the
shear stiffness $GA_s$  and for an extremely high axial stiffness (we use $EA=10^8$). Table~\ref{tab:batista} compares the
analytical values from \cite{Batista2016} with our numerical
results  for
values of $GA_s=5\times 10^{20}$ and 500.

\begin{table}[h]
    \centering
    \begin{tabular}{lllll}
    \hline
      force, $F$   & 20  & 20 & 200  & 200 \\
      shear stiffness, $GA_s$ & $5\times 10^{20}$ & 500& $5\times 10^{20}$ & 500\\
         \hline
    exact, Batista \cite{Batista2016}     & 0.301720774  & 0.317813874 \\
    Lyritsakis et al.~\cite{Lyritsakis2021} & 0.3017206 & 0.3178137 & 0.8106113 & 0.8536137 \\
       present, converged & 0.3017208 & 0.3178139 & 0.8106099 & 0.8539632 \\
    present, 16 segments & 0.3022736 & 0.3183590 & 0.8123628 & 0.8554802 \\
    \hskip 12mm relative error & $0.18~\%$  & $0.17~\%$ & $0.22~\%$ & $0.18~\%$
    \\ \hline
    \end{tabular}
    \caption{\new{Deflections reported in the literature and computed with the present approach using the Reissner model }}
    \label{tab:batista}
\end{table}

Numerical solutions of the same problem were reported in Table~1 of Lyritsakis et al.~\cite{Lyritsakis2021}, who used a hybrid beam element based on nonlinear programming. In addition to the moderate load level for which the deflection is about 30~\% of the cantilever length (corresponding to 15~\% of the span of a simply supported beam), they considered also 10 times higher load,
in our notation $F=200$. The comparison in Table~\ref{tab:batista}
indicates that, for load
$F=20$, our converged values of maximum deflection (obtained with a sufficiently high number of integration segments) perfectly agree with the
analytical solutions from \cite{Batista2016} as well as with the numerical values from \cite{Lyritsakis2021}. 
For the high load  $F=200$, there is a small difference between
our converged results and the values reported in \cite{Lyritsakis2021}, which were obtained using an integration
scheme with 6 quadrature points and probably were not fully converged. Our converged result for $F=200$ and extremely high $GA_s$
agrees in the first 6 digits with the exact solution 0.8106090 reported
in \cite{Lyritsakis2021}.
For simulations with 16 integration segments, the relative
errors reported in the last row of Table~\ref{tab:batista}
are around $0.2~\%$, even for the highly deformed state 
induced by force $F=200$.
}

\subsection{Clamped beam, uniformly distributed load}\label{sec:clamped}

\new{\subsubsection{Initial stiffness}}

The second simple benchmark is a beam of span $L$
clamped at both ends and subjected to a uniformly distributed vertical load of intensity $f$.
According to the linear Timoshenko theory, the \new{normalized} mid-span deflection would be
\old{
\beq 
w = \frac{fL^4}{384EI} + \frac{fL^2}{8GA_s}
\eeq 
which can be written in terms of the dimensionless quantities as
\beq 
\frac{w}{L} = 
\frac{fL^3}{EI}\frac{1}{384}\left(1 + \frac{144EI}{EAL^2}\right)
= 
\frac{fL^3}{EI}\frac{1}{384}\left(1 + \frac{12h^2}{L^2}\right)
\eeq 
}
\new{
\beq 
\frac{w}{L} = 
\frac{fL^3}{384EI}+ \frac{fL}{8GA_s}
= 
\frac{fL^3}{EI}\frac{1}{384}\left(1 + \frac{12h^2}{L^2}\right)
\eeq 
}
Here, the load is characterized by the
dimensionless variable $fL^3/EI$, \new{and we still assume a rectangular cross section and Poisson's ratio $\nu=0.25$.}

For a Kirchhoff beam, the dimensionless stiffness would be 384.
The effect of shear reduces it by the factor $(1+12h^2/L^2)$.
Let us consider two cases with the depth-to-span ratio 1/6 and 1/12.
The corresponding dimensionless structural stiffnesses are 288 and 354.4615, respectively.

\begin{figure}[h]
\hskip -20mm
    \begin{tabular}{cc}
    (a) & (b)\\ 
    \includegraphics{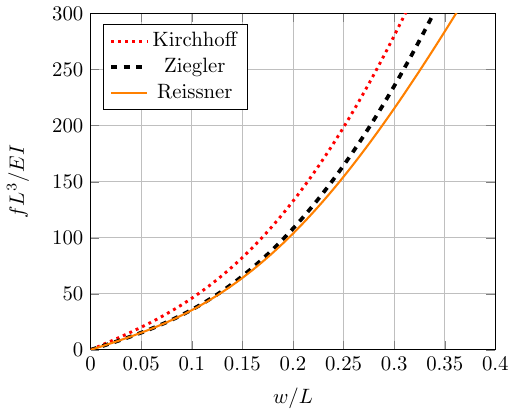}
    &
    \includegraphics{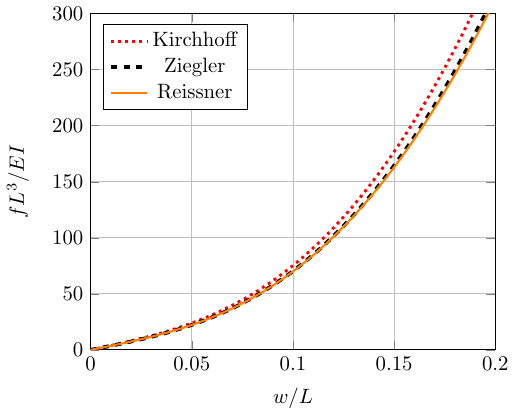}
    \end{tabular}
\caption{Clamped beam subjected to uniform load of intensity $f$: normalized load-displacement curves (a) for depth-to-span ratio 1:6, (b)~for depth-to-span ratio 1:12} 
\label{ex4_2}
\end{figure}

Since the applied load is distributed and applied at the beam level,
the initial stiffness cannot be directly deduced from the structural stiffness. A simple check consists in calculating the mid-span deflection 
caused by a very small load. For $h/L=1/6$, we set the dimensionless
load to $fL^3/EI=0.0288$. The deflection expected for a linear model
is then $w/L=10^{-4}$. Indeed, the numerical simulation with two elements
and 10,000 integration segments per element leads to a deflection that
differs from this value approximately by $2\cdot 10^{-11}$, so the relative 
error is only 0.00002~\% (in fact, this difference is probably caused by the 
nonlinearity). With 16, 32 and 64 integration segments, the relative errors
are found to be 0.59~\%, 0.15~\% and 0.037~\%, respectively, which confirms the
expected quadratic dependence of the error on the segment size. 
Similar results are obtained for the beam with $h/L=1/12$, where
a dimensionless load $fL^3/EI=0.03544615$ is expected to lead to 
$w/L=10^{-4}$. Numerical simulations with 16, 32 and 64 integration segments lead to relative errors of 0.72~\%, 0.18~\% and 0.045~\%, respectively. \old{These results} \new{The relative errors} do not depend on the adopted formulation
(Reissner or Ziegler).

\new{\subsubsection{Load-displacement diagrams}}

\begin{figure}[b!]
 \hskip -20mm
    \begin{tabular}{cc}
    \includegraphics[scale=1]{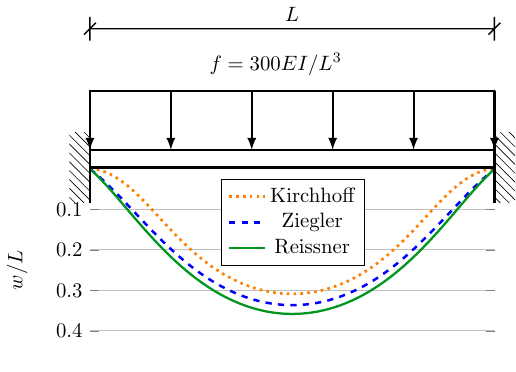}
    &
    \includegraphics[scale=1]{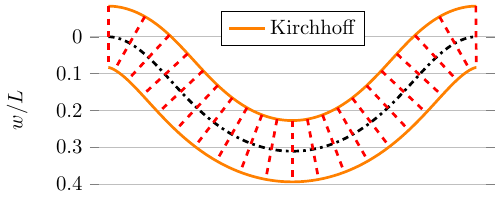}
    \\
    \includegraphics[scale=1]{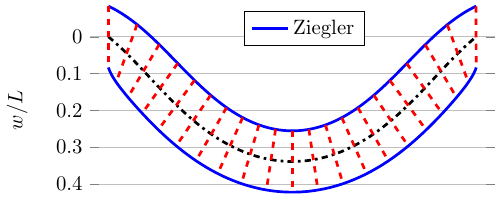}
    &
    \includegraphics[scale=1]{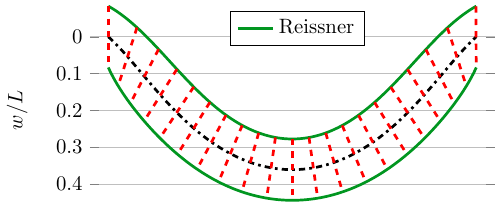}
    \end{tabular}
\caption{Clamped beam subjected to uniform load of intensity $f$: initial shape, boundary conditions and loading, and the deformed centerlines and rotated sections at load level $f=300\,EI/L^3$ obtained for depth-to-span ratio 1:6 using various models} 
\label{ex4_2_beams}
\end{figure}

The complete load-displacement diagrams that cover the highly nonlinear
range are plotted in   Fig.~\ref{ex4_2}. In addition to the
results obtained with the Reissner and Ziegler models, the figure
includes the response of the Kirchhoff model, i.e., a beam model without shear but with
axial extensibility. For the beam clamped at both ends, the axially inextensible Euler
model would not allow any changes of the centerline and thus
cannot be used (except if the changes are considered as infinitesimal, in which case both Euler and Kirchhoff formulations reduce to the standard linear beam theory). For the stocky beam with $h/L=1/6$, the difference between the Reissner and Ziegler formulations is more pronounced than in the previous example of a simply supported beam; see also Fig.~\ref{ex4_2_beams} with the deflected shapes
that correspond to the load level of $f=300\,EI/L^3$.
For $h/L=1/12$, both formulations yield very similar results,
but the effect of shear is still non-negligible.


\new{\subsubsection{Discretization error in nonlinear range}}

As an additional convergence check, Table~\ref{tab:3} shows the dimensionless mid-span deflections obtained numerically for the dimensionless
load $fL^3/EI=300$ and the corresponding relative
errors. The calculations are done for the stocky beam with $h/L=1/6$,
for which the Reissner and Ziegler models yield different results.
The results indicate that, for the same computational grid, 
the Reissner model gives typically a lower error than the Ziegler model. Asymptotically (for fine grids), 
the error is for both models proportional
to the square of the segment size. Some irregularity is observed
for very coarse grids (2 or 4 segments per element), for which the Reissner model yields
surprisingly high accuracy.

\begin{table}[htb]
     \centering
  \caption{\textcolor{black}{Clamped beam with $h/L=1/6$ subjected to uniform load: convergence of the mid-span displacement obtained for a given load $f=300\,EI/L^3$ as the number of integration segments per element is increased\\}}
\label{tab:3}
     \begin{tabular}{rllll}
     \hline
     
number of& \multicolumn{2}{c}{\textcolor{black}{Reissner model}} & \multicolumn{2}{c}{\textcolor{black}{Ziegler model}} \\
segments & $w/L$ & error $[\%]$ & $w/L$ & error $[\%]$\\
          \hline
2	 &  0.365004  &  1.09 &0.355012 & 4.75 \\
4	 & 0.362076  &  0.28 & 0.344361 &  1.61 \\
8	 &  0.361554  & 0.13  &  0.340514  & 0.48 \\
16	 &  0.361227  &  0.045 &  0.339328  & 0.13  \\
32	 &  0.361109  &  0.012 &  0.339011  & 0.032  \\
64	 &  0.361077  &  0.0032 & 0.338930   &  0.0078 \\
128	 &  0.3610683  &   0.0008 &  0.3389102  & 0.0019  \\
256	 &  0.3610662  & 0.0002  &   0.3389052  & 0.0005  \\
\new{converged} &  0.3610655  &  & 0.3389035  & \\
\hline
     \end{tabular}
 \end{table}

 \new{\subsubsection{Internal forces}}

\new{
The adopted numerical scheme provides good approximations
not only for kinematic variables such as displacements and
rotations, but also for static variables such as internal forces.
To demonstrate that and to illustrate the difference between the 
two considered formulations, Fig.~\ref{fig:intforces} shows the distribution
of the normal forces and shear force along the beam, plotted as 
functions of the initial distance from the left end, $\xi$.
For both internal forces, solid curves correspond to the almost converged
solution and dashed curves to the results obtained with
just 8 integration segments  at load level $f=300\,EI/L^3$. The forces are normalized by $EI/L^2$.
The values of internal forces are plotted at midpoints of integration segments and, for the accurate solution with 32 segments per element, connected by straight lines. Even for 8 segments per element, the error is seen to be quite small. On the other hand, the differences between the Reissner (blue) and Ziegler (red) internal forces are substantial,
especially near the supports, where the shear 
distortion is high (cf.\ Fig.~\ref{ex4_2_beams}). 
Recall that for the Reissner model, the normal force is interpreted as the component perpendicular to the section, while for the Ziegler model it is aligned with the tangent to the centerline).

\begin{figure}[]
 \hskip -20mm
    \begin{tabular}{cc}
    \includegraphics[scale=0.95]{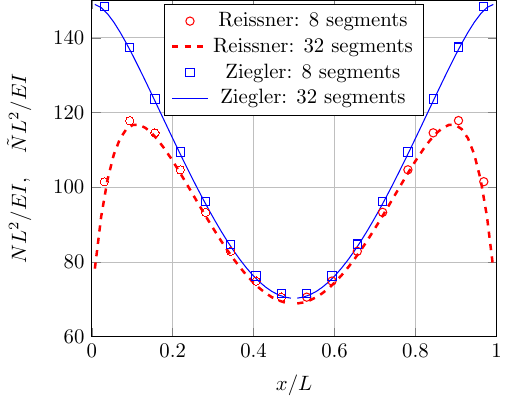}
    &
   \includegraphics[scale=0.95]{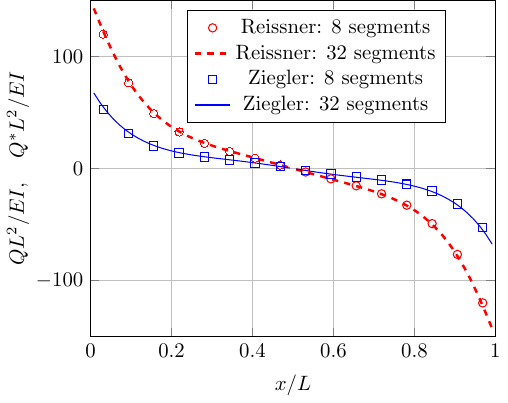}
    \\ (a) & (b)
    \end{tabular}
\caption{\new{Clamped beam subjected to uniform load of intensity $f$:
distribution of (a) normal force and (b) shear force at load level $f=300\,EI/L^3$ obtained for depth-to-span ratio 1:6 with the Reissner (red) or Ziegler (blue) model using 8 (points) or 32 (curves) integration segments per element}} 
\label{fig:intforces}
\end{figure}    

}

\subsection{Cantilever, uniformly distributed external moment}
\label{sec:example3}

The previous example demonstrated that the effect of distributed forces is properly captured by the developed algorithm. Let us now
investigate a somewhat artificial yet interesting case of a cantilever subjected to an externally applied distributed moment of uniform intensity $m$. The beam can be expected to fold
into a spiral-like shape. 
Since the normal and shear forces vanish
(independently of the specific definition),
shear plays no role here and the centerline
does not change its length. 
The problem can be treated analytically \new{(see~\ref{appB})}, which
provides a reference solution for the
numerical simulations. 

\begin{figure}[ht]
\hskip -5mm
    \begin{tabular}{cc}
    (a) & (b) 
    \\
\includegraphics[scale=1]{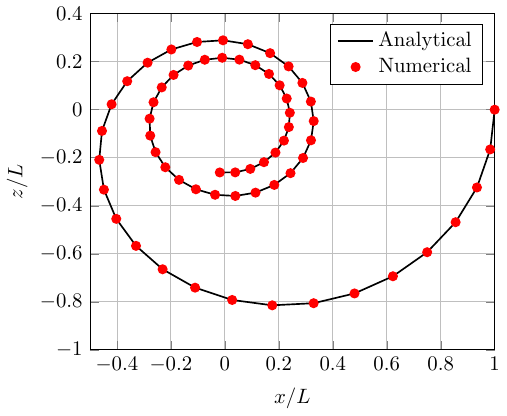}
    &  \includegraphics{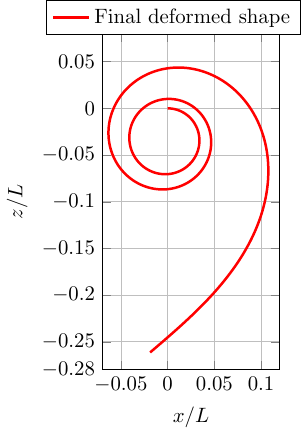}
    \end{tabular}
\caption{Cantilever loaded by a uniformly distributed external moment: (a) the trajectory of the cantilever tip when the applied moment intensity increases from 0 to $30\,EI/L^2$, (b) the
final deformed shape} 
\label{ex4_4}
\end{figure}

Numerical simulation of the cantilever with uniform distributed moment load can be done using a single beam element, which means that only
three global unknowns are introduced. The intensity of the moment load
is increased in 60 equal steps up to the final value of $m=30 EI/L^2$,
which corresponds to the dimensionless parameter $\mu=\sqrt{30/\pi}\approx 3.0902$.
The left end of the beam, placed at the origin, is clamped,
and the right end, placed at $(x,z)=(L,0)$, is free.  
Fig.~\ref{ex4_4}a shows the numerically evaluated position of the cantilever tip after
each step (filled circles) and the analytical solution 
parametrically described by (\ref{eq205})--(\ref{eq206}).
The final deformed shape is depicted in Fig.~\ref{ex4_4}b.
The numerical solution divided the element into 500 integration segments, and the calculated position of the cantilever tip
after the last step deviates from the exact value by $6.14 \cdot 10^{-6}$ times
the beam length $L$. 

\subsection{Shallow frame dome with snap-through}
\begin{figure}[ht]
\hskip -5mm
    \begin{tabular}{cc}
    (a) & (b) \\
    \includegraphics[scale=1.0]{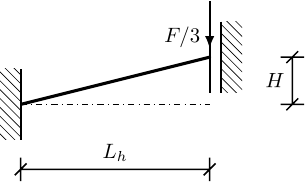}
    &
    \includegraphics[scale=1]{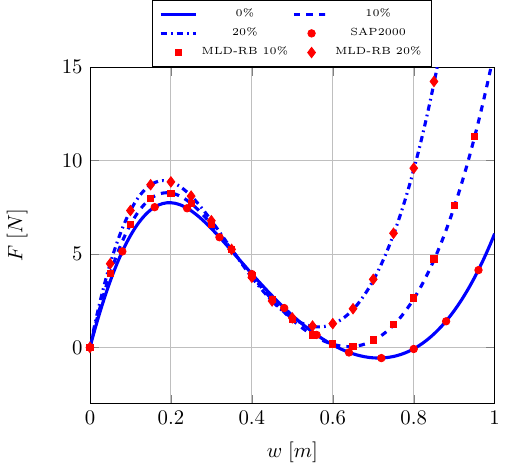}
    \end{tabular}
\caption{Shallow frame dome: (a) loading and supports of a planar beam
that replaces 1/3 of the spatial structure,
(b) load-displacement diagrams for standard beams (solid) and beams with 10~\% (dashed) and 20~\% (dash-dotted) rigid offsets; the isolated points correspond to the results obtained
by Di Re et al.\ using SAP2000 (circles) and their MLD-RB elements with offsets modeled as rigid beams (squares and diamonds)} 
\label{ex4_4b}
\end{figure}

For certain applications, it is useful to provide optional
connections of beam elements to the joints by rigid arms. 
For instance, if the nodes of the beam model are placed at
the centers of the actual physical joints of finite dimensions,
the theoretical length of the beam element measured as the
distance between nodes is longer than the actual beam segment
that has the given cross sectional characteristics, and the
response of the model may become more flexible than the actual
behavior. As demonstrated by Di Re et al.\cite{DiRe2024}, the experimental
load-displacement curve measured on a 3D-printed spatial lattice
structure by Jamshidian et al.\cite{JamBodRos20}
can be numerically reproduced by a beam finite element model
only if the end portions of the beams are made stiffer than the
regular part. 
Therefore, Di Re et al.\cite{DiRe2024} developed three-dimensional corotational 
finite elements with rigid joint offsets. 
Among other examples, they tested their formulation on a shallow
frame dome, previously studied by Battini et al.\cite{BatPacEri03}.
Even though the considered structure is spatial, its symmetry 
makes it possible to perform the analysis using a planar beam element.
The dimensions and supports are shown in Fig.~\ref{ex4_4b}a,
with $L_h=15$~m and $H=0.6$~m. The cross section is rectangular,
of width $b=0.14$~m and depth $h=0.17$~m, and the material is 
considered as linear elastic, with Young's modulus $E=10$~MPa and Poisson's ratio $\nu=0.3$. 
The frame dome consists of three
beams located in vertical planes rotated with respect to each
other by 120$^\circ$. If it is loaded by a vertical force at the 
common joint, the resulting diagram exhibits snap-through. 

To get a reference solution, Di Re et al.\cite{DiRe2024} analyzed the problem
using the SAP2000 FE package, with each beam divided into 64 finite
elements. Selected points in the load-displacement diagram
obtained in this way are shown as filled circles in Fig.~\ref{ex4_4b}b. 
The solid curve corresponds to our simulation with a single
beam element divided into 100 integration segments.
Very good agreement is observed. The largest deviation
in terms of force is $0.004$~N, while the snap-through force is $7.7484$~N. 

Subsequently, Di Re et al.\cite{DiRe2024} studied the effect of increased stiffness
near the beam ends. They considered rigid offsets over 
10~\% or 20~\% of the total beam length near both ends. 
For the present comparison, 
we use the results obtained with their MLD-RB;
see the squares and diamonds in Fig.~\ref{ex4_4b}. 
The continuous curves represent our simulations with a single
beam element divided into 100 integration segments.
The agreement is again very good.
The largest deviation for 10~\% offsets is $0.0255$~N, while the snap-through force is $8.2827$~N, and the largest deviation for 20~\% offsets is $0.00256$~N, while the snap-through force is $8.9288$~N.

Modeling of rigid offsets within the present framework
requires only minimal adjustments. It is sufficient to set
the flexibilities $1/EI$, $1/EA$, and $1/GA_s$ to zero in all
integration segments that are located in the rigid part of the beam.
This is indeed the simplest modification that requires minimal
programming effort. Of course, computational efficiency can be
increased by handling all rigid integration segments at once
and skipping all unnecessary operations.
The modified code replaces the algorithmic steps (\ref{e234mj3})--(\ref{e240mj3}) of the Reissner model or 
(\ref{ee215})--(\ref{ee231}) of the Ziegler model
by
\bea 
x_i &=& x_{i-1}+\Delta\xi_R \cos\varphi_{i-1}\\
z_i &=& z_{i-1}-\Delta\xi_R \sin\varphi_{i-1}\\
\varphi_i &=& \varphi_{i-1} 
\eea 
where $\Delta\xi_R$ is the total length of the rigid segment
(which does not need to be subdivided into integration
segments of the same size as the flexible part of the beam).
In the presence of loads acting directly on the beam,
it is also necessary to evaluate
\beq 
M_{pi} = M_{p,i-1}- \left(m(\xi_{i-1/2})+P_{x,i-1/2}\cos\varphi_{i-1}+ P_{z,i-1/2}
\sin\varphi_{i-1}\right)\Delta\xi_R
\eeq

\section{Examples---stability}

When developing the theoretical background in Section~\ref{sec:beam model}\new{ and \ref{appA}},
we presented two alternative versions of sectional equations---one corresponding
to the widely used Reissner approach (\ref{sec:2.2.2}\new{ and \ref{appA1}}) and the other that can be referred to as the Ziegler approach (Section~\ref{sec232}\new{ and \ref{appA2}}). In the geometrically linear approximation, both reduce to the same
standard Timoshenko theory. Also, when the shear deformation is neglected, or the sectional shear stiffness is set to infinity, they both reduce to what we consider as the Kirchhoff beam, even for large rotations.  However, differences in the
load-displacement curves are observed when large shear distortions occur; see e.g.\ Fig.~\ref{ex4_2}a.

Let us now investigate how the choice of a specific formulation of the sectional equations
affects the buckling load of an axially loaded straight column.
\old{It will be shown} \new{In~\ref{appC} it is shown} that the Reissner model is closely
related to the Haringx approach to stability analysis\cite{Haringx1942}, while the Ziegler model is an extension of the
Engesser approach\cite{Engesser1891} to the case of axially compressible columns. 
Analytical expressions for the critical load \old{will be} derived for both formulations \new{in~\ref{appC}} \old{, and the results will be compared with the known solutions described, e.g., in the Stability of Structures textbook by Ba\v{z}ant and Cedolin\cite{Bazant91Cedolin}.
This will generate} \new{can be used here as} additional benchmarks
for the numerical procedures.

\subsection{Critical state of a clamped column in compression}
\label{sec:ex6.1}

The ability of the developed formulations to capture compressive buckling will be demonstrated using the example of a column clamped on both sides, with depth-to-span ratio 1:6 or 1:12. 
These are the same geometries as used in the example in Section~\ref{sec:clamped},
but this time the lateral loading is replaced by an axial displacement
that is prescribed at one of the supports and initially induces
uniform compression. 
When the reaction force reaches the critical level, buckling is expected. The critical state can be numerically detected by
analyzing the eigenvalues of the structural stiffness matrix.

The described loading could in principle be simulated by a single
beam element, but then there would be no global unknowns and
the bifurcation would occur at the element level.
To capture it using traditional techniques,
let us model the column using two equally sized beam elements. 
The central node has three degrees of freedom that are used
as global unknowns. 

Since the column is clamped at both ends, its buckling length
is $L_b=L/2$, and the corresponding slenderness can be evaluated
as $s=L_b/i_y=L/(2\sqrt{I/A})=\sqrt{3}L/h^2$. 
Therefore, the square of slenderness is $s^2=3L^2/h^2= 108$
for the stocky case with $h/L=1/6$, and $s^2=432$ for the more
slender column with $h/L=1/12$. Moreover, the sectional stiffness
ratio $\Gamma=GA_s=EA=1/3$, provided that the Poisson ratio is set to 0.25
and the shear reduction coefficient is considered as $A_s/A=5/6$
for a rectangular section. Substituting into the analytical
formulae derived in the previous section,
we obtain the critical loads that correspond to various models.
The results are summarized in Table~\ref{tab:critloads} in the form
of critical strains, defined as $\eps_{crit}=P_{crit}/EA$.
The case of an incompressible model without shear corresponds
to the standard Euler solution, with critical stress $\sigma_{crit}=P_{E}/A=E\pi^2/s^2$.
The critical strain is formally obtained as the critical stress
divided by the elastic modulus, i.e., it is given by $\pi^2/s^2$, but should not be interpreted
as the actual strain at the onset of buckling, since this particular model
assumes that the centerline does not change its length.
For the Reissner model, the critical strain $\eps_{Rc}=\sigma_{Rc}/E$ is evaluated using formula (\ref{282e}) with $E$ removed. 
For the Ziegler model, formula (\ref{eq:epsZc}) directly provides
the critical strain
$\eps_{Zc}$.

\begin{table}[h]
    \caption{Critical strains (compressive) for the clamped-clamped column according to various models (numerical results computed
    with 2 beam elements and 32 integration segments per element, based on linear interpolation of the lowest eigenvalue over a prescribed strain increment of 0.001)\\}
    \centering
    \begin{tabular}{lllll}
    \hline
         &  $h/L=1/6$ && $h/L=1/12$ & \\
         & analytical & numerical  & analytical & numerical \\
      \hline
       Euler  & 0.091385 & 0.091312 & 0.022846 & 0.022828\\
       Kirchhoff ($GA_s\to\infty$)  & 0.101700 & 0.101643 & 0.023394 & 0.023374 \\
       Reissner ($GA_s/EA=1/3$)  & 0.078926 & 0.078870 & 0.021888 & 0.021871 \\
       Ziegler ($GA_s/EA=1/3$) & 0.077770 & 0.077716 & 0.021859 & 0.021842 \\ \hline
    \end{tabular}
    \label{tab:critloads}
\end{table}

In numerical simulations, the column is discretized into two beam elements of length 0.5 that connect nodes 1 and 2 and nodes 2 and 3. 
Node 3 is fixed (including rotation), and the axial displacement
at node 1 is increased in each incremental step by 0.001. 
For the fundamental solution, the axial displacement at node 2
increases in each step by 0.0005 while the lateral displacement
and the rotation remain zero. The structural tangent stiffness
is in this example diagonal, and the diagonal stiffness coefficient that corresponds
to the lateral displacement at node 2 gradually decreases.
The simulation is terminated when this coefficient becomes negative,
and the critical strain is evaluated by linear interpolation in the
last incremental step. For instance, for the Reissner model and
 $h/L=1/6$ (simulated with $EI=1$, $EA=432$ and $GA_s=144$), analyzed using 32 integration segments per element,
the diagonal stiffness coefficient is equal to 1.817168 after step 78
and to $-0.270419$ after step 79, and the critical strain is estimated
as $0.078+0.001\times 1.817168/(1.817168+0.270419)=0.078870$.
This is only by 0.07~\% below the analytical value of 0.078926.
For a simulation with 8 integration segments per element,
the numerically estimated critical strain would be 0.078037, which is by $1.1~\%$ below the analytical value.

Using a similar procedure, the numerical estimates have been constructed
 for the other cases (using 32 integration segments per element),
and the results are reported in Table~\ref{tab:critloads}.
Excellent agreement with the analytical results is observed.
For the incompressible model without shear (Euler), the numerical
solution cannot be driven by a prescribed axial displacement.
Instead, the applied axial load is incremented,
and the critical ``strain'' is evaluated as the critical load
divided by the real physical stiffness $EA$, while the numerical
simulation is performed with the axial sectional flexibility $1/EA$ set to zero.

\subsection{Post-buckling response in compression}

\begin{figure}[h]
\hskip -20mm
\begin{tabular}{cc}
 (a) & (b) \\
\includegraphics[scale=1]{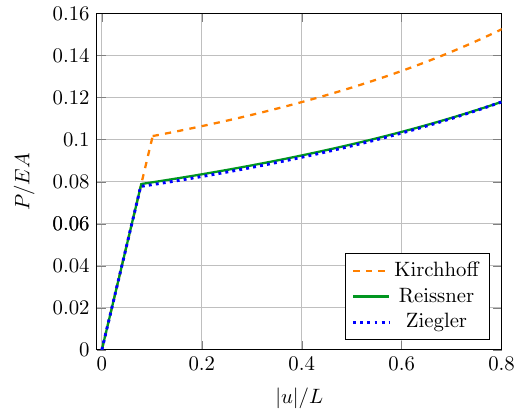}
&  
\includegraphics[scale=1]{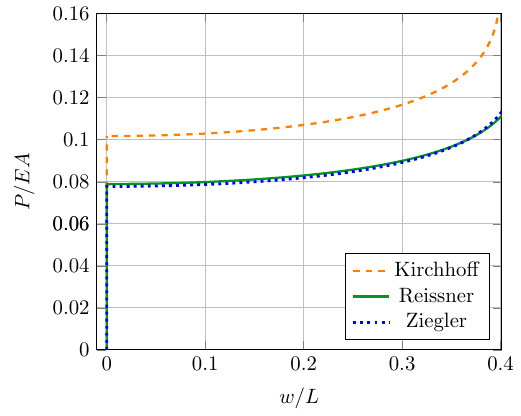}
\end{tabular}
\caption{Load-displacement diagrams for a compressed column clamped at both ends: (a) axial displacement versus applied force, 
(b) lateral displacement versus applied force
 }
    \label{fig:reissner-postcritical}
\end{figure}

Numerical simulations make it possible not only to detect the critical state,
but also to trace the bifurcated solution in the post-critical range.
The resulting diagrams in Fig.~\ref{fig:reissner-postcritical} show the 
evolution of the axial and lateral displacement
plotted against the loading force. 
We consider here a column clamped at both ends, 
with depth-to-length ratio $h/L=1/6$, rectangular cross section, and sectional stiffness ratio
$GA_s/EA=1/3$. The corresponding slenderness is $s=\sqrt{3}\,L/h\approx 10.392$.
The simulation has been done with 2 beam elements, each divided into
100 integration segments. In the graphs, the force is normalized
by the axial sectional stiffness $EA$, and the displacements are normalized
by the beam length $L$. The axial displacement $u$ is taken at the loaded
section and thus represents the axial shortening of the column.
The lateral displacement $w$ is taken at the computational node in the
middle of the column.
Before
bifurcation, the bar is compressed uniformly and the normalized
axial displacement, $|u|/L$, is equal to the normalized force, $P/EA$,
since both represent the strain magnitude. 
The transversal displacement remains equal to zero.

The theoretical value of the critical strain at which the bifurcation
is expected can be calculated using formula (\ref{282e}) or 
(\ref{eq:epsZc}), depending on which model is used. 
Since the considered geometry and stiffness ratio correspond to
the example in Section~\ref{sec:ex6.1}, the values of the critical strain can be found directly in Table~\ref{tab:critloads}. 
If the simulation is run without any perturbations, the solution
will follow the unstable fundamental solution. To trigger the bifurcation,
the solution is perturbed by a small lateral load applied at the 
computational node in the middle of the column. Once the equilibrium
iteration converges, the perturbation load is removed and the equilibrium iteration continues. If this is done in a subcritical state, the final
solution corresponds to the fundamental one, since this is the only solution
at the given load level. However, if the load is supercritical,
the process typically converges to one of the stable solutions on the
bifurcated branch.

 In this way, the equilibrium diagrams have been constructed for the
 Reissner and Ziegler models (for which the critical strains are
 $7.89~\%$ and $7.78~\%$, respectively). Only small differences
 in the response are observed in Fig.~\ref{fig:reissner-postcritical}. For comparison, the response has also been traced
 using the Kirchhoff model, for which the critical strain 
of 10.2~\% is substantially larger than for the models
that take into account compliance in shear. At comparable displacement levels in the post-critical range, the Kirchhoff model also requires higher applied loads.  
 
\begin{figure}[h]
\hskip -10mm
\begin{tabular}{ccc}
 \includegraphics[scale=1]{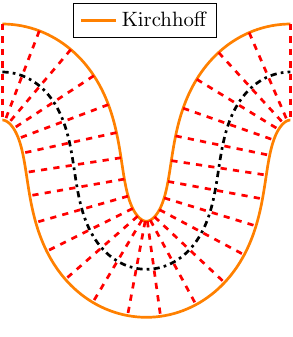}
 &
\includegraphics[scale=1]{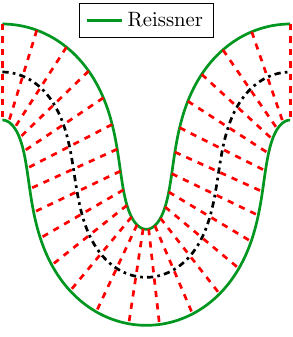}
&  
\includegraphics[scale=1]{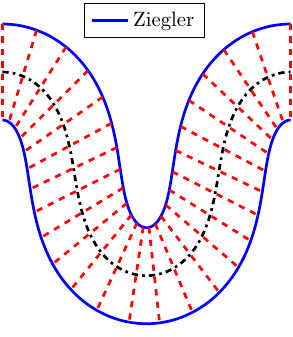}   
\end{tabular}
\caption{Deformed shapes in the post-critical state at $|u|/L=0.5$ obtained with various models }
    \label{fig:reissner-postcritical2}
\end{figure}

The deformed shapes in the post-critical range are visualized in
Fig.~\ref{fig:reissner-postcritical2}, which shows not only the deformed centerline
but also the inclination of selected cross sections, in order to
demonstrate the role of shear.
All the displayed shapes correspond to the state in which 
the axial displacement corresponds to one half of the initial column length (i.e., $|u|/L=0.5$ in Fig.~\ref{fig:reissner-postcritical}a).
The shapes of the deformed centerline are similar for all three models, but the sectional inclinations obtained with the Kirchhoff model are different
from those obtained with the models that account for shear, and the corresponding axial force is also different. 
The diagrams in Fig.~\ref{fig:reissner-postcritical}
cover the numerical solution of the governing equations in a wide range,
but extreme states with $|u|/L$ above 0.66 are nonphysical, because of internal contact that would occur in the folded column.

\subsection{\new{Critical state in tension}}

As shown in~\ref{sec:teninstab}, the Reissner model can lead to
bifurcations in tension, and the analysis needs to
be done separately for different types of supports. Therefore, 
illustrative numerical simulations of tensile bifurcations  will
be performed for the three combinations of boundary conditions analyzed in~\ref{sec:teninstab}.

First consider a bar of length $L$ {\bf clamped at one end}
and loaded at the other (free) end by a force that induces axial tension. 
This problem always has a fundamental solution that
corresponds to uniform stretching of the beam, with the
axis remaining straight. 
According to the theoretical analysis, the fundamental solution can be expected to lose stability when a certain
critical state is attained. The critical strain can be determined using formula (\ref{eq258}) and its value depends
on two dimensionless parameters: the sectional stiffness ratio $\Gamma$ and the slenderness $s$. 
For a cantilever of length $L$, the slenderness needs to be evaluated
using the buckling length $L_b=2L$. If the section is rectangular, of depth $h$, we get
\beq 
s = \frac{L_b}{\sqrt{I/A}} = \frac{2L}{\sqrt{h^2/12}} =4\sqrt{3}\,\frac{L}{h}
\eeq 
Let us examine two depth-to-span ratios, $h/L=1/4$ and $1/8$,
and three values of the stiffness ratio $\Gamma$,
with $\Gamma=1/3$ characterizing a homogeneous rectangular
section, and $\Gamma=0.1$ and 0.01 corresponding to sections
with low and very low shear stiffness.
The resulting analytical values of the critical strain
evaluated from (\ref{eq258}) are compared with the 
numerically obtained values in Table~\ref{tab:critstrainten}. 

The bar is modeled by 1 element with 32 integration segments. 
The difference between the numerical and analytical results in Table~\ref{tab:critstrainten} is partly due to the linear interpolation of the lowest eigenvalue of the stiffness matrix over the step in which this eigenvalue changes its sign from positive to negative. The size of this step
in terms of the prescribed strain increment has been fixed to 0.001. Refining the step (e.g., using the secant method for finding the root) could further reduce the error,
even at a fixed number of integration segments. 
However, the agreement is excellent even without this improvement.

\begin{table}[h]
   \caption{Critical strains in tension for a bar of depth $h$ and length $L$ according to the Reissner model with various stiffness ratios $\Gamma=GA_s/EA$   (numerical results are based on linear interpolation of the lowest eigenvalue over a prescribed strain increment of 0.001)\\}
    \label{tab:critstrainten}
    \centering
    \begin{tabular}{rllllll}
    \hline
    & \multicolumn{4}{c}{clamped at one end} & \multicolumn{2}{c}{clamped at both ends} \\ 
         &   \multicolumn{2}{c}{$h/L=1/4$} &  \multicolumn{2}{c}{$h/L=1/8$} &  \multicolumn{2}{c}{$h/L=1/6$} \\
     $GA_s/EA$    & analytical & numerical  & analytical & numerical& analytical & numerical \\
      \hline
       1/3  & 0.512537 & 0.512529 & 0.503192 & 0.503179 & 0.522385 & 0.522365\\
       1/10  & 0.122744 & 0.122736 & 0.114636 & 0.114219 & 0.132805 &  0.132788\\
       1/100  & 0.017513 & 0.017502 & 0.012664 & 0.012628 & 0.026976 &0.026969\\ \hline
    \end{tabular}
\end{table}

Next, let us look at a {\bf simply supported bar}. The theoretical analysis
indicates that, in this particular case, slenderness does not affect 
the critical strain given by formula (\ref{eq:epsrtss}).
For $\Gamma=1/3$, 0.1 and 0.01, the analytically evaluated critical
strain $\eps_{Rt}^{(ss)}$ is 1/2, 1/9 and 1/99, respectively. 
When the critical strain is imposed to the one-element numerical model of the bar, the tangent stiffness
matrix has a zero eigenvalue (up to the machine precision), independently
of the number of integration segments. The bifurcation is captured exactly,
because the bifurcation mode is very simple (zero deflection and uniform
rotation) and thus can be reproduced by the numerical approximation,
no matter how many integration segments are used.

Finally, consider a {\bf bar clamped at both ends}, this time with
$h/L=1/6$. Since the buckling length is now $L_b=L/2$, the slenderness
is evaluated as $s=\sqrt{3} L/h\approx 10.392$. The critical strains $\eps_{Rt}^{(dc)}$ obtained by solving equation (\ref{eqs247y}) are
listed in the penultimate column of Table~\ref{tab:critstrainten}.
Numerically, the bar is modeled by 2 elements, each divided into
16 integration segments (for a one-element model, no global degrees of freedom would be
introduced and the tangent stiffness could not be examined).
The imposed strain is again increased in increments of 0.001
and linear interpolation of the smallest eigenvalue is used.
The agreement is excellent.

\begin{figure}[h]
\hskip -20mm
\begin{tabular}{cc}
(a) & (b) \\
\includegraphics[scale=1]{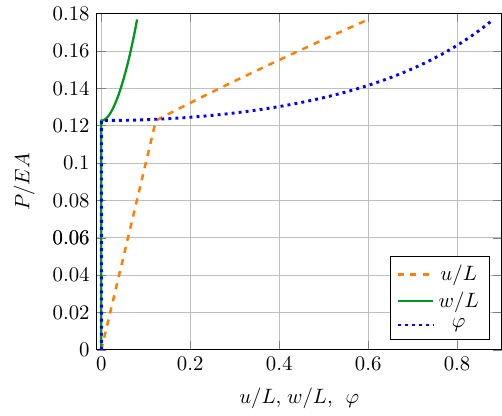}
& 
\includegraphics[scale=1]{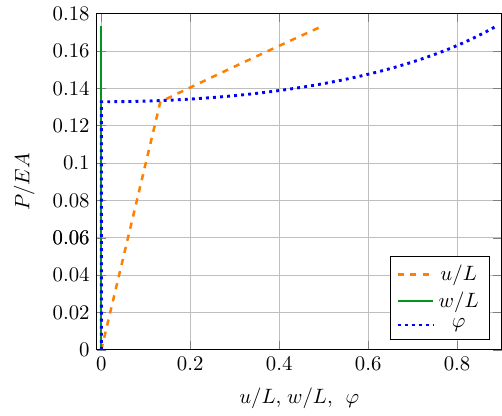}
\end{tabular}
\caption{Load-displacement diagrams under tension obtained with the Reissner model for a bar (a) clamped at the left end, $h/L=1/4$, (b) clamped at both ends, $h/L=1/6$}
    \label{f:ex6_4diagrams}
\end{figure}

\begin{figure}[h]
\centering
\begin{tabular}{ll}
(a) & \hskip 2mm\includegraphics[scale=1.2]{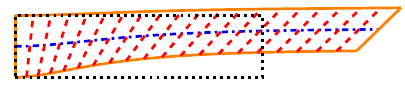}
\\ 
(b) & \includegraphics[scale=1.34]{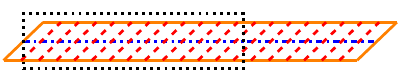}
\\ 
(c) & \hskip 2mm\includegraphics[scale=1.6]{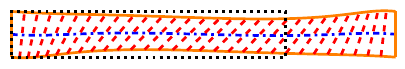}
\end{tabular}
\caption{Deformed shapes in post-critical states under tension obtained with the Reissner model for a bar (a) clamped at the left end, $h/L=1/4$, (b) simply supported, $h/L=1/4$, (c) clamped at both ends, $h/L=1/6$ (the dotted rectangles represent the undeformed shapes)}
    \label{f:ex6_4shapes}
\end{figure}

\subsection{\new{Post-critical response in tension}}

Numerical simulations make it possible not only to detect the critical state,
but also to trace the bifurcated solution in the post-critical range.
The resulting diagrams in Fig.~\ref{f:ex6_4diagrams} show the 
evolution of selected degrees of freedom (nodal displacements and rotation)
plotted against the applied loading force, $P$. By the adopted convention,
negative $P$ corresponds to tension, and thus the diagrams work with
the absolute value of $P$, which represents the axial tensile force
that would arise if the bar remained straight. On the vertical axis,
the force is normalized by the sectional stiffness, $EA$. Before
bifurcation, the bar is stretched uniformly and the normalized
axial displacement, $u/L$, is equal to the normalized force, $|P|/EA$,
since both represent the strain. Here, $u$ is the axial displacement of the
end section subjected to the external loading.  For the bar clamped
at the left end, the other quantities plotted in Fig.~\ref{f:ex6_4diagrams}a are the transversal 
displacement $w$ (also normalized by $L$) and the rotation $\varphi$
at the free right end. For the bar clamped at both ends, 
the transversal displacement and rotation plotted in Fig.~\ref{f:ex6_4diagrams}b are traced at the middle section.
In this case, the transversal displacement remains equal to zero even
after the bifurcation, because the deformed shape is skew-symmetric. 
The deformed shapes in the post-critical range are visualized in
Fig.~\ref{f:ex6_4shapes}, which shows not only the deformed centerline
but also the inclination of selected cross sections, in order to
demonstrate the predominant role of shear.
The shapes correspond to a stage at which the maximum rotation is 
about 0.8. This is quite far beyond the bifurcation point, but the
distribution of lateral displacement and rotation along the bar
is still similar to the predicted bifurcation modes analyzed
at the end of \ref{sec:teninstab}.

The visualization of deformed shapes makes it easier to understand the
nature of the tensile bifurcation. If the rotations are unconstrained,
the segments exhibit uniform shear distortion (Fig.~\ref{f:ex6_4shapes}b). As neighboring segments 
slide with respect to each other, the centerline gets extended and the applied load can supply work, which is transformed into strain energy.
If some constraints on rotations are imposed at one or both end sections,
the region of the bar that is far from the support(s) develops a pattern
that is close to uniform shear, but the boundary region(s) are affected
by the constraint(s) and the load needed to drive the deformation
is higher. For slender beams, the affected boundary regions are relatively small and the critical load is only slightly higher than for the 
uniform shear mode, while for stocky beams the boundary regions are
more important and the critical load can be substantially higher,
as seen in Fig.~\ref{6a}.

Another interpretation of the bifurcation is that if the sections rotate
while the bar axis remains straight,
the applied force acting along the bar axis is according to the Reissner model
decomposed into a normal force $N$ perpendicular to the section
and a shear force $Q$ parallel to the section; these internal forces are 
expressed as $N=|P|\cos\varphi$ and $Q=|P| \sin\varphi$, and 
they lead to deformation measures $\eps=N/EA$ and $\gamma=Q/GA_s$.
Since the centerline remains straight, the rotation $\varphi$ is 
equal to the shear angle $\chi$, which is linked to the deformation measures
by
\beq 
\tan\chi = \frac{\gamma}{1+\eps}
\eeq 
Making use of the relations mentioned above, we can set up the equation
\beq 
\tan\varphi = \frac{|P|\sin\varphi/GA_s}{1+|P|\cos\varphi/EA}
\eeq 
If the trivial case characterized by $\sin\varphi=0$ is excluded, the
equation can be rewritten as
\beq 
|P|\left(\frac{1}{GA_s}-\frac{1}{EA}\right)\cos\varphi =1
\eeq 
A non-trivial solution $\varphi\ne 0$ exists only if 
\beq 
|P|\left(\frac{1}{GA_s}-\frac{1}{EA}\right)> 1
\eeq 
which can be written as $|P|>P_{Rt}^{(ss)}$ where
\beq 
P_{Rt}^{(ss)} = \frac{1}{\dfrac{1}{GA_s}-\dfrac{1}{EA}} = EA \,\frac{1}{\dfrac{1}{\Gamma}-1}
= EA \,\frac{\Gamma}{1-\Gamma}
\eeq
This is an alternative derivation of the critical force for a simply
supported bar, or of the corresponding critical tensile strain $\eps_{Rt}^{(ss)}$ given by (\ref{eq:epsrtss}). In addition to the bifurcation condition,
we obtain an analytical description of the post-critical relation
between the (normalized) applied force and the finite rotation,
\beq 
 \frac{|P|}{EA} = \frac{\Gamma}{(1-\Gamma)\cos\varphi}
\eeq 
The total elongation of the bar can also be expressed analytically
in the form
\beq 
u = L\sqrt{(1+\eps)^2+\gamma^2} =  L\sqrt{\left(1+\frac{|P|}{EA}\cos\varphi\right)^2+\left(\frac{|P|}{GA_s}\sin\varphi\right)^2}
=\frac{|P|\,L}{GA_s}
\eeq

\section{Summary and conclusions}

We have proposed and tested a two-dimensional beam element that can accommodate arbitrary displacements and rotations in a geometrically exact manner. The sectional equations that describe the link between internal forces and deformation variables defined
at the level of an infinitesimal beam segment are considered as
linear, but extensions to material nonlinearity could be incorporated. Beam segments can deform by axial stretching,
shear and bending. These three modes are characterized by appropriate deformation variables and the corresponding 
work-conjugate internal forces. We have considered two
formulations that differ in the choice of the \old{deformation} variables
\old{for} \new{that characterize} axial and shear deformation.

\new{The following particular features of our approach can be emphasized:
\begin{itemize}
\item 
For both formulations, the main advantage is that that the minimum number of global degrees of
freedom is dictated only by the topology of the frame structure, and the accuracy can be
increased by refining the computational grid at the element level while the number of global
degrees of freedom is kept fixed.
\item 
Efficiency is also increased by directly addressing the global centerline coordinates and
sectional inclination with respect to the global axes as the primary unknowns at the element
level, thereby avoiding transformations between local and global coordinates.
\item
The algorithmic treatment has been extended to account for distributed loads and moments.
\item
Our model is not affected by shear locking, unlike displacement-based finite element
formulations of the Reissner–Simo beam model that account for transverse shear
deformation.
\item
Numerical simulations of both models confirm the results of the theoretical stability analysis
and efficiently allow tracing the entire post-buckling response. The Reissner model can
lead to bifurcations from a straight state under both compression and tension, whereas the Ziegler model, which is
more suitable for sandwich-like sections and built-up columns, leads only to
compressive buckling. 
\item 
The differences between the two models that arise from the choice of
primary deformation variables and the reference system are reflected in
the variation of internal forces along the beam. 
\end{itemize}

Two main limitations that we aim to address in our future studies concern (i) the
constitutive law, which is currently assumed to be linear elastic but, in general, should account
for material nonlinearities to accurately describe the mechanical behaviour of materials such as
rubbers and elastomers, and (ii) the two-dimensional formulation, which should
be extended to 3D.
}

\old{
In the theoretical part of the paper, {\bf stability} of a straight
beam loaded by an axial force has been investigated. It has been
shown that the traditional Reissner model has stability properties
that correspond to the Haringx approach, originally intended for springs, while the Ziegler model coincides,
 in the incompressible limit, with the more traditional Engesser
approach. For Reissner/Haringx, the compressive critical stress tends to infinity when the slenderness tends to zero,
and stability can also be lost in tension. On the other hand,
the Ziegler/Engesser formulation gives instability only in
compression and the critical stress remains bounded. 
In fact, for very stocky columns with not too low sectional
shear stiffness (as compared to the sectional axial stiffness),
bifurcation from the uniform compression state does not occur,
provided that the axial compressibility is taken into account
(which is not done by the original Engesser approach).
The results of the theoretical stability analysis have been
confirmed by numerical simulations, which can also efficiently trace
the entire post-buckling response.
In particular, it has been shown that the behavior of the
Reissner beam under high uniaxial tension is tricky and 
stability of the fundamental solution with uniform stretch and
zero rotations can be lost in modes that are strongly dependent on the boundary condition and their characterization requires
careful analysis, since it cannot be reduced to a unified 
description in which the support type only affects the buckling length.
}

\section*{Acknowledgments}
M.~Jir\'{a}sek and M.~Hor\'{a}k are grateful for the support of the Czech Science Foundation (project No.\ 19-26143X) received in 2023.
In 2024, their work was co-funded by the European Union under the ROBOPROX project (reg.~no.~CZ.02.01.01/00/22 008/0004590) and by the Czech Ministry of Education, Youth and Sports (ERC CZ project No.~LL2310).
\\
A grant from the Ministry for University and Research (MUR) for PRIN 2022, project No. 2022P7PF8J$\_$002, LAttice STructures for Energy aBsorption: advanced numerical analysis and optimal design (LASTEB), is gratefully acknowledged by C.~Bonvissuto and E.~La Malfa Ribolla. 

\appendix

\section{Two formulations of shear-flexible beams}\label{appA}

\subsection{Reissner model}\label{appA1}
\label{sec:2.2}

\subsubsection{Choice of deformation variables and strain energy density}

The sectional equations in the form originally suggested by Reissner \cite{reissner1972} and adopted in many other studies for two-dimensional beams \cite{BATISTA2016153, irschik2009continuum, humer2013exact} as well as their extension to 3D \cite{ibrahimbegovic1995finite, crisfield1999objectivity} are obtained if the strain energy density 
\beq \label{ee6}
{\cal E}_{int} = \half EA\eps^2 + \half GA_s\gamma^2 + \half EI \kappa^2
\eeq 
is postulated as
a quadratic function of the curvature, $\kappa$, and of two additional strain measures, 
\bea\label{ee7}
\varepsilon &=& \lambda_s \cos{\chi}-1\\
\gamma &=& \lambda_s \sin{\chi}\label{ee8}
\eea
which characterize the normal and shear deformation modes. 
Here, $\lambda_s$ and $\chi$ should be understood as intermediate variables with a clear geometric meaning, from which the ``true''
strain measures are derived. 
In the small-strain limit, the present description reduces
to the classical Timoshenko-Ehrenfest beam theory\cite{Tim21,Eli20}. 
In the undeformed
state, we have $\lambda_s=1$ and $\chi=0$. If $\lambda_s-1\ll 1$ and $\chi\ll 1$, the strain measures can be approximated
by $\eps\approx\lambda_s-1=$ relative change of length of the centerline and $\gamma\approx\chi=$ small shear angle (understood as the effective value, since the actual shear strain would vary
across the section).
Therefore, factors $EA$, $GA_s$ and $EI$
are standard sectional stiffnesses related to stretching, shear and bending that would be used by the linear Timoshenko model.
Each of them is the product of an elastic modulus (Young's modulus $E$, or shear modulus $G$) and a sectional geometry characteristic
(sectional area $A$, shear effective sectional area\cite{Cowper1966} $A_s$, or
sectional moment of inertia $I$).

In the variational approach,
the potential energy of the whole beam is considered as a functional dependent on the primary independent
kinematic fields, i.e., on the centerline displacement
components and the sectional rotation. For the centerline
stretch, we can directly use relation (\ref{eqn_lamsnew}). 
To express the shear angle $\chi$, we get back to the triangle in Fig.~\ref{T1}b and write
\bea \label{136}
\sin(\varphi-\chi)&=&-\frac{z'_s}{\lambda_s}\\
\cos(\varphi-\chi)&=&\frac{x'_s}{\lambda_s}\label{137}
\eea
For substitution into (\ref{ee7})--(\ref{ee8}), we do not really need 
$\chi$ itself but rather its sine and cosine. Rewriting
the left-hand sides of (\ref{136})--(\ref{137}) in
terms of sine and cosine of $\varphi$ and $\chi$ and
solving the resulting set of equations with the sine and
cosine of $\chi$ considered as unknowns, we obtain
\bea \label{ee11}
\cos\chi &=& \frac{x'_s}{\lambda_s}\cos\varphi - \frac{z'_s}{\lambda_s}\sin\varphi\\
\sin\chi &=&  \frac{x'_s}{\lambda_s}\sin\varphi  +\frac{z'_s}{\lambda_s}\cos\varphi
\label{ee12}
\eea 
Using this in (\ref{ee7})--(\ref{ee8}), we finally arrive at the strain-displacement equations
\bea\label{ee13}
\varepsilon &=&  x'_s\cos\varphi - z'_s\sin\varphi-1
\\
\gamma &=& x'_s\sin\varphi  +z'_s\cos\varphi
\label{ee14}
\eea

\subsubsection{Variational derivation of differential equilibrium equations}\label{sec:2.2.2}

Let us now focus on the potential energy and its variation. 
The potential energy stored in the elastic deformation of the whole beam is given by
\beq 
E_{int} = \int_0^L {\cal E}_{int}\,\dif\xi =
\int_0^L \left( \half EA\eps^2 + \half GA_s\gamma^2 + \half EI \kappa^2\right)\,\dif\xi
\eeq
and its first variation
\bea\nonumber
\delta E_{int} &=& \int_0^L \left(\frac{\partial{\cal E}_{int}}{\partial\eps}\delta\eps+\frac{\partial{\cal E}_{int}}{\partial\gamma}\delta\gamma+\frac{\partial{\cal E}_{int}}{\partial\kappa}\delta\kappa\right)\,\dif\xi =\\
\nonumber
&=&
\int_0^L \left(EA\eps\,\delta\eps + GA_s\gamma\,\delta\gamma + EI \kappa\,\delta\kappa\right)\,\dif\xi = \\
&=& \int_0^L \left(N\,\delta\eps + Q\,\delta\gamma + M\,\delta\kappa\right)\,\dif\xi
\label{ee16}
\eea 
can be interpreted as the internal virtual work done by
the internal forces
\bea \label{ee17}
N &=& \frac{\partial{\cal E}_{int}}{\partial\eps} = EA\eps \\
Q &=& \frac{\partial{\cal E}_{int}}{\partial\gamma} = GA_s\gamma \\
M &=& \frac{\partial{\cal E}_{int}}{\partial\kappa} = EI\kappa \label{ee19}
\eea 
on the virtual generalized strains $\delta\eps$, $\delta\gamma$ and $\delta\kappa$. Based on (\ref{ee13})--(\ref{ee14}) and (\ref{ee5}),
the virtual strains can be related to the virtual displacements and rotation:
\bea \label{ee20}
\delta\eps &=& \cos\varphi\,\delta x_s' - \sin\varphi\,\delta z_s' - \lambda_s\sin\chi\,\delta\varphi
\\ \label{ee21}
\delta\gamma &=& \sin\varphi\,\delta x_s' + \cos\varphi\,\delta z_s' + \lambda_s\cos\chi\,\delta\varphi
\\
\delta \kappa &=&\delta\varphi'
\label{ee22}
\eea 
In (\ref{ee20})--(\ref{ee21}), the factors 
multiplying $\delta\varphi$ have been rewritten
in terms of auxiliary variables $\lambda_s$ and $\chi$
based on relations (\ref{136})--(\ref{137}), which will later facilitate
the interpretation of equilibrium equation (\ref{eq28y}).
Substitution of (\ref{ee20})--(\ref{ee22}) into (\ref{ee16}) followed by integration by parts leads to
\bea\nonumber
\delta E_{int} &=& 
\left[(N\cos\varphi+Q\sin\varphi)\,\delta x_s+(-N\sin\varphi+Q\cos\varphi)\,\delta z_s+ M\,\delta\varphi\right]_0^L + \\ 
\nonumber
&&+
\int_0^L \left((-N\,\cos\varphi-Q\,\sin\varphi)'\delta x_s+(N\,\sin\varphi-Q\,\cos\varphi)'\delta z_s\right)\dif\xi +
\\
&&
+\int_0^L\left(-N\lambda_s\sin\chi+Q\lambda_s\cos\chi-M'\right)\delta\varphi \,\dif\xi
 \label{eq22xy}
\eea

Suppose that the beam is loaded by distributed external forces
described by global components $p_x$ and $p_z$,
and by a distributed external moment, $m$, all taken per unit length of the centerline in the undeformed state. The corresponding contribution to the potential energy is given by
\beq 
E_{ext} = -\int_0^L \left(p_x x_s+p_z z_s+m\varphi\right)\,\dif\xi
\eeq 
and its variation is
\beq \label{ee25}
\delta E_{ext} = -\int_0^L \left(p_x \,\delta x_s+p_z \,\delta z_s+m\,\delta \varphi\right)\,\dif\xi
\eeq 
Stationarity requirements for the total potential energy $E=E_{int}+E_{ext}$ correspond to a vanishing variation. Since the
variations of centerline coordinates and sectional rotation are
independent, condition $\delta E=\delta E_{int}+\delta E_{ext}=0$
leads to equations
\bea \label{eq26y}
-(N\,\cos\varphi+Q\,\sin\varphi)' &=& p_x
\\
(N\,\sin\varphi-Q\,\cos\varphi)' &=& p_z
\label{eq27y}
\\
-N\,\lambda_s\sin\chi+Q\,\lambda_s\cos\chi-M' &=& m
\label{eq28y}
\eea 
which can be interpreted as the equilibrium equations.
Based on (\ref{ee11})--(\ref{ee12}),
equation (\ref{eq28y})
could be rewritten in terms of the primary kinematic fields and the internal forces as
\beq
-\left(x_s'\sin\varphi+z_s'\cos\varphi\right)N+\left(x_s'\cos\varphi-z_s'\sin\varphi\right)Q-M' =m
\label{eq28y2}
\eeq

\subsubsection{Static boundary conditions}

Using the boundary terms in (\ref{eq22xy}), we could also
derive static (natural) boundary conditions at the beam ends. In the present context, the objective is not to solve a separate beam problem, but to develop 
the description of beam elements that will be assembled
into a frame structure. The coordinates and inclination
of the beam ends will be determined by the joint displacements
and rotations, which represent the global degrees of freedom, i.e., the basic unknowns of
the discretized frame model. 
In a typical incremental-iterative structural analysis, the basic
tasks on the beam element level are to evaluate the end forces and moments generated by a given increment of joint displacements and
rotations, and to evaluate the corresponding tangent element stiffness matrix.
Therefore, when a beam element
is processed, the coordinates and inclinations
of the beam ends can be considered as given.
The corresponding variations $\delta x_s$, $\delta z_s$ and $\delta\varphi$ at the beam ends are then zero, and the
boundary terms in (\ref{eq22xy}) vanish.
Nevertheless, for easier physical interpretation of the internal forces it is useful to look at the static boundary conditions
that would arise if the end variations are considered
as arbitrary. 
For instance, if the right beam end is considered
as free and is loaded by a given horizontal force $X_{ba}$,
vertical force $Z_{ba}$ and moment $M_{ba}$, 
then the contribution of external forces to the potential energy should be enriched by the expression
$E_{ext,L}=-X_{ba}x_s(L)-Z_{ba}z_s(L)-M_{ba}\varphi(L)$,
and the stationarity conditions deduced from the terms
that contain variations $\delta x_s(L)$, $\delta z_s(L)$
and $\delta\varphi(L)$ read
\bea \label{ee26}
N(L)\cos\varphi(L)+Q(L)\sin\varphi(L) &=& X_{ba} \\
-N(L)\sin\varphi(L)+Q(L)\cos\varphi(L) &=& Z_{ba} \label{ee27}\\
M(L) &=& M_{ba}\label{ee28}
\eea
Since $\varphi(L)$ is the inclination of the right end section with respect to the vertical global axis, it is clear that the internal forces $N$ and $Q$ have their standard meaning of the normal force, perpendicular to the section, and the shear force, parallel to the section, and $M$ is the bending moment. This conclusion may look trivial, but in Section~\ref{sec:modmod} it will be shown that, for a model based on 
another set of deformation variables, the meaning
of the work-conjugate quantities may be different.

In a completely analogous way, boundary conditions 
\bea \label{ee33}
-N(0)\cos\varphi(0)-Q(0)\sin\varphi(0) &=& X_{ab} \\
N(0)\sin\varphi(0)-Q(0)\cos\varphi(0) &=& Z_{ab} \label{ee34} \\
-M(0) &=& M_{ab}\label{ee35}
\eea
valid
at the left end can be derived. They differ only by the signs of
the terms on the left-hand side. Of course, the terms on the right-hand side are the forces and the moment acting on the left end.

\subsubsection{Conversion to a set of first-order equations}
\label{sec:2.2.4}

The complete Reissner beam model consists of equilibrium equations
(\ref{eq26y})--(\ref{eq27y}) and (\ref{eq28y2}), sectional equations (\ref{ee17})--(\ref{ee19}), and kinematic (strain-displacement) equations
(\ref{ee13})--(\ref{ee14}) and (\ref{ee5}), supplemented by
appropriate boundary conditions. 
Substituting the kinematic and sectional equations into the 
equilibrium equations, we could construct a set of three second-order
differential equations for unknown functions $x_s$, $z_s$ and $\varphi$. Instead of that, let us first integrate the equilibrium
equations analytically. Equations (\ref{eq26y})--(\ref{eq27y})
with boundary conditions (\ref{ee33})--(\ref{ee34}) have a closed-form solution
\bea \label{ee36}
N &=& -\left(X_{ab}+P_x\right)\cos\varphi + \left(Z_{ab}+P_z\right)\sin\varphi \\
Q &=& -\left(X_{ab}+P_x\right)\sin\varphi - \left(Z_{ab}+P_z\right)\cos\varphi \label{ee37}
\eea 
where
\bea 
P_x(\xi) &=& \int_0^\xi p_x(s)\,\mathrm{d}s
\\
P_z(\xi) &=& \int_0^\xi p_z(s)\,\mathrm{d}s
\eea
are functions that can be deduced from the prescribed distributed
load intensities $p_x$ and $p_z$. Recall that $X_{ab}$ and $Z_{ab}$
are constants that represent the global components of the force between the left end of the beam and the joint to which it is attached.
The values of functions $P_x$ and $P_z$ at $\xi$ correspond to the resultants of the
distributed loads between the left end and the section
at distance $\xi$ from the left end in the undeformed state.
Clearly, equations (\ref{ee36})--(\ref{ee37}) are force equilibrium equations
for the part of the beam between the left end and 
the section characterized by $\xi$, written for the directions
perpendicular to this section and parallel to this section.

An analogous formula can be constructed for the bending moment. Substituting (\ref{ee36})--(\ref{ee37})
into (\ref{eq28y2}) and rearranging the terms, we get
\beq 
M' = -m +(X_{ab}+P_x)z_s'-(Z_{ab}+P_z)x_s'
\eeq 
and, after integration,
\beq \label{ee41}
M = -M_{ab} + X_{ab}(z_s-z_a)-Z_{ab}(x_s-x_a) + M_p
\eeq 
where
\beq \label{e41}
M_p(\xi) = -\int_0^\xi m(s)\,{\rm d}s + \int_0^\xi P_x(s)\, z_s'(s)\, \textrm{d}s  -  \int_0^\xi P_z(s)x_s'(s)\, \textrm{d}s 
\eeq 
and symbols $x_a$ and $z_a$ denote the current coordinates of the centroid of the left end section, i.e., $x_a=x_s(0)$ and $z_a=z_s(0)$. We have also exploited the static boundary condition (\ref{ee35}).


The foregoing derivation may look somewhat complicated, but the result has a direct physical meaning. The first three terms on the right-hand side of (\ref{ee41}) represent
the contributions of a concentrated moment $M_{ab}$ and forces
$X_{ab}$ and $Z_{ab}$ applied at the left beam end to the 
bending moment at cross section $x$, taking into account the
deformed state of the beam, and the last term, $M_p$, 
represents the contribution of the distributed loads to the same
bending moment. 

In summary, equations (\ref{ee36})--(\ref{ee37}) and (\ref{ee41})
are equilibrium equations in their integrated form, as opposed
to the differential equations (\ref{eq26y})--(\ref{eq27y}) and (\ref{eq28y2}), which describe equilibrium
of an infinitesimal beam segment. 
The internal forces are easily transformed into deformation
variables using the inverted form of sectional equations (\ref{ee17})--(\ref{ee19}):
\bea 
\eps &=& \frac{N}{EA} \label{ee46}\\
\gamma &=& \frac{Q}{GA_s} \label{ee47}\\
\kappa &=& \frac{M}{EI} \label{ee48}
\eea 
Kinematic equations (\ref{ee13})--(\ref{ee14}) link deformation variables $\eps$ and $\gamma$ to the derivatives
of the displacement functions, and their inversion leads to
\bea \label{e46}
x_s' &=& (1+\eps)\cos\varphi+\gamma\sin\varphi\\
z_s' &=& \gamma\cos\varphi-(1+\eps)\sin\varphi
\label{e47}
\eea 
Substituting from (\ref{ee46})--(\ref{ee47}) and (\ref{ee36})--(\ref{ee37}) into (\ref{e46})--(\ref{e47}), we end up with two first-order differential  equations
\bea 
\nonumber
x_s' &=& \left(1+\frac{-(X_{ab}+P_x)\cos\varphi + (Z_{ab}+ P_z)\sin\varphi}{EA}\right)\cos\varphi\\
\label{ee51}
&&- \frac{\left(X_{ab}+P_x\right)\sin\varphi + \left(Z_{ab}+P_z\right)\cos\varphi}{GA_s}\sin\varphi \\
\nonumber
z_s' &=& -\frac{\left(X_{ab}+P_x\right)\sin\varphi + \left(Z_{ab}+P_z\right)\cos\varphi}{GA_s}\cos\varphi\\
&&-\left(1+\frac{-(X_{ab}+P_x)\cos\varphi + (Z_{ab}+ P_z)\sin\varphi}{EA}\right)\sin\varphi
\label{ee52}
\eea 
that contain three primary unknown functions, $x_s$, $z_s$ and $\varphi$.
The third equation of a similar kind,
\beq  \label{ee53}
\varphi' = \frac{-M_{ab} +X_{ab}(z_s-z_a)- Z_{ab}(x_s-x_a)
+  M_p }{EI}
\eeq
is obtained simply by combining (\ref{ee5}) with (\ref{ee48}) and (\ref{ee41}).
Equations (\ref{ee51})--(\ref{ee53}) form the theoretical basis
of the numerical procedure developed in this paper.
These first-order differential equations will be approximated using a finite difference scheme; see Section~\ref{sec:3.2.1}.

Similar, yet simpler equations were derived in \cite{JirMalHor21} for the Kirchhoff beam in a somewhat different setting---the deformed shape of the centerline was described by the displacements
(instead of current coordinates) with respect to a local coordinate system aligned with the beam. The present equations reduce to the 
formulation in  \cite{JirMalHor21} if the current coordinates
are expressed as $x_s=\xi+u_s$ and $z_s=w_s$, and the left-end coordinates $x_a$ and $z_a$ are set to $u_a$ and $w_a$. Moreover, fractions with $GA_s$ in the denominator are deleted,
because the shear stiffness is considered as infinite (no shear distortion is taken into account), and functions $P_x$, $P_z$ and $M_p$ that reflect the loads distributed along the beam are set to zero. After these modifications,
equations (\ref{ee51})--(\ref{ee53}) reduce to
\bea 
u_s' &=& \left(1+\frac{-X_{ab}\cos\varphi + Z_{ab}\sin\varphi}{EA}\right)\cos\varphi
-1\\ 
w_s' &=& -\left(1+\frac{-X_{ab}\cos\varphi + Z_{ab}\sin\varphi}{EA}\right)\sin\varphi
\\
\varphi' &=& \frac{-M_{ab} +X_{ab}(w_s-w_a)- Z_{ab}(\xi+u_s-u_a)}{EI}
\eea 
This exactly corresponds to equations (37)--(39) from \cite{JirMalHor21},
with only a formal difference---coordinate $\xi$ was replaced by $x$, because the coordinate system $(x,y,z)$ used  in \cite{JirMalHor21} was local, aligned with the beam.

\subsection{Ziegler model}\label{appA2}
\label{sec:modmod}
\subsubsection{Choice of deformation variables}

The standard Reissner model, presented in Section~\ref{sec:2.2}, 
characterizes the normal mode by the deformation
variable $\eps$ defined in (\ref{ee7}) and the shear mode by the
deformation variable $\gamma$ defined in (\ref{ee8}). The
work-conjugate quantities are then the normal force $N$ perpendicular to the deformed section, and the shear force $Q$
parallel to the deformed section.
This choice conceptually corresponds to the work of Haringx \cite{Haringx1942}, who studied the stability of helical springs and modified the previously used stability approach of Engesser \cite{Engesser1891}. 
Due to the particular definition of $\eps$, the axial stiffness $EA$ is activated by the change of distance between the neighboring sections. If the sections slide 
while keeping their distance constant, only the shear stiffness $GA_s$ is activated. As will be shown in Section~\ref{sec:5.1.2}, such a model predicts a certain finite critical load even in uniaxial tension, because the sliding mechanism can lead to an extension of the beam length measured along the centerline
even at vanishing ``normal strain'' $\eps$.  On the other hand,
in compression, the critical stress tends to infinity as the slenderness approaches zero. 

Alternatively, the variable characterizing the normal mode
can be derived from
the change of length of the centerline, which can occur as
a second-order effect also during parallel sliding of the sections.
This point of view could be more appropriate for sandwich-like sections
with top and bottom hard layers connected by a soft core, or 
for built-up columns modeled as
one single rod with equivalent sectional properties.
Such an approach corresponds to the original Engesser model for stability of a shear-deformable column\cite{Engesser1891}, which does not give any bifurcation in tension while in compression it yields a finite critical stress even in the limit of zero slenderness. 
In this extreme case, the buckling mode corresponds to a compressed layer in which individual fibers start rotating at a constant length, and the second-order work of the normal stress provides enough energy needed to induce shear deformation. 
In his simple stability analysis, Engesser\cite{Engesser1891}
took into account the effect of shear distortion, but he still
treated the beam axis as inextensible and incompressible.
The effect of axial deformation was added by Ziegler\cite{Zie82},
who demonstrated that it is of the same order of magnitude as the 
effect of shear. Ziegler's analysis was focused on stability
of an initially straight column, and his basic equations contain
some approximations based on the assumption of small rotations
and small shear angles. Therefore, his approach was not geometrically exact, but his sectional equations can be 
combined with kinematic equations
(\ref{ee13})--(\ref{ee14}) and (\ref{ee5}), leading to a model that can be considered as an
alternative formulation of a geometrically exact beam.
Adopting the terminology used in \cite{KimMehAtt20}, we will refer to this formulation as the Ziegler model,
even though it does not precisely correspond to the 
equations originally used in \cite{Zie82}.

The foregoing considerations motivate a modified choice of deformation
variables. Instead of $\eps$ defined in (\ref{ee7}) and related to the change in distance
between sections, let us use,
as the deformation variable characterizing the normal mode, the axial strain
\beq 
\teps = \lambda_s-1 = \sqrt{x_s'^2+z_s'^2} - 1
\eeq 
This variable corresponds to the Biot (engineering) strain \old{,
already mentioned in (\ref{ee7biot}),} and represents the relative change of length of the centerline.
For the shear deformation variable characterizing the shear mode, a natural choice is the shear angle $\chi$.
Based on (\ref{136})--(\ref{137}), the shear angle can be expressed as
\beq
\chi = \varphi + \arctan\frac{z_s'}{x_s'}
\eeq 
The flexural deformation is still described by the curvature $\kappa$
defined in (\ref{ee5}).

\subsubsection{Variational derivation of differential equilibrium equations and static boundary conditions}
\label{sec232}

The elastic energy density can be again considered in the simple quadratic and decoupled form
\beq \label{eq94}
{\cal E}_{int} = \half EA\teps^2 + \half GA_s\chi^2 + \half EI \kappa^2
\eeq 
which leads to linear sectional equations
\bea \label{ee17x}
\tN &=& \frac{\partial{\cal E}_{int}}{\partial\teps} = EA\teps \\
\label{ee18x}
\tQ &=& \frac{\partial{\cal E}_{int}}{\partial\chi} = GA_s\chi \\
M &=& \frac{\partial{\cal E}_{int}}{\partial\kappa} = EI\kappa \label{ee19x}
\eea 
The physical meaning of internal forces $\tN$ and $\tQ$ is different
from the previously used normal force $N$ and shear force $Q$,
and this is why they are denoted here by modified symbols with a tilde. This alternative form of the sectional equations has been debated and compared to the Reissner formulation by various authors\cite{KimMehAtt20,  genovese2019shear, 10.1115/1.3197142}.
It will be shown that $\tN$ is
the force tangent to the centerline and $\tQ$ is the force perpendicular
 to the centerline multiplied by $\lambda_s$.
On the other hand, the bending moment $M$ keeps the same meaning as before.

For the Ziegler model, a variational procedure
analogous to Section~\ref{sec:2.2.2}
leads to the variation of stored elastic energy in the form
\beq
\delta E_{int} = \int_0^L \left(\tN\,\delta\teps + \tQ\,\delta\chi + M\,\delta\kappa\right)\,\dif\xi
\label{ee16x}
\eeq
After the substitution of the variations
\bea 
\delta\teps &=&  \frac{x_s'\delta x_s' + z_s'\delta z_s'}{\lambda_s}
= \cos(\varphi-\chi)\,\delta x_s' -\sin(\varphi-\chi)\,\delta z_s'
\\
\delta\chi &=& \delta\varphi + \frac{x_s'\delta z_s' - z_s'\delta x_s'}{\lambda_s^2} =\delta\varphi + \frac{\cos(\varphi-\chi)\,\delta z_s' +\sin(\varphi-\chi)\,\delta x_s'}{\lambda_s} \\
\delta\kappa &=& \delta\varphi'
\eea 
and integration by parts, the expression for the variation is converted into
\bea\nonumber
\delta E_{int}&=& 
\left[
\left(\tN\cos(\varphi-\chi)+\frac{\tQ}{\lambda_s}\sin(\varphi-\chi)\right)\delta x_s + \right. 
\\ \nonumber
&&+\left.\left(-\tN\sin(\varphi-\chi)+\frac{\tQ}{\lambda_s}\cos(\varphi-\chi)\right)\delta z_s+M\delta\varphi\right]_0^L+
\\ \nonumber
&&+\int_0^L 
\left(-\tN\cos(\varphi-\chi)-\frac{\tQ}{\lambda_s}\sin(\varphi-\chi)\right)'\delta x_s\,\dif\xi + \\ \nonumber
&&+ \int_0^L\left(\left(\tN\sin(\varphi-\chi)-\frac{\tQ}{\lambda_s}\cos(\varphi-\chi)\right)'\delta z_s +\left(\tQ-M'\right)\delta\varphi \right)
\,\dif\xi \\
\label{ee16y}
\eea
The variation of the energy of external forces $\delta E_{ext}$ is still given by (\ref{ee25}), with possible enhancements by boundary terms,
and the resulting
stationarity
conditions derived from $\delta E_{int}+\delta E_{ext}=0$
read
\bea \label{ee80}
-\left(\tN\cos(\varphi-\chi)+\frac{\tQ}{\lambda_s}\sin(\varphi-\chi)\right)'&=& p_x
\\ \label{ee81}
\left(\tN\sin(\varphi-\chi)-\frac{\tQ}{\lambda_s}\cos(\varphi-\chi)\right)'&=& p_z
\\ \label{ee82}
\tQ-M' &=& m
\eea 
A comparison with equilibrium equations (\ref{eq26y})--(\ref{eq28y}) previously derived for the Reissner model shows that, in the force equilibrium equations
(\ref{ee80})--(\ref{ee81}), the internal forces $\tN$ and $\tQ$
are now rotated by $\varphi-\chi$, which is the rotation of the tangent
to the deformed centerline, while for the Reissner model the analogous
forces $N$ and $Q$ were rotated by the sectional inclination angle, $\varphi$. Furthermore, force $\tQ$ work-conjugate
with the shear angle $\chi$ is in the force equilibrium equations
(\ref{ee80})--(\ref{ee81}) divided by the centerline stretch, $\lambda_s$. To work with quantities that have a direct physical
meaning, it is preferable to introduce the ``true'' force
acting in the direction perpendicular to the deformed centerline,
defined as $Q^*=\tQ/\lambda_s$. 
The equilibrium equations then have the form
\bea \label{em1}
-\left(\tN\cos(\varphi-\chi)+Q^*\sin(\varphi-\chi)\right)'&=& p_x
\\ \label{em2}
\left(\tN\sin(\varphi-\chi)-Q^*\cos(\varphi-\chi)\right)'&=& p_z
\\ \label{em3}
\lambda_s Q^*-M' &=& m
\eea 
and the static boundary conditions, derived using the boundary terms
in (\ref{ee16y}) combined with appropriate boundary terms in the variation of the external potential energy, read
\bea \label{ee86}
-\tN(0)\cos(\varphi(0)-\chi(0))-Q^*(0)\sin(\varphi(0)-\chi(0)) &=& X_{ab} \\
\tN(0)\sin(\varphi(0)-\chi(0))-Q^*(0)\cos(\varphi(0)-\chi(0)) &=& Z_{ab} \\ \label{ee88}
-M(0) &=& M_{ab} \\
\tN(L)\cos(\varphi(L)-\chi(L))+Q^*(L)\sin(\varphi(L)-\chi(L)) &=& X_{ba} \\
-\tN(L)\sin(\varphi(L)-\chi(L))+Q^*(L)\cos(\varphi(L)-\chi(L)) &=& Z_{ba} \\
M(L) &=& M_{ba}
\eea 

Notice that, if $Q^*$ is used as one of the internal forces, 
sectional equation (\ref{ee18x}) needs to be rewritten as
\beq\label{eqQstar}
Q^* = GA_s\frac{\chi}{\lambda_s} 
\eeq

\subsubsection{Conversion to a set of first-order equations}

The subsequent derivation proceeds in a similar way as in Section~\ref{sec:2.2.4} for the Reissner model, but it is somewhat more complicated. First, equilibrium
equations in the differential form (\ref{em1})--(\ref{em3}) are integrated and
boundary conditions (\ref{ee86})--(\ref{ee88}) are exploited. The resulting 
expressions for the internal forces are
\bea \label{ee36z}
\tN &=& -\left(X_{ab}+P_x\right)\cos(\varphi-\chi) + \left(Z_{ab}+P_z\right)\sin(\varphi-\chi) \\
Q^* &=& -\left(X_{ab}+P_x\right)\sin(\varphi-\chi) - \left(Z_{ab}+P_z\right)\cos(\varphi-\chi) \label{ee37z}
\\
\label{ee41z}
M &=& -M_{ab} + X_{ab}(z_s-z_a)-Z_{ab}(x_s-x_a) + M_p
\eea
Equation (\ref{ee41z}) for the bending moment is exactly the same
as equation (\ref{ee41}) derived for the Reissner model.
However, equations (\ref{ee36z})--(\ref{ee37z}) for 
internal force components $\tilde N$ and $Q^*$ are
 different from   equations (\ref{ee36})--(\ref{ee37}) 
valid for the Reissner model, since  the expressions on their right-hand sides contain
not only the external forces and the primary kinematic variables $x_s$, $z_s$
and $\varphi$, but also the shear angle $\chi$, which is one of the deformation variables.
Therefore, when 
these expressions for internal forces are substituted into the inverted sectional equations 
\bea \label{ee103}
\teps &=& \frac{\tN}{EA} 
\\ \label{ee104}
\chi &=& \frac{\tQ}{GA_s} = \frac{\lambda_s Q^*}{GA_s}= \frac{(1+\teps)\, Q^*}{GA_s} 
\\ \label{ee105}
\kappa &=& \frac{M}{EI}
\eea 
we find the deformation variables not only on the left-hand side,
but two of them, namely $\teps$ and $\chi$, affect also the right-hand side.
 In particular, 
the relation obtained by substituting (\ref{ee37z}) into (\ref{ee104})
is not an explicit formula for the evaluation
of $\chi$ from other variables, but an implicit definition of $\chi$
as the solution of the nonlinear algebraic equation
\beq 
\chi = 
-\frac{1+\teps}{GA_s}\left( \left(X_{ab}+P_x\right)\sin(\varphi-\chi) + \left(Z_{ab}+P_z\right)\cos(\varphi-\chi)\right)
\eeq 
in which $\teps$ is given by (\ref{ee103}) combined with (\ref{ee36z}) and thus also depends
on $\chi$. This equation does not have a closed-form solution,
but in practical computations can be treated numerically.
Formally, we can define function $f_\chi(\varphi,P_1,P_2)$
such that its value $\chi$ for given arguments $\varphi$, $P_1$ and $P_2$
is the solution of the equation
\beq \label{eq92}
GA_s\chi +\left(1-\frac{P_1\cos(\varphi-\chi) + P_2\sin(\varphi-\chi)}{EA}\right)\left( P_1\sin(\varphi-\chi) + P_2\cos(\varphi-\chi)\right)=0
\eeq
that satisfies the constraint $-\pi/2<\chi<\pi/2$. 
For an arbitrary $\varphi$, we have $f_\chi(\varphi,0,0)=0$. This can be a starting point 
of an iterative scheme computing an approximation of $f_\chi(\varphi,P_1,P_2)$
for general arguments $P_1=X_{ab}+P_x$ and $P_2=Z_{ab}+P_z$, e.g., by the Newton method. 

The last step of the derivation exploits kinematic equations (\ref{136})--(\ref{137}), which provide  expressions for the coordinate derivatives 
\bea \label{171}
x_s' &=& \lambda_s \cos(\varphi-\chi)  = (1+\teps) \cos(\varphi-\chi) \\
z_s' &=& -\lambda_s \sin(\varphi-\chi) =-(1+\teps) \sin(\varphi-\chi)
\label{172}
\eea 
Replacing $\teps$ by the expression obtained from
(\ref{ee103}) and (\ref{ee36z}), and replacing $\chi$ by
$f_\chi(\varphi,X_{ab}+P_x,Z_{ab}+P_z)$,
 we construct two first-order differential equations
 \bea \nonumber
x_s' &=& 
\scalebox{0.8}{$\left(1-\dfrac{\left(X_{ab}+P_x\right)c_\varphi(\varphi,X_{ab}+P_x,Z_{ab}+P_z) -\left(Z_{ab}+P_z\right)s_\varphi(\varphi,X_{ab}+P_x,Z_{ab}+P_z)}{EA}\right) c_\varphi(\varphi,X_{ab}+P_x,Z_{ab}+P_z) $}
\\
\label{ee110}
\\ \nonumber
z_s' &=&
\scalebox{0.8}{$\left(-1+\dfrac{\left(X_{ab}+P_x\right)c_\varphi(\varphi,X_{ab}+P_x,Z_{ab}+P_z) - \left(Z_{ab}+P_z\right)s_\varphi(\varphi,X_{ab}+P_x,Z_{ab}+P_z)}{EA}\right) s_\varphi(\varphi,X_{ab}+P_x,Z_{ab}+P_z)$}
\\
\label{ee111}
\eea 
with unknown functions $x_s$, $z_s$ and $\varphi$.
To simplify the notation, we have introduced auxiliary functions 
\bea 
c_\varphi(\varphi,P_1,P_2) &=& \cos\left(\varphi-f_\chi(\varphi,P_1,P_2)\right) \\
s_\varphi(\varphi,P_1,P_2) &=& \sin\left(\varphi-f_\chi(\varphi,P_1,P_2)\right)
\eea 
The third differential equation, 
\beq \label{ee137}
\varphi' =  \frac{-M_{ab} + X_{ab}(z_s-z_a)-Z_{ab}(x_s-x_a) + M_p}{EI}
\eeq
is then obtained from the kinematic equation $\varphi'=\kappa$ combined with (\ref{ee105}) and (\ref{ee41z}).
This equation is exactly the same as equation (\ref{ee53}) for the Reissner model, but equations
(\ref{ee110})--(\ref{ee111}) are different. 
They look more complicated than the corresponding
equations (\ref{ee51})--(\ref{ee52}) obtained for the Reissner model,
but in Section~\ref{sec:3.4} it will be shown that the numerical
implementation is not much more demanding.


\section{Analytical solution for a cantilever subjected to uniformly distributed external moment}\label{appB}

This derivation leads to a reference solution for the numerical simulations in Section~\ref{sec:example3}.  An initially straight cantilever
of length $L$ is loaded by an externally applied distributed
moment of uniform intensity $m$.
If we consider the
left end as fixed and the right end as free,
the fundamental
equations (\ref{ee51})--(\ref{ee53}), with 
$M_{ab}=-mL$ and $M_p(\xi)=-m\xi$ evaluated from (\ref{e41}), reduce to
\bea 
x_s' &=& \cos\varphi \label{e200}\\
z_s' &=& -\sin\varphi \label{e201}\\
\varphi' &=& \frac{m}{EI}(L-\xi) \label{e202}
\eea 
Equation (\ref{e202}) with initial condition $\varphi(0)=0$ leads to
\beq 
\varphi(\xi) =  \frac{m}{2EI}(2L\xi-\xi^2)
\eeq 
Substituting this into (\ref{e200})--(\ref{e201}) and integrating, we obtain
\bea \nonumber
x_s(\xi) &=& \int_0^\xi \cos\left(\frac{m}{2EI}(2Ls-s^2)\right)\,{\rm d}s =\\ \nonumber
&=&\frac{L}{\mu}\left(\cos\frac{\pi\mu^2}{2}\left(C(\mu)-C\left(\mu-\frac{\mu\xi}{L}\right)\right)+\sin\frac{\pi\mu^2}{2}\left(S(\mu)-S\left(\mu-\frac{\mu\xi}{L}\right)\right)\right)\\
\\ \nonumber
z_s(\xi) &=&  -\int_0^\xi \sin\left(\frac{m}{2EI}(2Ls-s^2)\right)\,{\rm d}s = \\ \nonumber
&=& \frac{L}{\mu}\left(\sin\frac{\pi\mu^2}{2}\left(C(\mu)-C\left(\mu-\frac{\mu\xi}{L}\right)\right)-\cos\frac{\pi\mu^2}{2}\left(S(\mu)-S\left(\mu-\frac{\mu\xi}{L}\right)\right)\right)
\\
\eea
where
\beq 
\mu = L\sqrt{\frac{m}{\pi EI}}
\eeq 
is a dimensionless parameter introduced for convenience, and $C$ and $S$ are Fresnel integrals defined as
\bea 
C(x) = \int_0^x \cos s^2\,{\rm d}s, \hskip 10mm
S(x) = \int_0^x \sin s^2\,{\rm d}s 
\eea 
In particular, the coordinates of the cantilever tip in the deformed state are given by
\bea \label{eq205}
x_b &=& x_s(L) = \frac{L}{\mu}\left(\cos\frac{\pi\mu^2}{2}C(\mu)+\sin\frac{\pi\mu^2}{2}S(\mu)\right)\\ \label{eq206}
z_b &=& z_s(L) =  \frac{L}{\mu}\left(\sin\frac{\pi\mu^2}{2}C(\mu)-\cos\frac{\pi\mu^2}{2}S(\mu)\right)
\eea


\section{Bifurcation analysis}\label{appC}
\label{sec:5}\label{buckling}

This appendix investigates how the choice of a specific formulation of the sectional equations
affects bifurcations from a straight state of an axially loaded column. 
Analytical expressions for the critical load and buckling modes are derived for the Reissner and Ziegler models, and it is shown that the Reissner formulation can
lead to bifurcations not only in compression but also in tension.

\subsection{General setting}

Let us first recall that the adopted flexibility-based approach
to the description of a geometrically exact beam leads
to a set of three first-order differential equations,
which have the form (\ref{ee51})--(\ref{ee53}) if the Reissner sectional equations are selected, or (\ref{ee110})--(\ref{ee111}) and (\ref{ee137}) if the Ziegler sectional equations are used.
If an initially straight prismatic column of length $L$
is loaded axially by a force $P$
oriented such that $P>0$ leads to compression, and if no other
loading is present, there exists a fundamental solution 
that corresponds to uniform compression, without any bending
or shear. This solution is characterized by kinematic variables 
\bea \label{ee286}
x_{s0}(\xi) &=& \left(1-\frac{P}{EA}\right)\xi \\
z_{s0}(\xi) &=& 0 \\
\varphi_0(\xi) &=& 0\label{ee288}
\eea 
which describe the centerline coordinates in the deformed state and the sectional inclination.
Here we assume that the origin of the coordinate system is placed at the centroid of one end section, and that the $x$-axis 
passes through the centroid of the other end section and coincides with the local $\xi$-axis.

Functions (\ref{ee286})--(\ref{ee288}) satisfy the governing
differential equations and also the appropriate 
boundary conditions that correspond to various combinations
of the end supports. For the axial direction, we always
consider a blocked displacement at $\xi=0$ and a free displacement
at $\xi=L$, where the external force $P$ is applied.
This means that the kinematic boundary condition
\beq\label{eqs211}
x_s(0)=0
\eeq 
is imposed at the fixed end, while the free end is characterized
by the static boundary condition $X_{ba}=-P$. Using 
global equilibrium, the latter condition can be transformed into 
\beq 
X_{ab}=P
\eeq 
which is more handy for our purpose, since the governing differential
equations contain the left-end forces $X_{ab}$, $Z_{ab}$ and $M_{ab}$.

\begin{figure}
    \centering
    \begin{tabular}{ccc}
       (a) & (b)  & (c)  \\[3mm]
       \includegraphics[width=0.2\linewidth]{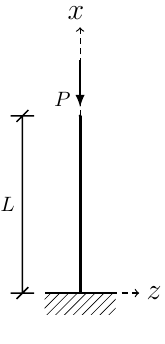}
         & 
           \includegraphics[width=0.2\linewidth]{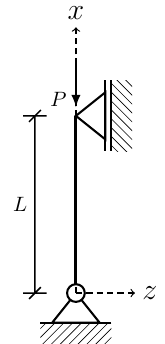}
           &
           \includegraphics[width=0.2\linewidth]{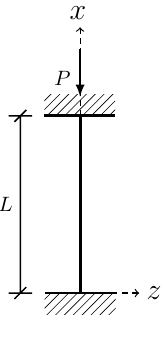}
    \end{tabular}
    \caption{Types of supports considered in the analysis: (a) a column clamped at one end, (b) a simply supported column, (c)~a column clamped at both ends}
    \label{fig:stab1}
\end{figure}

For constraints imposed on the lateral displacement and sectional rotation,
let us consider the following three options (Fig.~\ref{fig:stab1}):
\begin{enumerate}
    \item A {\bf column
clamped at one end}, with kinematic boundary conditions
\bea \label{eqs213}
z_s(0) &=& 0\\
\varphi(0) &=& 0
\eea 
and static boundary conditions $Z_{ba}=0$ and $M_{ba}=0$, which can be transformed via global equilibrium equations into
\bea 
Z_{ab} &=& 0\\
M_{ab} &=& P\,z_s(L) \label{eqs216}
\eea
\item 
A {\bf simply supported column}, with kinematic boundary conditions
\bea \label{eqs217}
z_s(0) &=& 0\\
z_s(L) &=& 0\label{eqs218}
\eea 
and static boundary conditions $M_{ab}=0$ and $M_{ba}=0$, which can be transformed via global equilibrium equations into
\bea \label{eqs219}
Z_{ab} &=& 0\\
M_{ab} &=& 0\label{eqs220}
\eea
\item A {\bf column clamped at both ends}, with kinematic boundary conditions
\bea \label{eqs221}
z_s(0) &=& 0\\
\label{eqs222}
\varphi(0) &=& 0\\
\label{eqs223}
z_s(L) &=& 0\\
\varphi(L) &=& 0\label{eqs224}
\eea 
\end{enumerate}

Let us now examine whether, for a given force $P$, there exists
a slightly perturbed solution that differs from the fundamental one
but also satisfies the governing equations and boundary conditions.
The perturbations are denoted by $\delta$ preceding the original symbol, and they are considered  as infinitely small, which makes
it possible to linearize the governing equations around the
fundamental solution. 
The fully general form of the linearized equations would be
quite lengthy, but it is greatly simplified when the
fundamental solution is substituted. For instance,
since $\varphi(\xi)=0$, all terms with $\sin\varphi$ vanish
and $\cos\varphi$ can be replaced by~1.
For the fundamental solution, $X_{ab}=P$, $Z_{ab}=0$ and
$M_{ab}=0$, but the perturbations $\delta Z_{ab}$ and $\delta M_{ab}$
can be nonzero while $\delta X_{ab}=0$. No distributed loads are considered,
and so functions $P_x$, $P_z$ and $M_p$ vanish and their perturbations as well. 

It should be emphasized that the mathematical analysis to be presented next is focused on the possible bifurcation from the fundamental solution. The lowest load at which such a bifurcation is detected will be considered as the critical
load. Before the critical load is attained, the solution
is unique and it could be shown that it is stable in the sense that an arbitrary infinitesimal perturbation would lead
to a (second-order) increase of the potential energy.
When the critical load is exceeded, the fundamental solution 
typically loses stability and the physical system follows
a bifurcated branch of the equilibrium diagram.

\subsection{Critical load -- Reissner model}
\label{sec:5.1.2}

\subsubsection{General analysis}\label{sec:reis-stab}

For the Reissner model, linearization of the governing equations (\ref{ee51})--(\ref{ee53}) around
the fundamental solution (\ref{ee286})--(\ref{ee288}) yields
\bea \label{ee291}
\delta x_s' &=& 0
\\ \label{ee292}
\delta z_s' &=& -\frac{1}{GA_s}\left(P\,\delta\varphi+\delta Z_{ab}\right) - \left(1-\frac{P}{EA}\right)\delta\varphi
\\ \label{ee293}
\delta\varphi' &=& \frac{1}{EI}\left(-\delta M_{ab}+P\,\delta z_s-\delta Z_{ab}\left(1-\frac{P}{EA}\right)\xi\right)
\eea 
Equation (\ref{ee291}) with boundary condition (\ref{eqs211}) leads to
$\delta x_s(\xi)=0$. Equations (\ref{ee292})--(\ref{ee293}) with unknown functions $\delta z_s$ and $\delta \varphi$ can be combined 
into a single second-order differential equation with one unknown
function. Differentiating (\ref{ee293}) and replacing $\delta z_s'$ by the right-hand side of (\ref{ee292}),
we obtain, after a rearrangement,
\beq \label{eqr1}
EI\,\delta\varphi'' + (1+\beta P)P\,\delta\varphi = - (1+\beta P)\,\delta Z_{ab} 
\eeq 
where
\beq \label{eq:beta}
\beta = \frac{1}{GA_s} - \frac{1}{EA}
\eeq 
is a newly introduced constant, which can be evaluated from the
given sectional stiffnesses.

In cases when $(1+\beta P)P>0$, the general solution
of equation (\ref{eqr1}) has the form
\beq \label{eqgensol}
\delta\varphi(\xi) = -\frac{\delta Z_{ab}}{P} + C_1\cos k\xi + C_2 \sin k\xi
\eeq 
where
\beq \label{eq:k}
k=\sqrt{\frac{(1+\beta P)P}{EI}}
\eeq 
and $C_1$ and $C_2$ are arbitrary integration constants.
The subsequent analysis depends on the chosen type of supports.

\begin{enumerate}
\item
For a {\bf column with one clamped end}, the solution of (\ref{ee292})--(\ref{ee293}) must
satisfy conditions $\delta z_s(0)=0$, $\delta\varphi(0)=0$,
$\delta Z_{ab}=0$ and $\delta M_{ab}=P\,\delta z_s(L)$,
which follow directly from (\ref{eqs213})--(\ref{eqs216}).
Equation (\ref{ee293}) evaluated at $\xi=L$ is then exploited to construct a second condition written
in terms of the primary unknown $\delta\varphi$, namely
$\delta\varphi'(L)=0$, which complements the kinematic boundary condition 
$\delta\varphi(0)=0$.
Substitution of the general solution (\ref{eqgensol}) 
with $\delta Z_{ab}$ set to zero
into these two 
conditions then gives 
\bea 
C_1&=&0\\
-C_1 k\sin kL + C_2 k \cos kL &=& 0
\eea 
A non-trivial solution exists only if $\cos kL=0$. Since $k>0$,
this condition is satisfied for $kL=(n-1/2)\pi$ where $n$ is a positive integer. For $n=1$, condition $kL=\pi/2$ can be rewritten
in terms of force $P$ as
\beq \label{ee314}
(1+\beta P)P = P_E
\eeq 
where
\beq 
P_E = \frac{EI\pi^2}{4L^2}
\eeq
\\is the
Euler critical load for a column of buckling length $L_b=2L$.
The solution of our problem would be equal to the Euler
critical load only if $\beta=0$, which happens for $GA_s=EA$.
For nonzero values of $\beta$, it is necessary to solve a quadratic
equation, which has two roots,
\beq \label{eqroots}
P_{1,2}=\frac{-1\pm \sqrt{1+4\beta P_E}}{2\beta}
\eeq 
\item
For a {\bf simply supported column}, the perturbations $\delta Z_{ab}$ and $\delta M_{ab}$ vanish by virtue of (\ref{eqs219})--(\ref{eqs220}), and
so the general solution (\ref{eqgensol}) reduces again to $\delta\varphi(\xi)=C_1\cos k\xi + C_2\sin k\xi$. Boundary conditions
(\ref{eqs217})--(\ref{eqs218}) lead to analogous conditions
for $\delta z_s(0)$ and $\delta z_s(L)$, which can be substituted into (\ref{ee293}) to obtain conditions $\delta\varphi'(0)=0$ and $\delta\varphi'(L)=0$ written in terms of the primary unknown
function. Substituting the general solution, we get linear
equations for the integration constants,
\bea 
C_2 k &=& 0 \\
-C_1 k\sin kL + C_2 k \cos kL &=& 0
\eea
A non-trivial solutions exists only if $\sin kL=0$, which 
holds for $k=n\pi/L$ where $n$ is a positive integer. 
In terms of the loading force $P$, this leads again
to equation (\ref{ee314}), in which 
\beq 
P_E = \frac{EI\pi^2}{L^2}
\eeq
is the
Euler critical load for a column of buckling length $L_b=L$.
For $\beta\ne 0$, the roots are again given by (\ref{eqroots}),
just the meaning of $P_E$ is different as compared to the 
column clamped at one end.
\item
Finally, for a {\bf column clamped at both ends}, 
the values of $\delta Z_{ab}$ and $\delta M_{ab}$ are not 
deducible from global equilibrium (since this case is statically 
indeterminate), but we can use four kinematic boundary conditions
(\ref{eqs221})--(\ref{eqs224}), rewritten in terms of the perturbations. Conditions (\ref{eqs222}) and (\ref{eqs224})
applied to the general solution (\ref{eqgensol}) 
lead to 
\bea \label{eqZ}
-\frac{\delta Z_{ab}}{P} + C_1 &=& 0 \\
-\frac{\delta Z_{ab}}{P} + C_1 \cos kL + C_2 \sin kL&=& 0
\eea 
while conditions (\ref{eqs221}) and (\ref{eqs223}) need to be combined
with (\ref{ee293}) written at $\xi=0$ and $\xi=L$, in order to
construct conditions for the primary unknown function $\varphi$.
The resulting equations read
\bea \label{eqM}
 C_2 k &=& -\frac{\delta M_{ab}}{EI}\\
-C_1 k \sin kL + C_2 k \cos kL &=& -\frac{\delta M_{ab}}{EI}
-\frac{\delta Z_{ab}}{EI}\left(1-\frac{P}{EA}\right)L
\eea 
Elimination of the static unknowns $\delta Z_{ab}$ and $\delta M_{ab}$ based on (\ref{eqZ}) and (\ref{eqM}) then yields a set of two homogeneous linear equations
for the integration constants:
\bea 
\left(
\begin{array}{cc}
    \cos kL-1 & \sin kL  \\
   \left(1-\frac{P}{EA}\right)\frac{PL}{EI} - k\sin kL   & k(\cos kL-1)
\end{array}
\right)
\left(\begin{array}{c}
C_1 \\ C_2 \end{array}\right) =
\left(\begin{array}{c}
0 \\ 0 \end{array}\right)
\eea 
A non-trivial solution exists only if the determinant of the 
square matrix vanishes. After simple rearrangements,
we obtain a nonlinear equation
\beq \label{eqs245}
2k\sin^2\frac{kL}{2}-\left(1-\frac{P}{EA}\right) \frac{PL}{EI}\sin\frac{kL}{2}\cos\frac{kL}{2}=0
\eeq 
which is satisfied if
\beq \label{eqs246}
\sin\frac{kL}{2} = 0
\eeq 
or if 
\beq \label{eqs247}
\tan\frac{kL}{2}=\left(1-\frac{P}{EA}\right) \frac{PL}{2kEI}
\eeq 
The first case is easy to handle and leads to $k=2n\pi/L$,
which for $n=1$ yields once again equation 
 (\ref{ee314}), in which 
\beq 
P_E = \frac{4EI\pi^2}{L^2}
\eeq
is the
Euler critical load for a column of buckling length $L_b=L/2$.
Analysis of condition (\ref{eqs247}) is postponed to the 
detained discussion of compressive and tensile bifurcations,
because its relevance is affected by the sign of $P$.
\end{enumerate}

The analysis has confirmed that the compressive critical load that corresponds
to various types of support arrangements can be evaluated
from a quadratic equation in a unified form (\ref{ee314}),
in which the constant $P_E$ on the right-hand side has the meaning
of the classical Euler buckling load $P_E=EI\pi^2/L_b^2$,
where $L_b$ is the buckling length with the same meaning
as in the classical analysis, i.e., $L_b=2L$ for a column
clamped at one end, $L_b=L$ for a simply supported (pin-ended) 
column, and $L_b=L/2$ for a column clamped at both ends.
Equation (\ref{ee314}) is equivalent with
a condition derived, in a somewhat different notation,
by Haringx in his  paper on helical springs (see the unnumbered equation between
equations (11) and (12) in Haringx\cite{Haringx1942},
in which $\xi$ corresponds to our $P_0/EA$, $l_0$ to our $L$, and $\alpha_0$, $\beta_0$ and $\gamma_0$ to our $EI$, $GA_s$ and $EA$, respectively).

For standard beams, the shear sectional stiffness, $GA_s$, can be expected to be smaller than the 
axial sectional stiffness, $EA$, unless the material and section have very
unusual properties (e.g., a strongly auxetic material with $\nu<-0.5$ and thus $G>E$). For $GA_s<EA$, parameter $\beta$ 
given by (\ref{eq:beta}) is positive,
both roots are real, and it is interesting to see that one
is positive and the other negative, indicating a potential bifurcation in tension.
Haringx\cite{Haringx1942} focused exclusively on the root that corresponds to compression. 
It is also worth noting that he provided closed-form
expressions for the equivalent sectional stiffnesses of a helical spring modeled as a straight beam. If these expressions are adopted, the resulting dimensionless ratio $GA_s/EA$ is $2(1+m)/m$, where 
$m$ was referred to as the Poisson ratio by Haringx, but 
its value was substituted as 10/3 and in today's terminology
it represents the reciprocal value $1/\nu$ of the actual
Poisson ratio $\nu$. This leads, for springs, 
to $GA_s/EA=2(1+\nu)$, which is typically somewhere near 2.5. From this point of view, it makes sense to 
consider cases with $GA_s>EA$ as well.

\subsubsection{Compressive buckling}

For easier interpretation of the results, it is useful to work with
a dimensionless stiffness ratio $\Gamma=GA_s/EA$, which is usually
smaller than 1 but for springs can be larger. For a typical material and section,
we can set for instance $E/G=2(1+\nu)=2.5$ (isotropic material with Poisson ratio 0.25) and $A/A_s=1.2$ 
(rectangular section), which gives $\Gamma=1/3$, as used in several examples in Section~\ref{sec:examples}.
According to the analysis presented in Section~\ref{sec:reis-stab},
the critical load in compression is for the Reissner model
given by the positive root of equation (\ref{ee314}), in which $P_E$
corresponds to the standard Euler critical load for the given
type of supports.
The positive root is obtained if the sign before the square root
in formula (\ref{eqroots}) is taken as positive. Since parameter
$\beta$ is given by (\ref{eq:beta}), 
the resulting formula for the critical load in compression
can be written as
\beq\label{eq:harourx}
P_{Rc} =\dfrac{\sqrt{1+4 \left( \dfrac{1}{GA_s}-\dfrac{1}{EA} \right)P_{E}}-1}{2 \left( \dfrac{1}{GA_s}-\dfrac{1}{EA} \right)}= \frac{EA}{2}\frac{\Gamma}{1-\Gamma}\left(\sqrt{1+\frac{4(1-\Gamma)}{\Gamma}\frac{P_E}{EA}}-1\right)
\eeq
For slender columns, the Euler critical load $P_E$ is much smaller than the sectional stiffness $EA$.
It is then possible to use the approximation
\beq\label{eq:harour}
P_{Rc} \approx
 \frac{EA}{2}\frac{\Gamma}{1-\Gamma}\left(1+\frac{1}{2}\cdot\frac{4(1-\Gamma)}{\Gamma}\frac{P_E}{EA}-1\right)= P_{E} = \frac{EI\pi^2}{L^2}
\eeq
Without the approximation, the exactly evaluated critical load will be somewhat smaller
than the Euler load (provided that $\Gamma<1$, i.e., $GA_s<EA$). Let us explore how the slenderness affects the critical
stress, defined as $\sigma_{Rc}=P_{Rc}/A$ (actually representing the magnitude of the compressive stress at the onset of buckling). Introducing the radius of inertia
$i_{y}=\sqrt{I/A}$ and the slenderness ratio $s=L_b/i_y$, we can express
$P_E/EA=\pi^2/s^2$
and 
\beq 
\sigma_{Rc} = \frac{P_{Rc}}{A} = \frac{E}{2}\frac{\Gamma}{1-\Gamma}\left(\sqrt{1+\frac{4\pi^2(1-\Gamma)}{\Gamma s^2}}-1\right)
\label{282e}
\eeq 
For large slenderness, $s\to\infty$, this gives 
$\sigma_{Rc} \approx E\pi^2/s^2$, 
which corresponds to the Euler hyperbola.
For small slenderness, $s\to 0^+$, we get 
\beq \label{ee324}
\sigma_{Rc} \approx \dfrac{E\pi}{s}\sqrt{\dfrac{\Gamma}{1-\Gamma}} 
\eeq 
So the critical stress tends to infinity, in inverse proportion to the slenderness.
Of course, this is just theoretical description of the asymptotic behavior 
of the solution. In reality, the solution becomes physically
meaningless for $\lambda_{s0}\le 0$, i.e., for a critical stress that
attains or exceeds $E$. For the already mentioned typical material and section
characterized by $\Gamma=1/3$, formula (\ref{ee324}) would give $\sigma_{Rc}\approx E$ 
for slenderness $s=\pi/\sqrt{2}\approx 2.22$. The sectional equation that postulates proportionality
between the normal force and the engineering strain is certainly
not realistic if the strain becomes truly large.
Nevertheless, the analysis of asymptotic properties is useful
for comparison of the main features of various models.

\begin{figure}[h]
\centering
\includegraphics[scale=1]{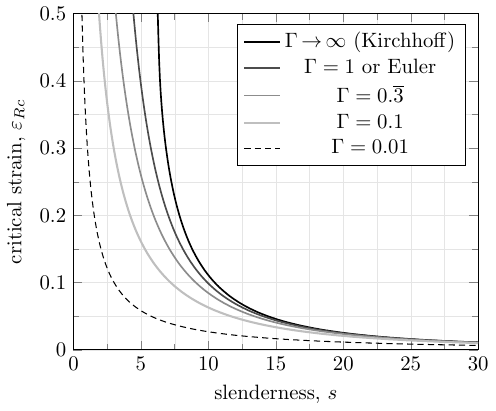}
    \caption{Reissner model: dependence of the critical compressive strain (i.e., normalized critical stress) $\eps_{Rc}=\sigma_{Rc}/E$ on the slenderness $s$ for selected
values of parameter $\Gamma=GA_s/EA$. The Euler hyperbola described by $\eps_{E}=\pi^2/s^2$ is obtained as a special case with $GA_s\to\infty$ and $EA\to\infty$, and it also corresponds to the case when $\Gamma=1$.}
    \label{5a}
\end{figure}

For illustration, the dependence of the normalized critical stress $\sigma_{Rc}/E$ on the slenderness $s$ is plotted in Fig. \ref{5a} for various
values of parameter $\Gamma$. The ratio $\sigma_{Rc}/E=\eps_{Rc}$ can be interpreted as the critical strain (in compression, but with a positive sign). For comparison, the graph also contains curves that are obtained if the shear compliance is neglected (Kirchhoff model), or if both the shear
and the normal compliance are neglected (Euler model). In these extreme cases,
it is not possible to use formula (\ref{eq:harourx}), which is valid only as long as $\Gamma<1$, i.e., $GA_s<EA$. For the limit cases, we need to
go back to condition (\ref{ee314}). If $GA_s\to\infty$ while $EA$ is finite, the equation remains quadratic and its solutions
are still given by (\ref{eqroots}), in which $\beta=-1/(EA)$ is negative. The denominator is then negative and both roots are positive,
the smallest one being given by
\beq 
P_{K} = \frac{EA}{2}\left(1-\sqrt{1-\frac{4P_E}{EA}}\right)
\eeq 
from which
\beq 
\sigma_{K} =\frac{P_K}{A} = \frac{E}{2}\left(1-\sqrt{1-\frac{4\pi^2}{s^2}}\right)
\eeq 
On the other hand, if both $GA_s$ and $EA$ are considered infinite,
condition (\ref{ee314}) reduces to a linear equation with only one solution,
$P=P_{E}$. The critical load is then the classical Euler buckling
load, $P_E=EI\pi^2/L_b^2$, and the corresponding graphical representation of the critical strain (or rather normalized critical stress) $\eps_{E}=\pi^2/s^2$ is  the
Euler hyperbola. Interestingly, exactly the same
dependence would be obtained for finite stiffnesses if $GA_s=EA$, i.e. $\Gamma=1$ and $\beta=0$.
For $\Gamma>1$, the corresponding curve is above the Euler hyperbola, which means that 
the properly evaluated critical load would be
larger than the Euler load. This is caused by
the stabilizing effect of the shortening of the
column due to compressive axial deformation.
The Euler solution refers to the buckling length evaluated
from the initial 
undeformed length, $L$, but the buckling load of 
an axially compressible column with zero
shear flexibility is $EI\pi^2/(\lambda_{crit}L_b)^2$ where
$\lambda_{crit}$ is the stretch at the onset of buckling.
Added shear flexibility has a destabilizing effect, but as long as $GA_s>EA$, the stabilizing effect of the axial shortening prevails. However, this is rather
an artificial case. As already mentioned, 
the shear stiffness is typically smaller than
the axial stiffness, and then the critical
load is below the Euler load. For low
shear stiffness, the difference can be dramatic.

\subsubsection{Tensile bifurcation}\label{sec:teninstab}

\begin{figure}[h!]
    \centering
  \includegraphics[scale=1]{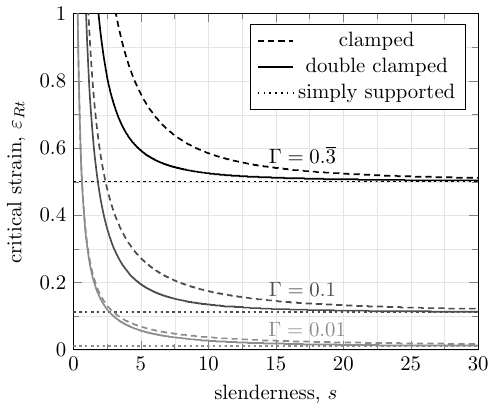}   
\caption{Reissner model: dependence of the critical tensile strain on the slenderness for three selected
values of parameter $\Gamma=GA_s/EA$, evaluated for a bar clamped at one end (dashed curves), 
 a bar clamped at both ends (solid curves), and
 a simply supported bar (dotted straight lines)} 
\label{6a}
\end{figure}

In addition to the positive root, equation (\ref{ee314}) has also a negative
root, which indicates a bifurcation in tension (in cases when $\Gamma<1$).
The critical load for tension is evaluated as
\beq\label{eq:285}
\vert P_{Rt}\vert = \dfrac{1+\sqrt{1+4P_{E} \left( \dfrac{1}{GA_s}-\dfrac{1}{EA} \right)}}{2 \left( \dfrac{1}{GA_s}-\dfrac{1}{EA} \right)} 
\eeq
and the corresponding critical strain in tension is 
\bea\label{eq258}
\eps_{Rt} = 
\frac{\vert P_{Rt}\vert}{EA} = \frac{\Gamma}{2(1-\Gamma)}\left(1+\sqrt{1+\dfrac{4(1-\Gamma)\pi^2}{\Gamma s^2}}\right)
\eea
For selected values of $\Gamma=1/3$ (homogeneous bar with $\nu=0.25$ and rectangular section), $\Gamma=0.1$ and 0.01 (sandwich bars with relatively low and extremely low shear
stiffness), the dependence of the critical strain on slenderness given by (\ref{eq258}) is visualized 
by the dashed curves in  Fig. \ref{6a} (labeled as 'clamped', which will be explained later).
Since the analyzed member is loaded by axial tension, we refer in this section to a bar instead of a column, but this is just a formal difference.
If the critical strain is too high,
the loss of stability would be predicted for a state in which the simple
linear model is not realistic, and refined analysis would be needed. 
For increasing slenderness, the critical strain approaches (from above) a finite limit, given by $\Gamma/(1-\Gamma)$. For our typical
case with $\Gamma=1/3$, the limit value would be 0.5, i.e., 
50\% strain. Physically, the solution is not totally meaningless,
but the range in which the relation between the nominal stress and 
engineering strain is approximately linear would most likely be exceeded. 
However, for sections with a low shear stiffness (compared to
the axial stiffness), the critical strain in tension could easily fall
into the almost linear range.
Therefore, the model is in principle capable of describing a certain type of bifurcation under axial tension, caused
by a sliding-type mechanism similar
to those observed in \cite{Zaccaria2011}.

Interestingly, it turns out that the tensile bifurcation cannot
be described in a unified format, covering various combinations
of supports by a single formula. 
If the loading leads to tension, which means that the force $P$
introduced in Fig.~\ref{fig:stab1} is considered as negative,
the analysis of the critical state presented in Section~\ref{sec:reis-stab}
needs to be refined. The first assumption that needs attention is that the general solution (\ref{eqgensol}) of differential equation (\ref{eqr1}) was written only for the
case when $(1+\beta P)P>0$. This is guaranteed if compression
is applied ($P>0$) and parameter $\beta$ is positive ($GA_s<EA$).
However, in general, one should also investigate cases
when $(1+\beta P)P$ is zero or negative. 

For $(1+\beta P)P<0$, functions sine and cosine in the general
solution are replaced by their hyperbolic counterparts. Detailed
analysis of various combinations of boundary conditions leads
to the conclusion that the only solution is the trivial one,
which means that the bifurcation is excluded. On the other hand,
the assumption that  $(1+\beta P)P=0$, which effectively means
that $1+\beta P=0$, reduces equation (\ref{eqr1}) to $EI\,\delta\varphi''=0$
and the general solution is an arbitrary linear function
$\delta\varphi(\xi)=C_1+C_2\xi$. For bars clamped at one or at both ends,
the boundary conditions do not allow for any nontrivial solution
of this kind, but for a {\bf simply supported bar} (pinned at both ends),
a constant function $\delta\varphi(\xi)=C_1$ satisfies conditions
$\delta\varphi'(0)=0$ and well as $\delta\varphi'(L)=0$, and a non-trivial
solution exists. A constant rotation function combined with
zero displacement functions is a possible bifurcation mode.
The corresponding critical force is determined from 
condition $1+\beta P=0$, which gives 
\beq 
P_{crit}=-\frac{1}{\beta} = \frac{GA_s\cdot EA}{GA_s-EA} = \frac{\Gamma}{\Gamma-1}EA
\eeq 
For $\Gamma>1$ (i.e., $GA_s>EA$), this force is positive and larger 
than $EA$, which is not admissible because the column length in such
a critical state would be negative. However, for $\Gamma<1$,
this critical force is negative and corresponds to a bar
under tension, with the critical strain equal to
$\Gamma/(1-\Gamma)$, which is always lower than
the critical strain from (\ref{eq258}). Consequently, the 
simply supported bar under tension 
(described by the Reissner model with $GA_s<EA$) can bifurcate into a special uniform mode with constant rotation
and zero displacements, and this happens at critical strain
\beq \label{eq:epsrtss}
\eps_{Rt}^{(ss)} = \frac{\Gamma}{1-\Gamma} = \frac{GA_s}{EA-GA_s}
\eeq 
Superscript 'ss' refers to 'simply supported'.
It is remarkable that this critical strain strain is independent of the slenderness,
and it is equal to 1/2 for $\Gamma=1/3$, to 1/9 for $\Gamma=0.1$,
and to 1/99 for $\Gamma=0.01$; see the dotted horizontal lines 
in Fig.~\ref{6a}.

Yet another refinement of the bifurcation analysis under tension
is needed for a {\bf bar clamped at both ends}. The reason is that 
condition (\ref{eqs245}) derived in Section~\ref{sec:reis-stab} can be satisfied in two ways,
described by (\ref{eqs246}) and (\ref{eqs247}), respectively. We have tacitly assumed that
(\ref{eqs246}) is satisfied for a lower load level than (\ref{eqs247}) and thus determines
the critical load. This is true for compression, when $P>0$
and the right-hand side of (\ref{eqs247}) is positive. The smallest positive
root of (\ref{eqs246}) is $kL/2=\pi$, while (\ref{eqs247}) with a positive right-hand side
cannot be satisfied by any positive value of $kL/2$ between 0 and $\pi$. 
This is why it was correct to describe the critical load in compression
by equation (\ref{ee314}), with $P_E$ set to the standard Euler critical load.
However, when we focus on tensile bifurcations,  
the load $P$ is assumed to be negative, which leads to a negative
right-hand side of (\ref{eqs247}), and the condition is satisfied by
a certain value of $kL/2$ smaller than $\pi$. The corresponding critical
load in tension needs to be determined numerically and it is 
different from the negative root of quadratic equation (\ref{ee314}).  

It is important to realize that $k$ and $P$ are not two independent
variables, since they are linked by (\ref{eq:k}). Therefore, if $P$ is considered
as the primary unknown, equation (\ref{eqs247}) should be written in its fully explicit form  as
\beq \label{eqs247x}
\tan\left(\frac{L}{2}\sqrt{\frac{(1+\beta P)P}{EI}}\right)=\left(1-\frac{P}{EA}\right) \frac{PL}{2EI}\sqrt{\frac{EI}{(1+\beta P)P}}
\eeq
We are looking for the negative root $P$ with the smallest
absolute value. To facilitate the analysis, let us convert
the problem into the dimensionless format by substituting
$P=-EA\eps$, $EI=EA L^2/(4s^2)$ and $\beta=(1-\Gamma)/(\Gamma EA)$ 
where $\eps>0$ is the tensile strain, playing the role of the 
primary unknown,  $s$ is the slenderness derived from the buckling length $L_b=L/2$, and $\Gamma=GA_s/(EA)$ is the sectional stiffness ratio.
The resulting dimensionless equation reads
\beq \label{eqs247y}
F(\eps)\equiv \sqrt{\frac{1-\Gamma}{\Gamma}-\frac{1}{\eps}}\,\tan\left(s\,\sqrt{\left(\frac{1-\Gamma}{\Gamma}\eps-1\right)\eps}\right) + s\,(1+\eps)=0
\eeq 
The expressions under the square root must 
be non-negative,  which is possible only if $\Gamma <1$ and $\eps>\Gamma/(1-\Gamma)\equiv \eps_{Rt}^{(ss)}$
where $\eps_{Rt}^{(ss)}$ is the strain defined in (\ref{eq:epsrtss}).
In terms of the sectional stiffness, the
first constraint is equivalent to
$GA_s<EA$.

For $s>0$ and $\Gamma<1$, function $F$ on the left-hand side of (\ref{eqs247y}) is increasing at all points $\eps\in[\Gamma/(1-\Gamma),\infty)$ at which it is defined (remains finite).
It is not defined for values of $\eps$ for which the argument
of the tangent function is equal to an odd multiple of $\pi/2$. 
The first of these special points is
\beq 
\eps_1 = \frac{\Gamma}{2(1-\Gamma)}\left(\sqrt{1+\frac{(1-\Gamma)\pi^2}{\Gamma s^2}}+1\right)
\eeq 
It is easy to show that $F(\Gamma/(1-\Gamma))>0$ (because the first term in (\ref{eqs247y}) vanishes and the second term is positive), and since $F$ is increasing on the interval
$[\Gamma/(1-\Gamma),\eps_1)$, equation (\ref{eqs247y}) has no solution in this interval. As $\eps_1$ is crossed, the value of 
$F$ jumps ``from plus infinity to minus infinity'', and then 
it increases continuously. At $\eps=\eps_{Rt}$, the argument
of the tangent function is $\pi$, and so the first term vanishes 
while the second term is positive. Consequently, $F(\eps_{Rt})>0$,
and the interval $(\Gamma/(1-\Gamma),\eps_{Rt})$ must contain 
exactly one solution of equation (\ref{eqs247y}).
Since this solution is smaller than $\eps_{Rt}$, the bifurcation
occurs at a lower strain level (and lower load level) than 
predicted by formula (\ref{eq258}), which was derived from condition
(\ref{eqs246}) instead of (\ref{eqs247}).
The dependence of the critical strain $\eps_{Rt}^{(dc)}$ 
that represents the smallest positive solution of equation
(\ref{eqs247y}) on the slenderness is indicated in Fig.~\ref{6a}
by the solid curves. 
They were constructed numerically
by solving equation (\ref{eqs247y}) on the interval $(\Gamma/(1-\Gamma),\eps_{Rt}(s))$ using the Newton method. Superscript 'dc' refers
to a 'double clamped' bar. For a given slenderness, the critical strain
evaluated for a bar clamped at both ends is always between the
critical strain $\eps_{Rt}^{(ss)}$ for a simply supported bar
and the critical strain $\eps_{Rt}$ evaluated from formula (\ref{eq258}),
which turns out to be valid only for a {\bf bar clamped at one end}.
Let us recall that the slenderness is understood as the ratio $s=L_b/i_y$
where $i_y=\sqrt{I/A}$ is the sectional radius of inertia
and $L_b$ is the buckling length, which is still taken
as $L_b=2L$ for a bar clamped at one end and $L_b=L/2$ for a bar
clamped at both ends, even though these standard expressions
lose their physical meaning---they do not represent here the distance
between the inflexion points of the bifurcated centerline. 

The foregoing observation motivates a deeper investigation of the ``buckling mode'' that corresponds to the tensile bifurcation. 
\begin{enumerate}
    \item For a {\bf bar clamped at one end}, the buckling mode is described
    by functions
\bea 
\delta\varphi(\xi) &=& \widehat{\delta\varphi}\, \sin\frac{\pi\xi}{2L} \\
\delta z_s(\xi) &=& \left(\sqrt{1+\frac{(1-\Gamma)\pi^2}{\Gamma s^2}}-1\right)\frac{L\,\widehat{\delta\varphi}}{\pi}\left(1-\cos\frac{\pi\xi}{2L}\right)
\eea 
The rotation function is obtained from the general solution (\ref{eqgensol}) by substituting 
$\delta Z_{ab}=0$ and $C_1=0$ and by formally replacing the arbitrary
constant $C_2$ by $\widehat{\delta\varphi}$, interpreted as the maximum
rotation found at the free end. The lateral displacement function is then
constructed by integrating the right-hand side of (\ref{ee292}) and imposing boundary
condition $\delta z_s(0)=0$.
\\
Recall that the tensile bifurcation occurs only if $GA_s<EA$, i.e., $\Gamma<1$.
Therefore, if constant $\widehat{\delta\varphi}$ is positive and all sections rotate counterclockwise, then the
deflection is also positive along the bar and the centerline moves downward. 
For an Euler or Kirchhoff beam with the same deformed shape of the centerline, the sign of the rotation would be negative.
This means that the sections rotate
in the opposite direction than the centerline, which amplifies the shear strains. 
\\
As an example, consider the stiffness
ratio $\Gamma=1/3$ and slenderness $s=20$. 
The corresponding critical strain is 0.5236, and
the ratio between the magnitudes of
the shear angle $\delta\chi=\delta\varphi+\delta z_s'/(1+\eps_{crit}^{ten})$
and of the angle $\delta\varphi_{Euler}=-\delta z_s'/(1+\eps_{crit}^{ten})$ by which the
centerline deviates from the undeformed axis is 33.3. 
Therefore, after the bifurcation,
the centerline remains ``almost straight'' and
the rotation of the sections is mainly due to shear distortion,
as will be illustrated by a numerically simulated post-critical
deformed shape in Fig.~\ref{f:ex6_4shapes}a.
\item For a {\bf simply supported bar}, the buckling mode is 
characterized by a constant rotation function and zero displacement function,
which can be formally described by
\bea 
\delta\varphi(\xi) &=& \widehat{\delta\varphi}  \\
\delta z_s(\xi) &=& 0
\eea 
This special case was mentioned in the discussion leading to formula (\ref{eq:epsrtss}). Here, the centerline remains straight while all sections
rotate by the same angle; see  Fig.~\ref{f:ex6_4shapes}b.
\item For a {\bf bar clamped at both ends}, the buckling mode is described
    by functions
\bea 
\delta\varphi(\tilde\xi) &=&  \frac{\widehat{\delta\varphi}}{1-\cos\kappa}\left(\cos\kappa\tilde\xi-\cos\kappa\right) \\
\delta z_s(\tilde\xi) &=& \frac{L\,\widehat{\delta\varphi}}{1-\cos\kappa}
\left(\frac{2\kappa}{ s^2\,\eps_{Rt}^{(dc)}}\left(\sin\kappa+\sin\kappa\tilde\xi \right)+ \left(1+\eps_{Rt}^{(dc)}\right)\frac{1+\tilde\xi}{2}\cos\kappa\right)
\nonumber \\
\eea 
where
\beq 
\tilde\xi = \frac{2\xi}{L}-1
\eeq 
is a normalized coordinate that varies from $-1$ at the left end to 1 at the right end,
\beq 
\kappa = \frac{s}{2}\sqrt{\left(\frac{1-\Gamma}{\Gamma}\eps_{Rt}^{(dc)}-1\right)\eps_{Rt}^{(dc)}}
\eeq 
is an auxiliary variable introduced for convenience, and $\eps_{Rt}^{(dc)}$
is the smallest positive solution of equation (\ref{eqs247y}).
Constant $\widehat{\delta\varphi}$ represents the rotation in the middle of the bar (at $\tilde\xi=0$). An example of a deformed shape in a post-critical
state will be shown in Fig.~\ref{f:ex6_4shapes}c.
\end{enumerate}

\subsection{Critical load -- Ziegler model}
\label{sec:5.1.3}

The fundamental solution (\ref{ee286})--(\ref{ee288}) satisfies not only the governing equations
of the Reissner model, but also equations 
(\ref{ee110})--(\ref{ee111}) and (\ref{ee137}), which describe
the Ziegler model. However, if the state is perturbed,
the linear equations that link the infinitesimal perturbations
are different. For the Ziegler model, we obtain
\bea 
\delta x_s'&=& 0 \\
 \delta z_s' &=& \left(1-\dfrac{P}{EA}\right)\frac{\left(1-\dfrac{P}{EA}\right)\delta Z_{ab}+GA_s\delta\varphi }{\left(1-\dfrac{P}{EA}\right)P-GA_s}\\
 \delta\varphi' &=&  \frac{1}{EI}\left(P\,\delta z_s - \delta M_{ab}-\delta Z_{ab}\left(1-\dfrac{P}{EA}\right)\xi\right)
\eea
Similar to the procedure used for the Reissner model, the second and third equations can be combined into a second-order
differential equation
\beq 
EI \left(\left(1-\frac{P}{EA}\right)P-GA_s\right) \,\delta\varphi''
-\left(1-\frac{P}{EA}\right)PGA_s\,\delta\varphi = \left(1-\frac{P}{EA}\right)GA_s\,\delta Z_{ab}
\eeq 
with a single unknown function, $\delta \varphi$. The general solution 
has the form 
\beq 
\delta\varphi(\xi) = -\frac{\delta Z_{ab}}{P} + C_1\cos k\xi + C_2 \sin k\xi
\eeq 
where
\beq \label{eq:kk}
k = \sqrt{\frac{P}{EI}\cdot\frac{1-\dfrac{P}{EA}}
{1-\left(1-\dfrac{P}{EA}\right)\dfrac{P}{GA_s}}}
\eeq
Of course, this is based on the assumption that the expression under the square
root is positive. 

Integration constants $C_1$ and $C_2$ and also unknowns $\delta Z_{ab}$
and $\delta M_{ab}$ have to be determined from the boundary conditions.
It is possible to treat various combinations of boundary conditions
in a unified way, simply by introducing the appropriate buckling length.
\begin{enumerate}
\item
For a column clamped at one end and free at the other, we set
$\delta Z_{ab}=0$, and the integration constants
are determined from conditions $\delta\varphi(0)=0$ and $\delta\varphi'(L)=0$. A nontrivial solution exists if $\cos kL=0$,
and the smallest positive value of $k$ satisfying this condition is
$\pi/(2L)$, which can be written as $\pi/L_b$ where $L_b=2L$ is the
buckling length.
\item 
For a simply supported column,  we set again
$\delta Z_{ab}=0$ and $\delta M_{ab}=0$, and the integration constants
are determined from conditions $\delta\varphi'(0)=0$ and $\delta\varphi'(L)=0$. A nontrivial solution exists if $\sin kL=0$,
and the smallest positive value of $k$ satisfying this condition is
$\pi/L$, which can be written as $\pi/L_b$ where $L_b=L$ is the
buckling length.
\item
For a column clamped at both ends, we have to treat
$\delta Z_{ab}=0$ and $\delta M_{ab}=0$ as unknowns and determine
them simultaneously with the integration constants
from conditions $\delta\varphi(0)=0$, $\delta z_s(0)=0$, $\delta\varphi(L)=0$ and $\delta z_s(L)=0$. A nontrivial solution exists if $\cos kL=1$,
and the smallest positive value of $k$ satisfying this condition is
$2\pi/L$, which can be written as $\pi/L_b$ where $L_b=L/2$ is the
buckling length.
\end{enumerate}
The bifurcation condition can thus be written in a unified way as $k=\pi/L_b$, and when $k$ is replaced by the expression in (\ref{eq:kk}), the condition can be
rewritten as
\beq\label{ee339} 
\left(\dfrac{1}{GA_s}+\frac{1}{P_{E}}\right)
\left(1-\dfrac{P}{EA}\right)P
 = 1
\eeq 
where $P_E=EI\pi^2/L_b^2$ is the Euler buckling load evaluated for the
given buckling length. In terms of $P$, condition (\ref{ee339}) is a quadratic equation with  roots
\beq \label{eee271}
P = \frac{EA}{2}\left(1\pm \sqrt{1-\frac{4}{EA}\frac{P_{E}GA_s}{P_{E}+GA_s}}\right)
\eeq 
If the discriminant (expression under the square root) is non-negative, both roots are positive.
The lower value is obtained with the negative sign in front of the
square root.
The critical load for the Ziegler model under compression is thus given by
\beq\label{eq:engour}
P_{Zc} =  \frac{EA}{2}\left(1- \sqrt{1-\frac{4}{EA}\frac{P_EGA_s}{P_E+GA_s}}\right) = \frac{EA}{2}\left(1- \sqrt{1-\frac{4\pi^2\Gamma}{\Gamma s^2+\pi^2}}\right)
\eeq 
which is clearly an expression different from formula (\ref{eq:harourx}) that was derived
for the Reissner model. Another difference is that, for the Ziegler model,
bifurcation in tension is not possible, because we have not obtained any negative root. 

In some cases the  discriminant of equation (\ref{ee339}) can be negative and 
the roots are then complex, which means that the fundamental solution
never bifurcates (even in compression). The condition of a nonnegative discriminant can be converted into
\beq 
\left(4GA_s-EA\right)P_E \le EAGA_s
\eeq 
For $\Gamma\equiv GA_s/(EA)\le 1/4$, this condition is always satisfied and the critical
load can indeed be evaluated from (\ref{eq:engour}). For $\Gamma\equiv GA_s/(EA)> 1/4$, real roots are obtained
only if the Euler load satisfies the condition
\beq 
P_E \le \frac{EAGA_s}{4GA_s-EA} =
EA\,\frac{\Gamma}{4\Gamma-1}
\eeq 
which can be expressed in terms of slenderness as
\beq\label{eq:limslen}
s \ge \pi \sqrt{4-\dfrac{1}{\Gamma}} 
\eeq 
It is interesting to explore what would happen if the condition is satisfied as an equality. For the typical case with $\Gamma=1/3$, this corresponds to slenderness
$s=\pi$ and Euler load $P_E=EA$. The corresponding Ziegler critical load is 
$P_{Zc}=EA/2$. If the slenderness is reduced just a little bit,
the critical load jumps to infinity. This can be interpreted as a consequence of the fact that by compressing the column,
we make it shorter and thus reduce the buckling length, which has a stabilizing effect. On the other hand, the increasing compressive stress
has a destabilizing effect. For very stocky columns, the stabilizing effect may 
dominate. Let us emphasize once again that analysis of such extreme cases
is meant to elucidate the asymptotic properties of the mathematical model,
while from the physical point of view it is not realistic to assume
that the linear stress-strain law would still be applicable.

\begin{figure}[h!]
  \hskip -20mm
    \begin{tabular}{cc}
    (a) & (b) 
    \\
\includegraphics[scale=1]{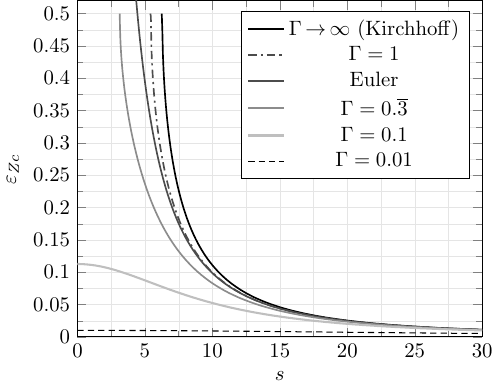}
    &  \includegraphics[scale=1]{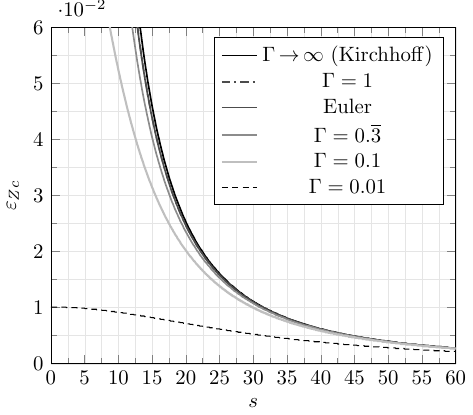}
    \end{tabular}
\caption{Ziegler model: (a) dependence of the compressive critical strain $\eps_{Zc}=P_{Zc}/EA$ on the slenderness $s$ for selected
values of parameter $\Gamma=GA_s/EA$, (b) close-up view for small critical strain values.} 
    \label{7a}
\end{figure}

In summary, a bifurcation in compression is possible when $\Gamma\le 1/4$ (i.e., $4GA_s\le EA$) or the slenderness satisfies condition (\ref{eq:limslen}).
The corresponding critical strain is given by
\beq \label{eq:epsZc}
\eps_{Zc} = \frac{P_{Zc}}{EA} =  \frac{1}{2}\left(1- \sqrt{1-\frac{4\pi^2\Gamma}{\Gamma s^2+\pi^2}}\right)
\eeq
Its  dependence on the slenderness is illustrated by the graphs in Fig.~\ref{7a}, where part (a) provides the global
view and part (b) reveals more details in the range of
low critical strains.  For $\Gamma=1$, the critical strain
is given by
\beq 
\eps_{Zc} = \frac{1}{2}\left(1- \sqrt{\frac{s^2-3\pi^2}{s^2+\pi^2}}\right), \hskip 10mm \mbox{for }s\ge \sqrt{3}\,\pi
\eeq 
and in the limit of $\Gamma\to\infty$, it is given by
\beq 
\eps_{Zc} = \frac{1}{2}\left(1- \sqrt{1-\frac{4\pi^2}{s^2}}\right), \hskip 10mm \mbox{for }s\ge 2\pi
\eeq
For comparison, the Euler hyperbola is plotted in black.

For large slenderness $s$ and general stiffness ratio $\Gamma$, one can use the approximation
\beq 
\eps_{Zc} \approx \frac{\pi^2\Gamma}{\pi^2+\Gamma s^2}
\eeq
which shows that the Euler hyperbola is approached (from below) as $s\to\infty$.
In the opposite extreme, when the slenderness tends to zero, the limit behavior of the curve
depends on the stiffness ratio $\Gamma$.
For $\Gamma\le 1/4$,
the critical strain approaches a finite limit
\beq 
\eps_{Zc}\vert_{s=0} =  \frac{1}{2}\left(1- \sqrt{1-4\Gamma}\right)
\eeq 
which remains below 50 \%.
For $\Gamma> 1/4$, the critical strain
attains its maximum value of 50 \% when the slenderness becomes
equal to $\pi\sqrt{4-1/\Gamma}$, and there is no bifurcation (and thus no critical strain)
for still smaller slenderness values. 

The effect of shear on stability was in a simplified form taken into account by Engesser\cite{Engesser1891},
who proceeded in analogy to
the Euler solution and preserved the assumption of axial incompressibility. Formula (9) from Engesser\cite{Engesser1891}  would
be in the present notation written as
\beq \label{eq:eng}
P_{Zc} = \frac{P_E}{1+\dfrac{P_E}{GA_s}}
\eeq 
This agrees with the solution of the reduced version of 
equation (\ref{ee339}), in which $EA$ is considered as infinite
and thus the terms with $EA$ in the denominator are deleted. The resulting equation is linear and has a unique solution (\ref{eq:eng}).

An extension of Engesser's analysis to the case of
axially compressible columns was developed by
Ziegler\cite{Zie82}, who also critically compared both 
approaches to the stability of shear-flexible columns
(Haringx versus Engesser/Ziegler)
and identified the definition of internal forces
as the reason why the results differ.
Formula (5.14) from Ziegler\cite{Zie82} agrees
with our (\ref{eee271}).
Ziegler restricted attention to the initial deviation
from a straight shape and used certain approximations.
Nevertheless, his sectional equations (4.5)
exactly correspond to our (\ref{ee17x}), (\ref{eqQstar}) and (\ref{ee19x}), 
and they can form a basis
of a geometrically exact beam theory, which is why
we refer to the formulation presented in Section~\ref{sec:modmod}
as the Ziegler model. 
The link between the choice of sectional
equations for beam models and 
large-strain continuum elasticity was also 
discussed by Ba\v{z}ant\cite{Bazant1971ACS}
and further developed
in the Stability of Structures textbook by Ba\v{z}ant and Cedolin\cite{Bazant91Cedolin}.

\bibliographystyle{elsarticle-num} 
\bibliography{biblio_shearbeam}

\end{document}